\documentclass[reqno, a4paper, 10pt,oneside]{amsart}

\makeatletter
\def\specialsection{\@startsection{section}{1}%
	\z@{\linespacing\@plus\linespacing}{.5\linespacing}%
	{\normalfont}}
\def\section{\@startsection{section}{1}%
	\z@{.7\linespacing\@plus\linespacing}{.5\linespacing}%
	{\normalfont\scshape\bfseries}}

\makeatother

\usepackage[top=1in, bottom=1.4in, inner=1in, outer=1in, includehead, showframe=false]{geometry}
\usepackage[algo2e, linesnumbered, noline, noend, ruled]{algorithm2e}
\usepackage[foot]{amsaddr}
\usepackage{amsfonts}
\usepackage{amsmath}
\usepackage{amssymb}
\usepackage{amsthm}
\usepackage[english]{babel}
\usepackage{booktabs}
\usepackage{caption}
\usepackage{braket}
\usepackage{float}
\usepackage[T1]{fontenc}
\usepackage{graphicx}
\usepackage[utf8]{inputenc}
\usepackage{lmodern}
\usepackage{lipsum}
\usepackage{mathtools}
\usepackage{nicefrac}
\usepackage[defaultlines=2,all]{nowidow}
\usepackage[per-mode=symbol, group-digits=integer]{siunitx}
\usepackage[table,dvipsnames]{xcolor}
\usepackage{enumitem}
\usepackage{subcaption}
\usepackage{stmaryrd}
\usepackage{textcomp}

\usepackage{cite}
\usepackage{tikz}
\usepackage[allcolors=cyan,colorlinks]{hyperref}

\usetikzlibrary{calc, patterns}

\definecolor{c0}{RGB}{0, 107, 164}
\definecolor{c1}{RGB}{255, 128, 14}
\definecolor{c2}{RGB}{171, 171, 171}
\definecolor{c3}{RGB}{89, 89, 89}
\definecolor{c4}{RGB}{200, 82, 0}
\definecolor{c5}{RGB}{44,160,44}
\definecolor{c6}{RGB}{148,103,189}

	\theoremstyle{plain}
	\newtheorem{Theorem}{Theorem}[section]

	\theoremstyle{definition}
	\newtheorem{Definition}[Theorem]{Definition}
	\newtheorem{Remark}[Theorem]{Remark}
	
	\newtheorem*{Example*}{Example}
	\newtheorem*{Remark*}{Remark}
	
	\setenumerate{label=(\arabic*), ref=(\arabic*)}
	
	\newcommand{\qe}[1]{``#1''}

	\newcommand{\R}{\mathbb{R}}
	
	\newcommand{\transposed}{^\top}
	
	\newcommand{\norm}[2]{\left\lvert \left\lvert #1 \right\rvert \right\rvert_{#2}}
	\newcommand{\abs}[1]{\left\lvert #1 \right\rvert}
	\newcommand{\dmeas}[1]{\ \mathrm{d}#1}
	
	\newcommand{\grad}{\nabla}
	
	\newcommand{\integral}[1]{\int_{#1}}
	\newcommand{\dual}[3]{\left\langle #1 , #2 \right\rangle_{#3}}

	\newcommand{\admissiblegeom}{\mathcal{A}}

	\newcommand{\holdall}{\mathrm{D}}
	\newcommand{\vectorfield}{\mathcal{V}}
	\newcommand{\flow}{\Phi}

	\makeatletter
	\DeclareOldFontCommand{\rm}{\normalfont\rmfamily}{\mathrm}
	\DeclareOldFontCommand{\sf}{\normalfont\sffamily}{\mathsf}
	\DeclareOldFontCommand{\tt}{\normalfont\ttfamily}{\mathtt}
	\DeclareOldFontCommand{\bf}{\normalfont\bfseries}{\mathbf}
	\DeclareOldFontCommand{\it}{\normalfont\itshape}{\mathit}
	\DeclareOldFontCommand{\sl}{\normalfont\slshape}{\@nomath\sl}
	\DeclareOldFontCommand{\sc}{\normalfont\scshape}{\@nomath\sc}
	\makeatother

	\newcommand{\gradientdefo}{\mathcal{G}}

	\newcommand{\iidx}[1]{_{#1}}

	\numberwithin{equation}{section}
	
	\newcommand*\closure[1]{%
		\hbox{%
			\vbox{%
				\hrule height 0.5pt 
				\kern0.5ex
				\hbox{%
					\ensuremath{#1}%
				}%
			}%
		}%
	} 
	
\SetCommentSty{mycommfont}
\SetKwComment{Comment}{\# }{}

\begin{document}

\title[Enforcing Mesh Quality in Shape Optimization]{Enforcing Mesh Quality Constraints in Shape Optimization with a Gradient Projection Method}
\author{Sebastian Blauth$^{*,1}$}
\address{$^*$ Corresponding Author}
\address{$^1$ Fraunhofer Institute for Industrial Mathematics ITWM, Kaiserslautern, Germany}
\email{\href{mailto:sebastian.blauth@itwm.fraunhofer.de}{sebastian.blauth@itwm.fraunhofer.de}}
\author{Christian Leith\"auser$^1$}
\email{\href{mailto:christian.leithaeuser@itwm.fraunhofer.de}{christian.leithaeuser@itwm.fraunhofer.de}}

\begin{abstract}
	For the numerical solution of shape optimization problems, particularly those constrained by partial differential equations (PDEs), the quality of the underlying mesh is of utmost importance. Particularly when investigating complex geometries, the mesh quality tends to deteriorate over the course of a shape optimization so that either the optimization comes to a halt or an expensive remeshing operation must be performed before the optimization can be continued. 
	In this paper, we present a novel, semi-discrete approach for enforcing a minimum mesh quality in shape optimization. Our approach is based on Rosen's gradient projection method, which incorporates mesh quality constraints into the shape optimization problem. The proposed constraints bound the angles of triangular and solid angles of tetrahedral mesh cells and, thus, also bound the quality of these mesh cells. The method treats these constraints by projecting the search direction to the linear subspace of the currently active constraints. Additionally, only slight modifications to the usual line search procedure are required to ensure the feasibility of the method. We present our method for two- and three-dimensional simplicial meshes. We investigate the proposed approach numerically for the drag minimization of an obstacle in a two-dimensional Stokes flow, the optimization of the flow in a pipe governed by the Navier-Stokes equations, and for the large-scale, three-dimensional optimization of a structured packing used in a distillation column. Our results show that the proposed method is indeed capable of guaranteeing a minimum mesh quality for both academic examples and challenging industrial applications. Particularly, our approach allows the shape optimization of complex structures while ensuring that the mesh quality does not deteriorate.
	
	
\end{abstract}

{\noindent\footnotesize This is a post-peer-review, pre-copyedit version of an article published in Computer Methods in Applied Mechanics and Engineering. The final version is available online at \url{https://doi.org/10.1016/j.cma.2025.118451}.
}

\maketitle


\vspace{-0.5cm}

\section{Introduction}
\label{sec:introduction}

Shape optimization constrained by partial differential equations (PDEs) is concerned with the optimization of some objective functional depending on both the solution of a PDE as well as the shape of the domain on which the former is posed. For analyzing the sensitivity of such a functional w.r.t.\ infinitesimal changes of the shape, shape calculus (see, e.g., \cite{Delfour2011Shapes}) is used. This gives rise to the so-called shape derivative which enables the use of gradient-based optimization methods and, thus, facilitates the efficient optimization of relevant practical problems. The field of shape optimization has many interesting practical applications, e.g., the optimization of electrolysis cells \cite{Blauth2024Multi}, microchannel systems \cite{Blauth2021Model}, aircraft design \cite{Schmidt2013Three}, electromagnetic devices \cite{Gangl2022Multi}, or electrical impedance tomography \cite{Hintermueller2008Electrical}. 

In this paper, we consider free-form shape optimization based on shape calculus, which does not involve any low-dimensional parametrizations of the shape and does, thus, not limit the reachable shapes. During the shape optimization, the nodes of the mesh are moved to change the shape of the discretized domain. This deformation of the mesh often leads to a deterioration of the mesh quality, at least for some cells. Moreover, this can even cause a breakdown of the optimization algorithm, as discussed below, unless an expensive re-meshing is performed, which might not even always be possible. For these reasons, the task of computing mesh deformations which preserve the mesh quality and do not lead to too much degradation of the mesh quality has received great interest in the recent literature: In \cite{Schulz2016Efficient}, Steklov-Poincar\'e-type metrics for shape optimization have been introduced, which provide a mathematical justification for using the linear elasticity equations for computing mesh deformations. In \cite{Iglesias2018Two} and \cite{Etling2020First} nearly conformal and restricted mesh deformations, respectively, for shape optimization are investigated and analyzed numerically. In \cite{Herzog2023manifold, Herzog2024discretize} a novel discrete Riemannian metric for triangular meshes is introduced and applied to solve two-dimensional shape optimization with a regularization for the mesh quality. In \cite{Luft2021Simultaneous,Luft2021Pre} so-called pre-shape calculus is proposed and used to simultaneously optimize shapes and mesh quality. Finally, in \cite{Deckelnick2022novel} and \cite{Mueller2021novel} the use of $W^{1,\infty}$ and $p$-harmonic mesh deformations, respectively, is proposed to ensure a better mesh quality in the context of shape optimization. Other recent developments in the field of shape optimization consider the development and improvement of solution algorithms: In \cite{Blauth2021Nonlinear} and \cite{Schulz2016Efficient} nonlinear conjugate gradient methods and (limited memory) BFGS methods, respectively, are investigated for shape optimization. In \cite{Schmidt2023linear}, Newton methods for shape optimization are used. Finally, in \cite{Blauth2023Space}, the space mapping technique is extended to the field of shape optimization. 

For the numerical solution of shape optimization problems, the quality of the underlying mesh is crucial as this determines how accurately the involved PDEs can be solved. In particular, only a few bad mesh cells can lead to a degradation of the solution or can even make the numerical solution of the PDEs impossible. Even if the solution quality is not affected by the degradation of some mesh cells, this can still cause problems for the shape optimization: Often, the quality of the mesh cells is deteriorating as they tend to self-intersect or invert due to the deformations of the computational mesh. In this case, the line search procedure employed for gradient-based shape optimization algorithms needs to find a suitable step size which does not lead to inverted mesh cells. Consequently, this leads to a reduction of the step size and can finally lead to a failure of the entire optimization. We note that the problem of computing high quality mesh deformations does also arise, e.g., in the context of so-called arbitrary Lagrangian-Eulerian (ALE) methods for solid mechanics \cite{Gadala2004Recent, Donea2004Arbitrary, Boman2012Efficient} and fluid structure interaction \cite{Stein2003Mesh, Zuijlen2007Higher}. Moreover, we refer the reader to \cite{Sieger2015shape, Selim2016Mesh, Upadhyay2021Numerical} for an overview of mesh deformation approaches.

In this paper, we present a novel approach for preserving the mesh quality during shape optimization. We use a semi-discrete approach which enforces constraints on the quality of the mesh cells and, thus, guarantees that the mesh quality is sufficiently high throughout the entire optimization. We do so by formulating constraints on the minimum and maximum (solid) angles of triangular and tetrahedral mesh cells. After computing a mesh deformation, the nonlinear mesh quality constraints are incorporated into the optimization by projecting the mesh deformation according to the gradient projection method \cite{Rosen1960gradient,Rosen1961gradient} on the linear subspace of the currently active constraints, deforming the mesh along this subspace, and performing a correction step which ensures that all previously active constraints are active again. Additionally, the usual line search procedure is slightly modified to ensure that no new constraints are violated to ensure the feasibility of the method. Our method is semi-discrete in the sense that, first, a deformation is computed based on the underlying infinite-dimensional, continuous setting, and afterwards, this deformation is interpreted as a finite-dimensional deformation of the mesh and modified according to the mesh quality constraints. This makes our approach very flexible in the sense that it is completely independent from the way a mesh deformation is computed. Previous approaches, although they have been shown to enhance the mesh quality during shape optimization, cannot, to the best of our knowledge, guarantee a lower bound on the mesh's quality. In contrast, our method, which incorporates the mesh quality into the optimization in the form of constraints, gives a guarantee that the mesh quality is always above a user-specified threshold during the entire optimization. To the best of our knowledge, such an approach has not yet been considered in the literature.

We investigate our proposed method numerically for three examples. The first is an academic example of drag minimization of an obstacle in Stokes flow, which has been used extensively as benchmark problem in the previous literature \cite{Iglesias2018Two,Mueller2021novel,Blauth2021Nonlinear}. The second is concerned with the optimization of a flow in a pipe governed by the Navier-Stokes equations and has also been considered in the previous literature \cite{Schmidt2010Efficient,Ham2019Automated,Blauth2021Nonlinear}. The final one is a large-scale industrial shape optimization problem for improving the separation efficiency of a structured packing for distillation which is taken from our previous work \cite{Blauth2025CFD}. The numerical results show that our novel approach of incorporating the mesh quality directly into the shape optimization yields promising results and a significant improvement of the mesh quality throughout the shape optimization. 

This paper is structured as follows: In Section~\ref{sec:preliminaries}, we recall the gradient projection method for finite-dimensional nonlinear optimization problems and briefly present some basic results from shape calculus. In Section~\ref{sec:gradient_projection}, we propose our novel gradient projection method for shape optimization which incorporates the mesh quality constraints. Finally, in Section~\ref{sec:numerics} we analyze the proposed method numerically and show its great performance for solving both academic and industrial shape optimization problems.

\section{Preliminaries}
\label{sec:preliminaries}

In this section, we first recall the gradient projection method for finite-dimensional optimization problems, which is the basis of our gradient projection method for the mesh quality constraints. Additionally, we present some basic results from shape calculus.

\subsection{The Gradient Projection Method for Nonlinear Optimization}
\label{ssec:rosen_finite}

Let us start with briefly presenting the gradient projection method of Rosen (cf. \cite{Rosen1960gradient,Rosen1961gradient}), following \cite{Luenberger2021Linear}. Note that we only present an overview over the method, for a more detailed analysis we refer the reader, e.g, to \cite{Luenberger2021Linear,Rosen1960gradient,Rosen1961gradient}. 

We consider a constrained, nonlinear, and finite-dimensional optimization problem of the form
\begin{equation}
	\label{eq:constrained_problem}
	\begin{aligned}
		&\min_{x\in \R^d} f(x)\\
		\text{subject to } 
		&\begin{alignedat}[t]{2}
			g_i(x) &= 0 \quad &&\text{ for all } i \in \mathcal{E},\\
			g_j(x) &\leq 0 \quad &&\text{ for all } j \in \mathcal{I},
		\end{alignedat} 
	\end{aligned}
\end{equation}
where $f \in C^1(\R^d; \R)$ with $d\in \mathbb{N}_{>0}$ is some cost functional to be optimized, $\mathcal{E}$ and $\mathcal{I}$ are (finite and disjoint) sets of indices, $g_i \in C^1(\R^d; \R)$ for $i \in \mathcal{E}$ denote the equality, and $g_j \in C^1(\R^d; \R)$ for $j \in \mathcal{I}$ denote the inequality constraints. Throughout this section, we use the index $i$ for the equality and $j$ for the inequality constraints. We call a point $x \in \R^d$ feasible if $g_i(x) = 0$ for all $i \in \mathcal{E}$ and $g_j(x) \leq 0$ for all $j \in \mathcal{I}$. Moreover, at some feasible point $x\in \R^d$, we define the active set $\mathcal{A} = \mathcal{A}(x)$ as the set of all indices whose constraints are active at $x$, i.e.,
\begin{equation*}
	\mathcal{A}(x) = \Set{k \in \mathcal{E} \cup \mathcal{I} | g_k(x) = 0} = \mathcal{E} \cup \Set{j \in \mathcal{I} | g_j(x) = 0}.
\end{equation*}

Starting from some feasible point $x_0$, the iterates $x_n$ of the gradient projection method are generated as follows. First, the gradient $\grad f(x_n)$ as well as the active set $\mathcal{A}(x_n)$ are computed and we denote the cardinality of $\mathcal{A}(x_n)$ by $q = \abs{\mathcal{A}(x_n)}$. Then, we compute the gradient of all currently active constraints at $x_n$, i.e., $\grad g_k(x_n)$ for $k \in \mathcal{A}(x_n)$ and create the matrix $A = A(x_n) \in \R^{q \times d}$, which is the Jacobian of the active constraints, i.e.,
\begin{equation}
	\label{eq:matrix_A}
	A_{i,j} = \frac{\partial g_i}{\partial x_j}(x_n) \quad \Leftrightarrow \quad A = [\grad g_k(x_n) \transposed]_{k\in \mathcal{A}(x_n)}.
\end{equation}
The main idea of the gradient projection method is to project some search direction $s \in \R^d$, which is usually given by $s = -\nabla f(x_n)$, to the tangent space of the currently active constraints, i.e., to compute a projected search direction $p \in \R^d$ so that $A p = 0$. This projection can be computed as (cf. \cite{Luenberger2021Linear})
\begin{equation}
	\label{eq:projected_direction}
	p = s - A\transposed \lambda,
\end{equation}
where $\lambda$ solves
\begin{equation}
	\label{eq:lagrange_multiplier}
	A A\transposed \lambda = A s.
\end{equation}
Hence, the projected search direction can be formally written as
\begin{equation}
	\label{eq:rosen_projection_matrix}
	p = s - A\transposed \left( A A\transposed \right)^{-1} A s = \left( I - A\transposed \left( A A\transposed \right)^{-1} A \right) s = P s,
\end{equation}
where $I \in \R^{d\times d}$ denotes the identity matrix in $\R^d$ and $P = I - A\transposed \left( A A\transposed \right)^{-1} A \in \R^{d\times d}$ is the projection matrix at $x_n$. It is easy to see that $A p = 0$ holds. 

\begin{Remark}
	For the simplicity of presentation we assume that the gradients of all active constraints are linearly independent, which is also known as the linear independence constraint qualification (LICQ) (see, e.g., \cite{Nocedal2006Numerical}). This ensures that the projection matrix $P = I - A\transposed (A A\transposed)^{-1} A$ is well-defined. In case the gradients of the active constraints are linearly dependent, the Moore-Penrose pseudoinverse can be used instead \cite{Klooster2017Optimizing} which leads to the modified projection matrix $\hat{P} = I - A^+ A$, where $A^+$ is the Moore-Penrose pseudoinverse of $A$. Finally, we note that for the numerical experiments in Section~\ref{sec:numerics}, the LICQ did hold for all considered examples.
\end{Remark}

The projected search direction is used to compute a new trial iterate $y = x_n + t p$, where $t > 0$ is some trial step size. As the gradient projection is an active set method, it ensures that all constraints which were active at $x_n$ are also active at the next iterate. If all constraints were linear, then this would indeed hold for $y$. However, due to the nonlinearity of the constraints, this is, in general, not true. For this reason, $y$ has to be projected back to the surface of the active constraints, which can be done with the help of Newton's method (with a frozen Jacobian) and works as follows. Starting with $y_0 = y$ we compute iterates $y_l$ as follows
\begin{equation*}
	y_{l+1} = y_l - A\transposed \left( A A\transposed \right)^{-1} h(y_l),
\end{equation*}
where $h \colon \R^d \to \R^q$ with $h_k = g_k$ for all $k \in \mathcal{A}(x_n)$ denotes the function of all active constraints at $x_n$. This method converges to some $y^*$ for a sufficiently close initial guess, i.e., if the step size $t$ is sufficiently small. Moreover, all constraints that were active at $x_n$ are also active at $y^*$. When needed, we explicitly denote the dependence of $y^*$ on $x_n$ and $t$ as $y^* = y^*(x_n, t)$. In Figure~\ref{fig:projection}, the projection of the negative gradient and back-projection onto the set of active constraints is illustrated.

However, the point $y^*$ computed this way may still not be feasible since it could violate some inequality constraints which were inactive at $x_n$ but are now violated by $y^*$, as the constraints that were inactive previously have not been considered at all so far. Due to the continuity of the constraint functions, there exists some step size $t^*>0$ so that $y^*(x_n, t^*)$ does not violate any previously inactive constraint and at least one additional constraint is active at $y^*(x_n, t^*)$. Hence, the step size $t^*$ can be computed, e.g., with a bisection approach. Afterwards, a suitable step size for the next iterate $x_{n+1}$ can be computed, e.g., with an Armijo line search (see, e.g., \cite{Hinze2009Optimization}). This procedure is illustrated in Figure~\ref{fig:interpolation}.

\begin{figure}[t]
	\centering
	\newcommand{\lw}{1.25}
	\newcommand{\radius}{3}
	\renewcommand{\angle}{30}
	\newcommand{\opac}{0.33}
	\begin{subfigure}[t]{0.48\textwidth}
		\centering
		\begin{tikzpicture}
			\coordinate (xn) at (0,\radius);

			\fill[c6, opacity=\opac, pattern=north west lines, pattern color=c6] (xn) -- (0,{\radius*sin(\angle)}) -- ({\radius*cos(\angle)},{\radius*sin(\angle)}) arc[start angle=\angle, end angle=90, radius=3cm];
			\fill[c6, opacity=\opac, pattern=north west lines, pattern color=c6] (xn) -- (0,{\radius*sin(\angle)}) -- ({-\radius*cos(\angle)},{\radius*sin(\angle)}) arc[start angle=180-\angle, end angle=90, radius=3cm];
			
			\draw [c6, line width=\lw] (xn) arc (90:30:3);
			\draw [c6, line width=\lw] (xn) arc (90:150:3);
			\node [c6, rotate=40] at (-2.4, -0.75+\radius) {$g_1(x) = 0$};
			\node [color=c6, below=0.75cm] at (xn) {$g_1(x) \leq 0$};

			\draw [-latex, c0, line width=\lw] (xn) -- +(2,2) node[midway,sloped,above,align=center] {$-\nabla f(x_n)$};
			\draw [-latex, c1, line width=\lw] (xn) -- +(2,0) node[midway,sloped,above,align=center] {$p$};
			\draw [dotted, c1, line width=\lw] (-3, \radius) -- (3, \radius);
			\draw [dashed, line width=\lw] (2,2+\radius) -- (2,\radius);
			
			\coordinate (y) at (2.5, \radius);
			
			\draw [dotted, line width=\lw] (2.5,\radius) -- (2.5,-1.3+\radius);
			
			\fill (y) circle (2.5pt);
			\node [above right] at (y) {$y = x_n + tp$};
			
			\coordinate (ystar) at (2.5, -1.3+\radius);
			\fill (ystar) circle (2.5pt);
			\node [right] at (ystar) {$y^*(x_n, t)$};
			
			\fill (xn) circle (2.5pt);
			\node [below] at (xn) {$x_n$};
		\end{tikzpicture}
		\caption{Projection onto the active constraints: At $x_n$, constraint $g_1$ is active. The negative gradient $-\nabla f(x_n)$ is first projected to the tangent space of $g_1$ to produce the search direction $p$. A tentative step $y$ is generated as $y=x_n + tp$ for some step size $t$. This is then projected back to the active set with Newton's method, generating $y^*$.}
		\label{fig:projection}
	\end{subfigure}%
	\hfill%
	\begin{subfigure}[t]{0.48\textwidth}
		\centering
		\begin{tikzpicture}
			
			\coordinate (xn) at (0,\radius);
			\coordinate (y) at (2.5, \radius);
			\coordinate (ystar) at (2.5, \radius - 1.3);
			\coordinate (ybar) at (1.7, \radius);
			\coordinate (ybarstar) at (1.7, \radius-0.5);
			
			\fill[c6, opacity=\opac, pattern=north west lines, pattern color=c6] (xn) -- (0,{\radius*sin(\angle)}) -- ({\radius*cos(\angle)},{\radius*sin(\angle)}) arc[start angle=\angle, end angle=90, radius=3cm];
			\fill[c6, opacity=\opac, pattern=north west lines, pattern color=c6] (xn) -- (0,{\radius*sin(\angle)}) -- ({-\radius*cos(\angle)},{\radius*sin(\angle)}) arc[start angle=180-\angle, end angle=90, radius=3cm];
			
			\draw [line width=\lw, c6] (xn) arc (90:\angle:3);
			\draw [line width=\lw, c6] (xn) arc (90:180-\angle:3);
			
			\fill [c5, opacity=\opac, pattern=north east lines, pattern color=c5] (-3,5) -- plot [smooth] coordinates {(0.5, 5) (0.5,4) (1.75,2.25) (0, {\radius*sin(\angle)})} -- (-3, {\radius*sin(\angle)}) -- cycle;
			\draw [line width=\lw, c5] plot [smooth] coordinates {(0.5, 5) (0.5,4) (1.75,2.25) (0, {\radius*sin(\angle)})};
			
			\draw [-latex, c0, line width=\lw] (xn) -- +(2,2) node[near end,sloped,above] {$-\nabla f(x_n)$};
			\draw [dotted, c1, line width=\lw] (-3, \radius) -- (3, \radius);
			\draw [dotted, line width=\lw] (y) -- (ystar);
			\draw [dotted, line width=\lw] (ybar) -- (ybarstar);
			
			\fill (xn) circle (2.5pt);
			\fill (y) circle (2.5pt);
			\fill (ystar) circle (2.5pt);
			\fill (ybar) circle (2.5pt);
			\fill (ybarstar) circle (2.5pt);
			
			\node [below] at (xn) {$x_n$};
			\node [above right] at (y) {$y$};
			\node [right] at (ystar) {$y^*$};
			\node [above right] at (ybar) {$\bar{y}$};
			\node [below left] at (ybarstar) {$\bar{y}^*$};
			
			\node [color=c5, above left=1cm] at (xn) {$g_2(x) \leq 0$};
			\node [color=c6, below=0.75cm] at (xn) {$g_1(x) \leq 0$};
		\end{tikzpicture}
		\caption{Finding a feasible step size: At $x_n$, only constraint $g_1$ is active. The projection $y^*$ of a tentative step $y$ violates $g_2$. A bisection approach can be used to find a stepsize $\bar{t}$ so that $\bar{y} = x_n + \bar{t}p$ is projected to $\bar{y}^*$, where both $g_1$ and $g_2$ are active.}
		\label{fig:interpolation}
	\end{subfigure}
	\caption{Illustration of the gradient projection method.}
	\label{fig:rosen_illustration}
\end{figure}

\begin{algorithm2e}[!b]
	\KwIn{Initial feasible point $x_0$, initial step size $t_0$, stopping tolerance $\tau \in (0,1)$, maximum number of iterations $n_{\mathrm{max}} \in \mathbb{N}$, parameters $\sigma \in (0,1)$ and $\omega \in (0,1)$ for the Armijo rule \label{line:input}}
	\For{$n=0,1,2,\dots, n_\mathrm{max}$}{
		Evaluate $f(x_n)$ and compute its gradient $\grad f(x_n)$\\
		Compute some gradient-based search direction $s$, e.g., $s = -\grad f(x_n)$ \\
		Evaluate the constraints $g_k(x_n)$ for $k\in \mathcal{E} \cup \mathcal{I}$ and compute the set of active constraints $\mathcal{A}(x_n)$ \\
		Compute the approximate Lagrange multipliers $\lambda$ by solving \eqref{eq:lagrange_multiplier} \\
		Compute the projected search direction $p = s - A\transposed \lambda$ with \eqref{eq:projected_direction}\\
		\If{$\abs{p} < -\min \Set{0, \lambda_j}$}{
			Drop constraint $j$ from the set of active constraints for which $\gamma = -\lambda_j$ \\
			Set $p$ to the projected direction according to the modified active set
		}
		\ElseIf{$p = 0$ \textup{ and } $\lambda_j > 0$ \textup{ for all } $j \in \mathcal{A}(x_n) \cap \mathcal{I}$}{
			Stop with minimizer $x_n$\\	
		}
		Compute $y^*(x_n, t)$ \\
		\If{\textup{additional constraints are violated at} $y^*(x_n, t)$}{
			Compute the maximum feasible step size $t^*$ at $x_n$ with a bisection \\
			Set $t = t^*$ \\
		}
		\While{$f(y^*(x_n, t)) > f(x_n) + \sigma t (p, \grad f(x_n)) $}{
			Decrease the step size $t = \omega t$ \\
		}
		Compute the new iterate $x_{n+1} = y^*(x_n, t)$ \\
		Increase the step size for the next iteration: $t = \nicefrac{t}{\omega}$
	}
	\caption{Gradient projection method.}
	\label{algo:rosen}
\end{algorithm2e}

Up until now, we have only added additional constraints to the active set. However, there are some circumstances under which constraints can also be dropped. After having computed $\lambda$ for projecting the search direction in \eqref{eq:lagrange_multiplier}, we define $\gamma = -\min \Set{\lambda_j | j \in \mathcal{A}(x_n) \cap \mathcal{I}} \cup \Set{0}$. If $\abs{Ps} \geq \gamma$, we proceed as detailed above and use $p = Ps$ as search direction. On the other hand, if $\abs{Ps} < \gamma$, then $\gamma > 0$ holds if $s \neq 0$ and we denote by $j \in \mathcal{A}(x_n) \cap \mathcal{I}$ the index for which $-\lambda_j = \gamma$. In this case, we define a modified active set $\bar{\mathcal{A}}(x_n)$ by dropping constraint $j$, i.e., $\bar{\mathcal{A}}(x_n) = \mathcal{A}(x_n) \setminus \{j\}$. With this modified active set, we define $\bar{A} = [\nabla g_k(x_n)\transposed]_{k\in \bar{\mathcal{A}}(x_n)}$. Then, we use as search direction $\bar{p} = s - \bar{A}\transposed \left( \bar{A} \bar{A}\transposed \right)^{-1} \bar{A}s$, i.e., the search direction $\bar{p}$ is the one obtained from dropping constraint $j$ from the active set. In particular, when choosing $s = -\grad f(x_n)$, this constraint dropping procedure ensures that if $p = Ps = 0$ at $x_n$, then we have $\grad f(x_n) + \lambda\transposed A = 0$ due to \eqref{eq:projected_direction} and, additionally, $\lambda_j > 0$ for all $j \in \mathcal{A}(x_n) \cap \mathcal{I}$. These are just the first order necessary conditions for a minimizer of the constrained problem \eqref{eq:constrained_problem} if the vector of $\lambda$ is padded by zero for all inactive inequality constraints, as $x_n$ is also feasible by construction.

A brief summary of the gradient projection method can be found in Algorithm~\ref{algo:rosen}. As we use the algorithm for the numerical solution of constrained optimization problems, all equations are only solved approximately up to a specified tolerance. Moreover, also constraints can only be satisfied up to a tolerance $\varepsilon$ and constraint $k$ is considered active at $x$ if $\abs{g_k(x)} \leq \varepsilon$.

\subsection{Shape Calculus}
\label{ssec:shape_calculus}

As we want to apply the previously described gradient projection method to solve shape optimization problems, we now briefly recall some basic results in this field. A shape optimization problem constrained by a PDE can be written as
\begin{equation*}
	\begin{aligned}
		&\min_{\Omega \in \admissiblegeom} \mathcal{J}(\Omega, u) \quad \text{ subject to } \quad e(\Omega,u) = 0,
	\end{aligned}
\end{equation*}
where $\mathcal{J}$ is a cost functional depending on the shape of the domain $\Omega$, which is sought in the set of all admissible shapes $\admissiblegeom$, and the state variable $u$. The PDE constraint is denoted by $e(\Omega, \cdot) \colon U(\Omega) \to V(\Omega)^*$ and is usually considered in the weak form
\begin{equation*}
	\text{Find } u \in U(\Omega) \text{ such that } \quad \dual{e(\Omega, u)}{v}{V(\Omega)^*,V(\Omega)} = 0 \quad \text{ for all } v \in V(\Omega).
\end{equation*}
Here, $V(\Omega)^*$ denotes the topological dual space of $V(\Omega)$ and $\dual{\varphi}{v}{V(\Omega)^*,V(\Omega)}$ is the duality pairing between $\varphi \in V(\Omega^*)$ and $v\in V(\Omega)$. We assume that the state equation has a unique solution $u = u(\Omega)$ for all $\Omega \in \admissiblegeom$ which satisfies $e(\Omega, u(\Omega)) = 0$. Hence, we can introduce the equivalent reduced problem
\begin{equation*}
	\min_{\Omega \in \admissiblegeom} J(\Omega) = \mathcal{J}(\Omega, u(\Omega)),
\end{equation*}
where the PDE constraint is formally eliminated due to the reduced cost functional $J$. 


Shape calculus can be used to analyze such shape optimization problems. For a detailed introduction to this topic, we refer the reader, e.g., to \cite{Delfour2011Shapes}. We use the so-called speed method to derive the variation of the cost functional w.r.t.\ infinitesimal changes of the domain's shape. Let $\holdall \subset \R^d$ be an open and bounded hold-all domain and consider some $\Omega \subset \holdall$. Let $\vectorfield \in C^k_0(\holdall;\R^d)$ for some $k\geq 1$, i.e., the space of all $k$ times continuously differentiable vector fields with compact support in $\holdall$. We consider the evolution of some point $x_0 \in \Omega$ under the flow defined by the solution of the initial value problem
\begin{equation}
	\label{eq:ivp}
	\dot{x}(t) = \vectorfield(x(t)), \quad x(0) = x_0.
\end{equation}
Due to the regularity of $\vectorfield$, the above problem has a unique solution $x(t)$ for $t\in [0, \tau]$ if $\tau > 0$ is sufficiently small (cf.\ \cite{Delfour2011Shapes}). Hence, we define the flow $\flow^\vectorfield_t$ of $\vectorfield$ as follows
\begin{equation*}
	\flow^\vectorfield_t \colon \holdall \to \holdall; \quad x_0 \mapsto \flow^\vectorfield_t x_0 = x(t).
\end{equation*}
Now, we can define the shape derivative of a shape functional.
\begin{Definition}
	\label{def:shape_derivative}
	Let $\mathcal{S} \subset \Set{\Omega | \Omega \subset \holdall}$, $J \colon \mathcal{S} \to \R$ be a shape functional, $\Omega \in \mathcal{S}$, and $\vectorfield \in C^k_0(\holdall; \R^d)$ for $k\geq 1$ with associated flow $\flow^\vectorfield_t$. Additionally, we assume that $\flow^\vectorfield_t(\Omega) \in \mathcal{S}$ for all $t \in [0, \tau]$ with some sufficiently small $\tau > 0$. The Eulerian semi-derivative of $J$ is given by the following limit
	\begin{equation*}
		dJ(\Omega)[\vectorfield] := \lim\limits_{t\searrow 0} \frac{J(\flow^\vectorfield_t (\Omega)) - J(\Omega)}{t} = \left. \frac{d}{d t} J(\flow^\vectorfield_t(\Omega)) \right|_{t=0^+}
	\end{equation*}
	if the limit exists. Moreover, let $\Xi$ be a topological vector subspace of $C^\infty_0(\holdall;\R^d)$. Then, the shape functional $J$ is called shape differentiable at $\Omega \subset \holdall$ if it has an Eulerian semi-derivative at $\Omega$ in all directions $\vectorfield \in \Xi$ and, additionally, the mapping $dJ(\Omega) \colon \Xi \to \R;\ \vectorfield \mapsto dJ(\Omega)[\vectorfield]$ is linear and continuous. In this case, we call $dJ(\Omega)$ the shape derivative of $J$ at $\Omega$ (w.r.t.\ $\Xi$).
\end{Definition}

%
\begin{Remark}
	In general, the shape derivatives of PDE constrained optimization problems can be calculated efficiently using the so-called adjoint approach which entails the solution of an adjoint equation. For this reason, to compute the shape derivative for such problems, both the state and adjoint equation have to be solved. For more details regarding the adjoint approach we refer the reader, e.g., to \cite{Hinze2009Optimization}.
\end{Remark}

The sensitivity information in the shape derivative can now be used to compute the so-called gradient deformation which can be used as part of gradient-based shape optimization algorithms. To do so, we consider the shape derivative acting on vector fields in $H^1_0(\holdall; \R^d)$. This is reasonable as the shape derivative is defined for vector fields in $C^\infty_0(\holdall; \R^d)$ and this space is dense in $H^1_0(\holdall; \R^d)$. Hence, the shape derivative can be uniquely extended to vector fields in $H^1_0(\holdall; \R^d)$, see, e.g., \cite{Alt2016Linear}. For a more detailed discussion of this topic we refer the reader to \cite[Chapter~5.2]{Allaire2021Shape}. In particular, we consider the case where some parts of the boundary $\Gamma = \partial \Omega$ of the initial domain $\Omega$ are fixed during the optimization and denote these fixed parts by $\Gamma_\mathrm{fix}$. We consider the subspace 
\begin{equation*}
	H^1_{\Gamma_\mathrm{fix}}(\holdall;\R^d) = \Set{\vectorfield \in H^1_0(\holdall;\R^d) | \vectorfield = 0 \text{ on } \Gamma_\mathrm{fix}}.
\end{equation*}
\begin{Definition}
	\label{def:gradient_deformation}
	Let $\Omega \subset \holdall$ and $J$ be shape differentiable at $\Omega$. Let $a_\Omega\colon H^1_{\Gamma_\mathrm{fix}}(\holdall;\R^d) \times H^1_{\Gamma_\mathrm{fix}}(\holdall;\R^d) \to \R$ be a symmetric, continuous, and coercive bilinear form. The gradient deformation $\gradientdefo \in H^1_{\Gamma_\mathrm{fix}}(\holdall;\R^d)$ (w.r.t. $a_\Omega$) is defined as the unique solution of the variational problem
	\begin{equation}
		\label{eq:gradient_deformation}
		\text{Find } \gradientdefo \in H^1_{\Gamma_\mathrm{fix}}(\holdall;\R^d) \text{ such that } \quad a_\Omega(\gradientdefo, \vectorfield) = dJ(\Omega)[\vectorfield] \quad \text{ for all } \vectorfield \in H^1_{\Gamma_\mathrm{fix}}(\holdall;\R^d),
	\end{equation}
	which exists due to the Lax-Milgram Lemma.
\end{Definition}
In particular, if $\gradientdefo \neq 0$, it holds that $dJ(\Omega)[-\gradientdefo] = -a_\Omega(\gradientdefo, \gradientdefo) < 0$ due to the coercivity of $a_\Omega$. This implies that the negative gradient deformation is a direction of descent in the sense that an infinitesimal perturbation of $\Omega$ in direction of $-\gradientdefo$ yields a decrease of the shape functional $J$. This fact is used for gradient-based shape optimization algorithms as described, e.g., in \cite{Schulz2016Efficient, Blauth2021Nonlinear}.

\begin{Remark}
	The choice of the bilinear form $a_\Omega$ is crucial for the mesh quality when the mesh is deformed with the gradient deformation. We follow the widely-used approach of using the linear elasticity equations for the bilinear form $a_\Omega$ (see, e.g., \cite{Blauth2021Nonlinear,Etling2020First,Schulz2016Efficient}), which are given by
	\begin{equation}
		\label{eq:linear_elasticity}
		\begin{aligned}
			&a_\Omega \colon H^1_{\Gamma_\mathrm{fix}}(\holdall; \R^d) \times H^1_{\Gamma_\mathrm{fix}}(\holdall; \R^d) \to \R;\\
			&(V, W) \mapsto a_\Omega(V, W) = \integral{\holdall} 2\mu_\mathrm{elas}\ \varepsilon(V):\varepsilon(W) + \lambda_\mathrm{elas}\ \mathrm{div}(V)\mathrm{div}(W) + \delta_\mathrm{elas} V\cdot W \dmeas{x},
		\end{aligned}
	\end{equation}
	where $\varepsilon(V) = \nicefrac{1}{2}(DV + DV\transposed)$ is the symmetric part of the Jacobian $DV$, $\lambda_\mathrm{elas}$ and $\mu_\mathrm{elas}$ are the so-called Lam\'e parameters, for which we assume $\mu_\mathrm{elas} > 0$ and $2\mu_\mathrm{elas} + d \lambda_\mathrm{elas} > 0$, and $A:B$ denotes the Frobenius inner product between $A, B \in \R^{d\times d}$. Moreover, $\delta_\mathrm{elas} \geq 0$ is a damping parameter, which is required to be positive in case $\Gamma_\mathrm{fix} = \emptyset$, where we have a pure Neumann problem. However, we note that there exist alternative approaches, such as using the (nonlinear) p-Laplacian to define the mesh deformations \cite{Mueller2021novel}.
\end{Remark}

\subsection{Gradient-Based Descent Algorithms for Shape Optimization}

\begin{algorithm2e}[t]
	\KwIn{Initial geometry $\Omega\iidx{0}$, initial step size $t\iidx{0}$, stopping tolerance $\tau \in (0,1)$, maximum number of iterations $n_{\text{max}}\in  \mathbb{N}$, parameters $\sigma \in (0, 1)$ and $\omega \in (0,1)$ for the Armijo rule}
	\For{n=0,1,2,\dots, $n_\textrm{max}$}{
		Compute the solution of the state and adjoint systems \label{lin:pde_solves} \\
		Compute the gradient deformation $\gradientdefo\iidx{n}$ by solving \eqref{eq:gradient_deformation} \label{lin:gradient_defo} \\
		\If{$\norm{\gradientdefo_n}{a_{\Omega_n}} \leq \tau \norm{\gradientdefo_{0}}{a_{\Omega_{0}}}$ \label{lin:stop_test}}{
			Stop with approximate solution $\Omega_{n}$
		}
		Compute a search direction $\mathcal{S}_{n}$, e.g., $\mathcal{S}_{n} = -\gradientdefo_{n}$ \label{lin:search_direction} \\
		\While{$J((I + t\mathcal{S}_n)\Omega_n) > J(\Omega_{n}) + \sigma t\ a_{\Omega_{n}}\left( \gradientdefo_{n}, \mathcal{S}_{n} \right)$ \label{lin:armijo}}{
			Decrease the step size: $t = \omega t$ \label{lin:armijo_fail} \\
		}
		Set $t_{n} = t$ and update the geometry via $\Omega_{n+1} = (I + t_n \mathcal{S}_n) \Omega_n$ \label{lin:update_geometry} \\
		Increase the step size for the next iteration: $t = \nicefrac{t_{n}}{\omega}$ \label{lin:increase_size} \\
	}
	\caption{Gradient-based descent algorithm for shape optimization.}
	\label{algo:shape_descent}
\end{algorithm2e}

Let us recall a gradient-based descent algorithm for shape optimization, following \cite{Blauth2021Nonlinear}. To do so, we first discuss how to update the shape of a domain with some vector field. While the speed method has nice theoretical properties, it is not so easy to use for the numerical solution of shape optimization problems. Instead, we consider the so-called perturbation of identity (cf.\ \cite{Blauth2021Nonlinear,Delfour2011Shapes}) which is defined as follows. Let $\Omega \subset \holdall$ be some domain and let $\vectorfield \in C^k_0(\holdall;\R^d)$. Then, we define the perturbed domain $\Omega_t$ with the perturbation of identity as
\begin{equation}
	\label{eq:perturbation_of_identity}
	\Omega_t = \left( I + t\vectorfield \right)\Omega = \Set{x + t\vectorfield(x) | x\in \Omega}.
\end{equation}
As the perturbation of identity presents an equivalent alternative to the speed method for computing first order shape derivatives \cite{Delfour2011Shapes}, we observe that a perturbation of identity in direction of the negative gradient deformation yields a decrease of the shape functional, analogous to our previous discussion. However, we note that the speed method and perturbation of identity do, in general, not coincide for higher order shape deriavtives, cf.~\cite[Chapter 4, Section 3.2]{Delfour2011Shapes}. 
With this, we have everything at hand and can present a general, gradient-based descent algorithm for shape optimization in Algorithm~\ref{algo:shape_descent} (cf.\ \cite{Blauth2021Nonlinear}). Note that in Line~\ref{lin:search_direction} of Algorithm~\ref{algo:shape_descent} a search direction is computed. For the classical method of steepest descent, one can choose $\mathcal{S}_n = -\gradientdefo_n$, i.e., use the negative gradient deformation as search direction. Other choices are possible and we refer the reader to \cite{Blauth2021Nonlinear} and \cite{Schulz2016Efficient}, where nonlinear conjugate gradient (CG) methods and limited-memory BFGS methods for shape optimization are introduced.

\section{Mesh Quality Constraints for Shape Optimization}
\label{sec:gradient_projection}

In this section, we combine the two topics presented in the previous section: The gradient projection method and shape optimization based on shape calculus. We do so to introduce constraints on the mesh quality which are then incorporated into the shape optimization.

\subsection{Discretization of the Shape Optimization Problem}
\label{ssec:discretization_shape_optimization}

To formulate constraints on the mesh quality, the shape optimization approach presented in the previous section has to be discretized. We assume that this is done with the finite element method, but other approaches, such as finite volume methods, are also available. To discretize the geometry, a computational mesh is used on which the corresponding PDE constraints are solved. We denote this mesh by $\Omega_h$. For the shape optimization approach presented previously, the domain $\Omega$ plays the role of the optimization variable. After discretization, $\Omega$ is represented by the mesh $\Omega_h$ and, in particular, the nodes of the computational mesh act as the optimization variables. 

Throughout this paper, we restrict our attention to simplicial meshes, i.e., two-dimensional triangular as well as three-dimensional tetrahedral meshes. An extension to other mesh types, such as two-dimensional quadrilateral or three-dimensional hexahedral meshes, is straightforward. The reason for this restriction is that the finite element software FEniCS \cite{Alnaes2015FEniCS, Logg2012Automated}, where we implement our method, does, in general, only support simplicial meshes. In the following, let $N \in \mathbb{N}$ be the number of mesh nodes and $M \in \mathbb{N}$ be the number of mesh cells. Moreover, we denote by $v \in \R^{dN}$, where $d=2,3$ is the dimensionality of the problem, the vector that contains the coordinates of the mesh nodes, i.e., $v = [x_1, y_1, x_2, y_2, \dots, x_N, y_N]\transposed$ in 2D and $v = [x_1, y_1, z_1, x_2, y_2, z_2, \dots, x_N, y_N, z_N]\transposed$ in 3D. To discretize the shape optimization approach, the computation of the gradient deformation \eqref{eq:gradient_deformation} is discretized with linear Lagrange elements. As the degrees of freedom (DoFs) of linear Lagrange elements are the point evaluations at the nodes of the mesh, the discretized gradient deformation $\gradientdefo_h$ can be identified with its DoF vector. This vector is denoted by $G_h \in \R^{dN}$ and contains the point evaluation of $\gradientdefo_h$ at the nodes of the mesh. For the sake of simplicity, we assume w.l.o.g.\ that the DoFs of the linear Lagrange elements and the nodes of the mesh are ordered in the same manner. If this is not the case, a suitable permutation for this identification has to be considered. Moreover, any discretized search direction $\mathcal{S}_h$ considered in Algorithm~\ref{algo:shape_descent} is (also) discretized with linear Lagrange elements and can be identified analogously with its DoF vector $S_h \in \R^{dN}$. If the shape $\Omega$ is updated with the perturbation of identity \eqref{eq:perturbation_of_identity} of some vector field $\mathcal{V}$ according to $\Omega_t = (I + t\mathcal{V})\Omega$, then the discretized analogue is given by the following: Let $\mathcal{V}_h$ be the discretized vector field corresponding to $\mathcal{V}$ and let $V_h$ be its DoF vector. The discretized perturbation of identity updates the mesh $\Omega_h$, which is identified by the vector its node coordinates $v$ as follows
\begin{equation*}
	(\Omega_h)_t = (I + t\mathcal{V}_h) \Omega_h \qquad \text{ where } \qquad v_t = v + t V_h
\end{equation*}
and $v_t \in \R^{dN}$ is the coordinate vector of the updated mesh $(\Omega_h)_t$. For this reason, the gradient deformation can be interpreted as the gradient of the cost functional w.r.t.\ the mesh coordinates in the discretized setting.

\begin{Remark}
	\label{rem:gradient_based_methods}
	For the discretized search direction $\mathcal{S}_h$ an obvious choice is given by $\mathcal{S}_h = -\mathcal{G}_h$, which leads to the gradient projection method. However, other choices for the search direction, such as shape deformations computed with the BFGS \cite{Schulz2016Efficient} or nonlinear CG \cite{Blauth2021Nonlinear} methods, are also possible and, usually, lead to more efficient solution algorithms. We refer to \cite{Blauth2021Nonlinear}, where the efficiency of the methods is compared numerically.
\end{Remark}


\subsection{Constraints on the Mesh Quality}
\label{ssec:mesh_quality_constraints}

Let us now define the mesh quality constraints we consider for the shape optimization.

\subsubsection{Two-Dimensional Triangular Meshes}

In the following, let $\tau$ be a non-degenerate triangle with nodes $v_i$, $v_j$, and $v_k$, where $v_m = [x_m, y_m]\transposed$ for $m=i,j,k$. Here, non-degenerate means that the triangle does have a non-zero area. Moreover, we denote by $\alpha_i$, $\alpha_j$, and $\alpha_k$ the triangles' angles at node $v_i$, $v_j$, and $v_k$, respectively. In the literature, there exist many quality measures for triangular cells, e.g., the aspect ratio, radius ratio, edge ratio or condition number. For an overview and analysis of these measures we refer the reader, e.g., to \cite{Pebay2003Analysis}. In this paper, we focus our attention on the angles of the triangle and want to ensure that these are bounded from below and above to avoid degeneration of the triangle. Note that this, implicitly, also establishes bounds on many triangular quality measures, as discussed in \cite{Pebay2003Analysis}. For this reason, we now propose the following mesh quality constraints:
\begin{equation}
	\label{eq:triangle_angle_constraint}
	\alpha_m \geq \alpha_\mathrm{thr} \qquad \text{ for } m = i,j,k,
\end{equation}
which is equivalent to saying that the minimum angle of the triangle is bounded from below by $\alpha_\mathrm{thr}$. Here, $\alpha_\mathrm{thr} \in (0, \nicefrac{\pi}{3})$ is a constant minimum angle threshold which is independent of the node coordinates but may be different for each mesh cell. We discuss the choice of $\alpha_\mathrm{thr}$ later in Section~\ref{ssec:choice_of_angle}. This constraint ensures that all triangular angles are bounded from below and, thus, prevents the degradation of the triangle. As we have that $\alpha_i + \alpha_j + \alpha_k = \pi$, it follows that
\begin{equation*}
	\alpha_m \leq \pi - 2\alpha_\mathrm{thr} \qquad \text{ for } m = i,j,k,
\end{equation*}
i.e., the maximum angle is bounded by $\pi - 2\alpha_\mathrm{thr}$ from above if \eqref{eq:triangle_angle_constraint} holds. For this reason, the constraints in \eqref{eq:triangle_angle_constraint} are sufficient to bound the minimum and maximum angle of the triangle and, therefore, they establish bounds on the mesh quality as discussed in \cite{Pebay2003Analysis}. In the framework of the gradient projection method from Section~\ref{ssec:rosen_finite}, the constraints \eqref{eq:triangle_angle_constraint} can be written in the form
\begin{equation}
	\label{eq:triangle_constraint_functions}
	g_m = \alpha_\mathrm{thr} - \alpha_m \leq 0 \qquad \text{ for } m=i,j,k.
\end{equation}
Hence, if our mesh consists of $M$ triangular cells, we consider a total of $3M$ mesh quality constraints.

Let us now compute the gradient of $g_m$ w.r.t.\ the coordinates of the mesh's nodes. Let the edge between $v_i$ and $v_j$ be denoted by $e_{ij}$, i.e., $e_{ij} = [x_i - x_j, y_i - y_j]\transposed$, then the angle $\alpha_i$ can be computed as
\begin{equation*}
	\alpha_i = \arccos\left( \frac{e_{ij} \cdot e_{ik}}{\abs{e_{ij}} \abs{e_{ik}}} \right).
\end{equation*}
Using elementary calculations, the gradient of the constraint functions $g_m$ of \eqref{eq:triangle_constraint_functions} w.r.t.\ the mesh coordinates is given by
\begin{equation*}
	\begin{aligned}
		\frac{\partial g_i}{\partial v_i} &= -\frac{\partial \alpha_i}{\partial v_i} = -T(e_{ij}, e_{ik}) - T(e_{ik}, e_{ij}), \\
		\frac{\partial g_i}{\partial v_j} &= -\frac{\partial \alpha_i}{\partial v_j} = T(e_{ij}, e_{ik}), \\
		\frac{\partial g_i}{\partial v_k} &= -\frac{\partial \alpha_i}{\partial v_k} = T(e_{ik}, e_{ij}),
	\end{aligned}
\end{equation*}
where we use the notation
\begin{equation}
	\label{eq:derivative_notation_2D}
	\frac{\partial g_i}{\partial v_j} = \left[ \frac{\partial g_i}{\partial x_j}, \frac{\partial g_i}{\partial y_j} \right]\transposed
\end{equation}
and the vector $T$ is given by
\begin{equation}
	\label{eq:triple_product}
	T(a,b) = \frac{1}{\abs{a}} \frac{a \times (a\times b)}{\abs{a \times (a\times b)}}.
\end{equation}
Here, the vectors $a, b \in \R^2$ are padded with a zero in the third component for the cross products. In particular, it is easy to see that $a\times (a\times b)$ has a zero in the third component, so that $T(a,b)$ can be interpreted as a vector in $\R^2$, again. Obviously, it holds that $\frac{\partial g_i}{\partial v_m} = 0$ for all nodes $v_m$ which are not part of the triangle $\tau$. Hence, the gradient of the triangular mesh quality constraints given above is sparse and only has 6 nonzero entries for each constraint.

\subsubsection{Three-Dimensional Tetrahedral Meshes}

Let $\tau$ be a non-degenerate tetrahedron whose nodes are denoted by $v_i$, $v_j$, $v_k$, and $v_l$, where $v_m = [x_m, y_m, z_m]\transposed$ for $m=i,j,k,l$. Again, non-degenerate means that the volume of the tetrahedron is non-zero. We denote by $\alpha_m$ the solid angle at $v_m$, which is given by the surface area of the spherical triangle that is formed by projecting each of the remaining nodes to the unit sphere around $v_m$ (cf.\ \cite{Liu1994Relationship}). In the literature, there exist plenty of quality measures for tetrahedra, e.g., the aspect ratio, radius ratio or mean ratio. For an overview, the reader is referred, e.g., to \cite{Liu1994Relationship}. Particularly, we note that the minimum solid angle of a tetrahedron can be used as quality measure and in \cite{Liu1994Relationship} it is shown that this is, in a weak sense, equivalent to the radius and mean ratio. For this reason, we restrict our attention to the solid angles of the tetrahedron. As before, our aim is to bound these from above and below to avoid degeneration of the tetrahedral mesh cells. To do so, we propose the following mesh quality constraints
\begin{equation}
	\label{eq:tetrahedral_mesh_quality_constraints}
	\alpha_m \geq \alpha_\mathrm{thr} \qquad \text{ for } m = i,j,k,l,
\end{equation}
where $\alpha_\mathrm{thr} \in (0, \alpha^*)$ is a constant minimum solid angle treshold which might be different for each element, analogously to before, and $\alpha^* = \arccos(\nicefrac{23}{27}) \approx \qty{0.55129}{\radian}$ is the solid angle of a regular tetrahedron. Similar to the triangular case, it holds that $\sum_{m\in {i,j,k,l}} \alpha_m \leq 2\pi$ \cite{Liu1994Relationship}, so that \eqref{eq:tetrahedral_mesh_quality_constraints} also gives an upper bound on the solid angles, namely
\begin{equation*}
	\alpha_m \leq 2\pi - 3\alpha_\mathrm{thr} \qquad \text{ for } m = i,j,k,l,
\end{equation*}
i.e., the maximum solid angle is bounded by $2\pi - 3\alpha_\mathrm{thr}$ if \eqref{eq:tetrahedral_mesh_quality_constraints} holds. For the gradient projection method, the constraints \eqref{eq:tetrahedral_mesh_quality_constraints} can be written as
\begin{equation}
	\label{eq:tetrahedron_quality_constraint_functions}
	g_m = \alpha_\mathrm{thr} - \alpha_m \leq 0 \qquad \text{ for } m = i,j,k,l.
\end{equation}
Using \eqref{eq:tetrahedral_mesh_quality_constraints} in a mesh consisting of $M$ tetrahedrons yields a total of $4M$ mesh quality constraints.

To evaluate this numerically, we make use of the following relation from \cite{Liu1994Relationship}. Let $\delta_{i,j}$ be the dihedral angle corresponding to edge $e_{ij} = [x_i - x_j, y_i - y_j, z_i - z_j]\transposed$. Then, it holds that 
\begin{equation*}
	\alpha_i = \delta_{ij} + \delta_{ik} + \delta_{il} - \pi. 
\end{equation*}
The dihedral angles can be computed with the formula
\begin{equation*}
	\delta_{ij} = \arccos\left( \frac{(e_{ij}\times e_{ik})\cdot (e_{ij}\times e_{il}) }{\abs{e_{ij} \times e_{ik}}  \abs{e_{ij}\times e_{il}}} \right) = \arccos \left( \frac{n_{ijk} \cdot n_{ijl}}{\abs{n_{ijk}} \abs{n_{ijl}}} \right),
\end{equation*}
where $n_{ijk} = e_{ij} \times e_{ik}$ is the (unsigned and unscaled) normal vector of the face spanned by $v_i$, $v_j$, and $v_k$. Note that this formulation yields $0 \leq \delta_{ij} \leq \pi$ as well as $\delta_{ij} = \delta_{ji}$. Now, using the results for the triangular case with an application of the chain rule yields the following derivative of the constraint functions w.r.t.\ the node coordinates
\begin{equation*}
	\begin{aligned}
		\frac{\partial g_i}{\partial v_i} = -\frac{\partial \alpha_i}{\partial v_i} &= \left( \mathrm{T}(n_{ijk}, n_{ijl}) - \mathrm{T}(n_{ikj}, n_{ikl}) \right) \times e_{jk} + \left( \mathrm{T}(n_{ijl}, n_{ijk}) - \mathrm{T}(n_{ilj}, n_{ilk}) \right) \times e_{jl} \\
		&\quad + \left( \mathrm{T}(n_{ikl}, n_{ikj}) - \mathrm{T}(n_{ilk}, n_{ilj}) \right) \times e_{kl} \\
		%
		%
		\frac{\partial g_i}{\partial v_j} = -\frac{\partial \alpha_i}{\partial v_j} &= \left( \mathrm{T}(n_{ikj}, n_{ikl}) - \mathrm{T}(n_{ijk}, n_{ijl}) \right) \times e_{ik} + \left( \mathrm{T}(n_{ilj}, n_{ilk}) - \mathrm{T}(n_{ijl}, n_{ijk}) \right) \times e_{il} \\
		\frac{\partial g_i}{\partial v_k} = -\frac{\partial \alpha_i}{\partial v_k} &= \left( \mathrm{T}(n_{ijk}, n_{ijl}) - \mathrm{T}(n_{ikj}, n_{ikl}) \right) \times e_{ij} + \left( \mathrm{T}(n_{ilk}, n_{ilj}) - \mathrm{T}(n_{ikl}, n_{ikj}) \right) \times e_{il} \\
		\frac{\partial g_i}{\partial v_l} = -\frac{\partial \alpha_i}{\partial v_l} &= \left( \mathrm{T}(n_{ijl}, n_{ijk}) - \mathrm{T}(n_{ilj}, n_{ilk}) \right) \times e_{ij} + \left( \mathrm{T}(n_{ikl}, n_{ikj}) - \mathrm{T}(n_{ilk}, n_{ilj}) \right) \times e_{ik}
	\end{aligned}
\end{equation*}
where $T$ is defined as before in \eqref{eq:triple_product}.
Analogously to before, we use the notation
\begin{equation}
	\label{eq:derivative_notation_3D}
	\frac{\partial g_i}{\partial v_j} = \left[ \frac{\partial g_i}{\partial x_j}, \frac{\partial g_i}{\partial y_j}, \frac{\partial g_i}{\partial z_j} \right]\transposed.
\end{equation}
Similarly to before we observe that $\frac{\partial g_i}{\partial v_m} = 0$ if $v_m$ is not part of the tetrahedron $\tau$. Again, this means that the gradient of the tetrahedral mesh quality constraints is sparse and has only 12 nonzero entries for each constraint.

Allthough it might seem simpler to directly use the dihedral angle for the tetrahedral mesh quality constraints, we note that this is, in fact, not a suitable mesh quality measure as discussed in, e.g., \cite{Dompierre1998Proposal}. There, it is shown that there exist certain configurations where a tetrahedron can degenerate while all its dihedral angles are bounded away from zero. Such problems do not occur for the solid angle as this is a proper mesh quality measure \cite{Liu1994Relationship}. For this reason, we restrict our attention to the solid angle for the three-dimensional case.

\begin{Remark}
	For the sake of simplicity, we propose to bound all (solid) angles instead of only considering the minimum (solid) angle of each cell. This facilitates a slightly easier implementation and avoids technical difficulties associated with the non-differentiability of the minimum-function.
\end{Remark}

\begin{Remark}
	\label{rem:number_of_active_constraints}
	The presented mesh quality constraints consider a total of $3M$ or $4M$ constraints, respectively, for the discretized shape optimization problem, which, however, only has $2N$ or $3N$ optimization variables, respectively, where $M \gg N$. From a formal standpoint this leads to an over-constrained optimization problem which is, in general, not solvable. However, as we will see in our numerical experiments in Section~\ref{sec:numerics}, only a small portion of the considered constraints actually become active during the optimization so that our proposed approach works very well in practice. 
\end{Remark}

\subsection{Implementation Details}

Let us now discuss some implementation details of our method. First, we discuss how to combine the gradient projection method with shape optimization to incorporate the mesh quality constraints proposed previously. Afterwards, we introduce additional constraints for fixed boundaries and discuss the choice of the minimum angle threshold.

\subsubsection{The Gradient Projection Method for Mesh Quality Constraints}

To incorporate the previously defined constraints \eqref{eq:triangle_constraint_functions} and \eqref{eq:tetrahedron_quality_constraint_functions} into our shape optimization procedure of Algorithm~\ref{algo:shape_descent}, we use the gradient projection method described in Section~\ref{ssec:rosen_finite}. As discussed in Section~\ref{ssec:discretization_shape_optimization}, the optimization variables of the discretized shape optimization problem are the nodes of the mesh and the discretized search directions and gradient deformation of Algorithm~\ref{algo:shape_descent} can be identified with their respective DoF vectors. In particular, at first, we follow the procedure outlined in Algorithm~\ref{algo:shape_descent} to compute a suitable (discretized) search direction $\mathcal{S}_h$ for our problem, which we identify with its DoF vector $S_h$. As stated in Remark~\ref{rem:gradient_based_methods}, suitable gradient-based search directions can be computed with the gradient descent, nonlinear CG or BFGS methods. Then, we proceed as discussed in Section~\ref{ssec:rosen_finite} and project the search direction using \eqref{eq:projected_direction} and \eqref{eq:lagrange_multiplier}. To form the matrix $A$ of \eqref{eq:matrix_A}, the derivatives of the constraint functions w.r.t.\ the mesh nodes, which have been derived in the previous section, are used. In particular, the projection step consists of computing
\begin{equation}
	\label{eq:project_search_direction}
	P_h = \left( I - A\transposed \left( A A\transposed \right)^{-1} A \right) S_h.
\end{equation}
Note that $P_h$ can be interpreted as the DoF vector of some projected search direction $\mathcal{P}_h$, analogously to before. We note that the projection step is carried out by using \eqref{eq:lagrange_multiplier} and \eqref{eq:projected_direction}. In particular, the projected search direction $\mathcal{P}_h$ is feasible for the (linearized) constraints. 

To deal with the nonlinearity of mesh quality constraints, the line search procedure of Algorithm~\ref{algo:shape_descent} has to updated slightly in analogy to the strategy presented in Section~\ref{ssec:rosen_finite} to find a step which is feasible for all previously active constraints and which does not violate any new ones. The proposed steps of back-projecting the step and using a bisection method for these issues works just as detailed in Section~\ref{ssec:rosen_finite}. Note that after the modified line search procedure, we have obtained a new design which is represented by the updated computational mesh. The latter is, by construction, feasible for the mesh quality constraints. Hence, we can continue with Algorithm~\ref{algo:shape_descent} after a suitable step is found in the line search. In particular, this application of the gradient projection method guarantees that the mesh quality constraints formulated in Section~\ref{ssec:mesh_quality_constraints} cannot be violated so that the computational mesh is guaranteed to have a minimum, user-specified quality during the entire shape optimization. The minimum mesh quality only depends on the angle threshold $\alpha_\mathrm{thr}$ and its choice is discussed later in Section~\ref{ssec:choice_of_angle}.

\begin{Remark}
	
	Regarding the computational effort of our approach we note the following: The modifications introduced by the gradient projection method only affect the optimization variables and how they are treated. In particular, the method does not require additional solutions of the state or adjoint equations, which is the most expensive part of Algorithm~\ref{algo:shape_descent}. Additionally, as discussed in Remark~\ref{rem:number_of_active_constraints}, usually only a small number of constraints are actually active during the optimization so that the linear system solved in \eqref{eq:lagrange_multiplier} is comparatively small. For these reasons, our proposed approach is computationally inexpensive and does not require a lot of additional time or computational resources.
\end{Remark}

\subsubsection{Additional Constraints for Fixed Boundaries}

The framework for mesh quality constraints presented above works just fine in case that no boundaries are fixed. However, in the case that some boundaries are fixed for the optimization, the gradient deformation, or more generally, the search direction for gradient-based shape optimization, usually take these into account via homogeneous Dirichlet boundary conditions for the gradient deformation as discussed in Definition~\ref{def:gradient_deformation}. Suppose that some mesh quality constraint is active in a cell, where one face (or side) of the cell is part of the fixed boundary. To compute a direction which is feasible w.r.t.\ the linearized constraints, the search direction is projected as discussed above. However, this could yield a projected search direction which no longer vanishes on the fixed boundary. To avoid this, we have to impose additional constraints for these fixed boundaries: Each node, which is part of a fixed boundary, must not be moved. Suppose that some node $v_m = [x_m, y_m]\transposed$ (or $v_m = [x_m, y_m, z_m]$ for three-dimensional problems) is part of the fixed boundary. Then, the desired behavior is achieved with the following equality constraint
\begin{equation}
	\label{eq:fixed_boundary_constraint}
	g_m = v_m - v_m^0 = 0,
\end{equation}
where $v_m^0$ denotes the location of node $v_m$ in the initial mesh. Its derivative is given by
\begin{equation*}
	\frac{\partial g_m}{\partial v_m} = \begin{cases}
		[1,1]\transposed \quad &\text{for two-dimensional problems},\\
		[1,1,1]\transposed \quad &\text{for three-dimensional problems,}
	\end{cases}
\end{equation*}
where we use the notation \eqref{eq:derivative_notation_2D} and \eqref{eq:derivative_notation_3D}, respectively, for the derivatives. Similarly to before, the derivative of $g_m$ w.r.t.\ any other node is zero. Hence, the derivative of the equality constraint \eqref{eq:fixed_boundary_constraint} is sparse and has only two nonzero entries in 2D and three nonzero entries in 3D.

\subsubsection{Choice of the Minimum Angle Threshold}
\label{ssec:choice_of_angle}

Let us now investigate the choice of the minimum angle threshold. To do so, we note that the gradient projection method presented earlier is a feasible method in the sense that the initial guess has to be feasible and all iterates it produces are feasible w.r.t.\ the constraints. In particular, the initial mesh has to satisfy the posed mesh quality constraints. This limits the possible choices of the minimum angle threshold $\alpha_\mathrm{thr}$. 

We consider two approaches for the choice of $\alpha_\mathrm{thr}$. The first is to specify a single minimum angle threshold $\alpha_\mathrm{thr}$ which should be used for all mesh cells. This works well if the initial mesh is, e.g., very regular and contains cells with approximately the same quality. If we denote by $\alpha^0_\mathrm{min}$ the minimum (solid) angle in the initial mesh, then $\alpha_\mathrm{thr}$ has to satisfy $\alpha_\mathrm{thr} \leq \alpha^0_\mathrm{min}$ so that the initial mesh is feasible w.r.t.\ the mesh quality constraints.

However, for practical applications there might be some areas where the mesh generation is problematic and it might not be possible to obtain a good quality initial mesh. Another possible issue is the use of mesh cells which are stretched to resolve boundary layers, in particular, if triangular or tetrahedral elements are used. These scenarios have in common that there are usually just a few bad cells in the initial mesh, either in problematic areas of the geometry or at regions, where a boundary layer is resolved, but the majority of mesh cells is regular and of good quality. Using a single global minimum angle threshold does not work well in these cases: For the few bad cells, the lower bound posed this way is sensible. For the remaining mesh cells with good quality, however, this lower bound is too small so that it allows for a deterioration of the mesh quality. To circumvent this issue, a second approach for choosing $\alpha_\mathrm{thr}$ is to use a different threshold for each cell, based on the initial angles of the repsective cell. To do so, we consider some mesh cell $\tau$ (either a triangle or tetrahedron) on the initial mesh and denote its minimum (solid) angle by $\alpha^0_{\mathrm{min}, \tau}$. Moreover, let $\nu \in (0,1)$ be some relative tolerance. Then, we propose to use the minimum angle threshold 
\begin{equation*}
	\alpha_\mathrm{thr} = \nu \alpha^0_{\mathrm{min}, \tau}
\end{equation*}
for mesh cell $\tau$. This choice of $\alpha_\mathrm{thr}$ is relative to the quality of the cells in the initial mesh and, thus, well-suited for the case of mesh with very inhomogeneous quality.

Finally, we note that these choices can also be combined as follows. First, a global minimum angle threshold $\alpha_\mathrm{thr}$ and a relative tolerance $\nu$ is chosen as discussed above. For all mesh cells with $\alpha^0_{\mathrm{min}, \tau} > \alpha_\mathrm{thr}$, the threshold $\alpha_\mathrm{thr}$ is used. For all other cells, the minimum angle threshold is given by $\nu \alpha^0_{\mathrm{min}, \tau}$. This procedure ensures that the initial mesh is feasible w.r.t.\ the mesh quality constraints and combines the advantages of both approaches.

\subsubsection{Numerical Implementation}

We have implemented the proposed gradient projection method and the mesh quality constraints as extension to our open-source software package cashocs \cite{Blauth2021cashocs, Blauth2023Version}, which is based on the finite element software FEniCS \cite{Alnaes2015FEniCS,Logg2012Automated}. Our software package cashocs implements and automates the adjoint approach for solving PDE constrained optimization problems with a particular focus on shape optimization problems. To do so, cashocs utilizes the automatic differentiation capabilities of the so-called unified form language (UFL) \cite{Alnaes2014Unified, Ham2019Automated}, which is a part of FEniCS. We note that the required adjoint systems and shape derivatives for solving the optimization problems are computed automatically by cashocs and do not have to be derived or implemented by the user. For this reason, we also do not state these for the subsequent numerical examples but only refer to the relevant literature.

To facilitate a fast numerical evaluation of the constraint functions and their derivatives, these are implemented as custom C++ code, which is compiled at run-time with FEniCS. The linear algebra sub-problems of projecting a direction to the tangent space of the currently active constraints as well as the back-projection discussed in Section~\ref{ssec:rosen_finite} are solved with PETSc \cite{Balay2025PETSc/TAO, Balay2025PETSc, Balay1997Efficient} using its Python wrapper petsc4py \cite{Dalcin2011Parallel}. As linear solver a conjugate gradient (CG) method preconditioned by BoomerAMG \cite{Henson2002BoomerAMG}, an algebraic multigrid method that is part of hypre \cite{hypre} and interfaced through PETSc, is used. Moreover, our method is parallelized with MPI to facilitate the solution of large-scale problems arising from industrial applications.

\section{Numerical Experiments}
\label{sec:numerics}

Let us now consider our proposed method numerically. First, we investigate two academic model examples that highlight our approach's superior performance w.r.t.\ the preservation of mesh quality compared to the classical approach. To show that our proposed method works well with both the BFGS and gradient descent methods, the first problem is solved with the former while the second is solved with the latter method. Afterwards, we consider a large-scale three-dimensional shape optimization problem for increasing the separation efficiency of a structured packing in a distillation column, which demonstrates the capabilities of our proposed method for industrial applications.

\subsection{Shape Optimization of an Obstacle in Stokes Flow}

\begin{figure}[b]
	\centering
	\begin{subfigure}{0.6\textwidth}
		\centering
		\includegraphics[width=\textwidth]{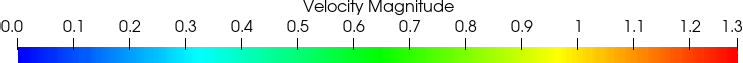}
	\end{subfigure}\\
	\vspace{0.1cm}
	\begin{subfigure}{0.475\textwidth}
		\centering
		\includegraphics[width=\textwidth]{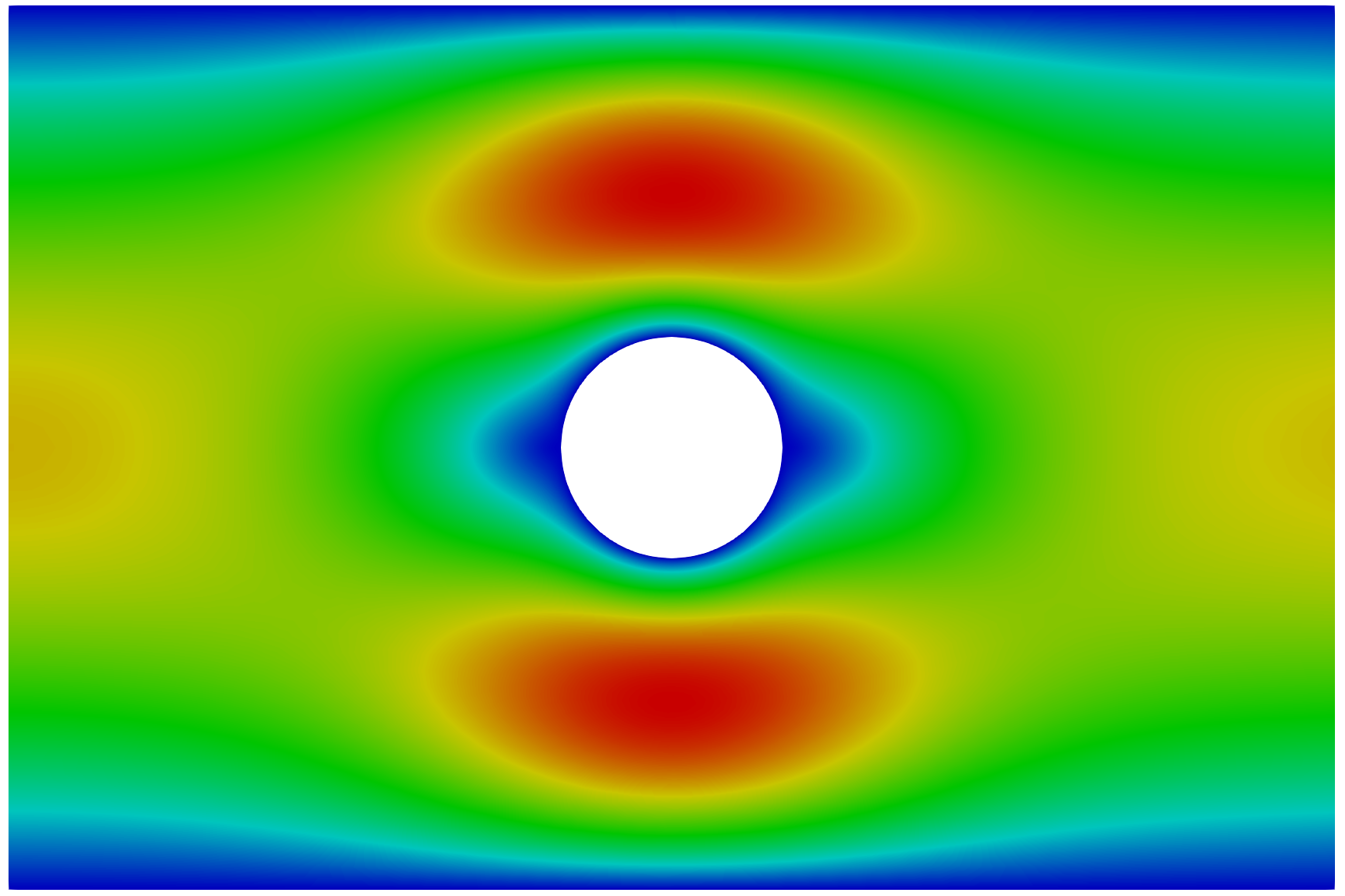}
		\caption{Initial geometry.}
	\end{subfigure}
	\hfil
	\begin{subfigure}{0.475\textwidth}
		\centering
		\includegraphics[width=\textwidth]{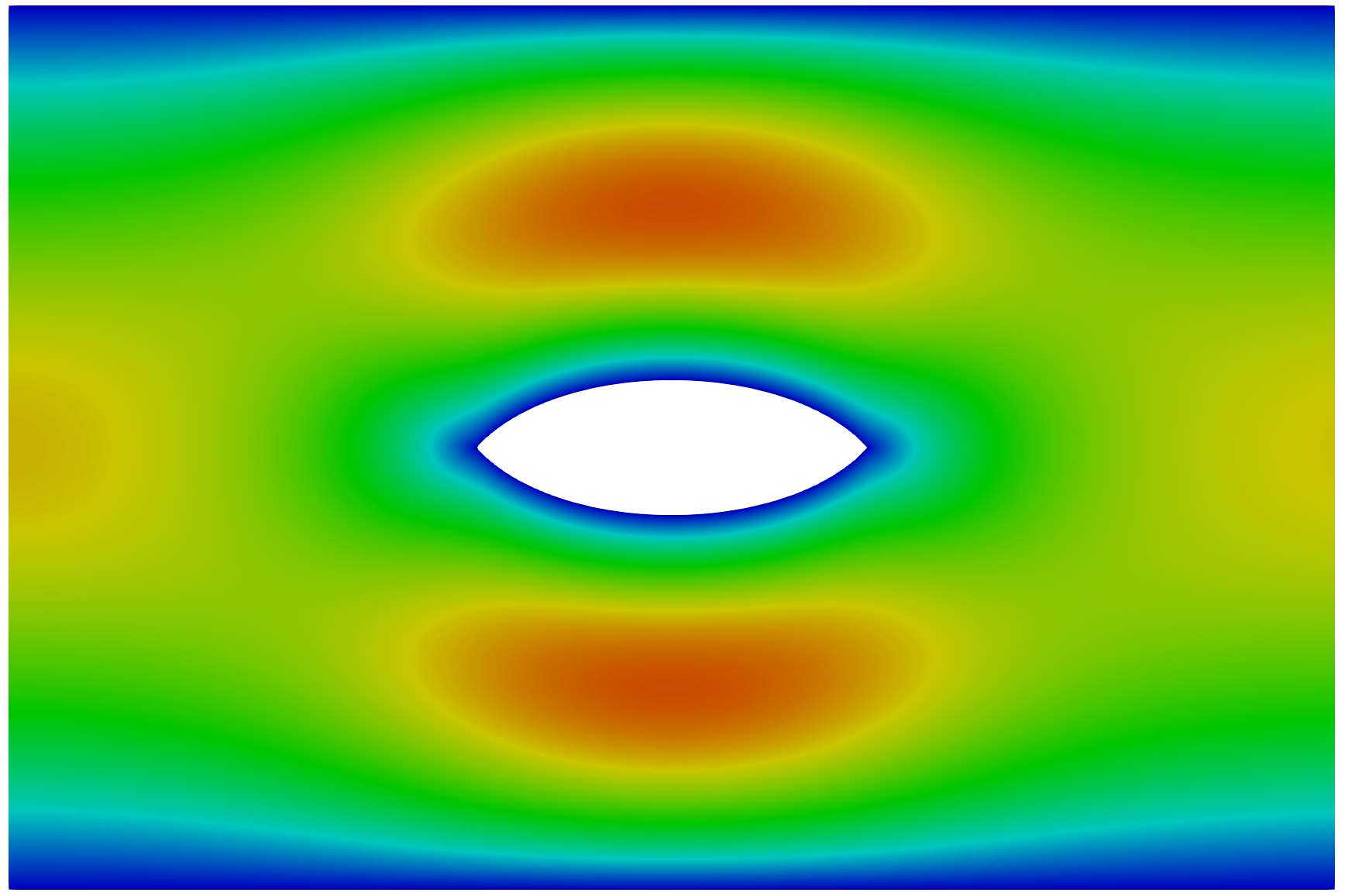}
		\caption{Optimized geometry.}
	\end{subfigure}
	\caption{Velocity magnitude for the Stokes problem \eqref{eq:stokes}, computed with the BFGS method with mesh quality constraints.}
	\label{fig:comparison_stokes_velocity}
\end{figure}

Let us consider the shape optimization of an obstacle in Stokes flow, which is a popular benchmark problem for shape optimization algorithms and has been used, e.g., in \cite{Blauth2021Nonlinear,Iglesias2018Two, Mueller2021novel}. For this, we consider the minimization of the energy dissipated by the flow. Additionally, the obstacle's volume and barycenter are fixed to ensure a non-trivial solution to the problem.
Let the flow domain be denoted by $\Omega$ and its boundary by $\Gamma = \partial \Omega$. The exterior boundary is divided into the flow inlet $\Gamma_\mathrm{in}$, the wall boundary $\Gamma_\mathrm{wall}$, and the outlet $\Gamma_\mathrm{out}$. Its interior boundary is the boundary of the obstacle $\Gamma_\mathrm{obs}$. The non-dimensionalized Stokes system is given by
\begin{equation}
	\label{eq:stokes}
	\begin{alignedat}{2}
		-\Delta u + \nabla p &= 0 \quad &&\text{ in } \Omega,\\
		\nabla \cdot u &= 0 \quad &&\text{ in } \Omega,\\
		u &= u_\mathrm{in} \quad &&\text{ on } \Gamma_\mathrm{in},\\
		u &= 0 \quad &&\text{ on } \Gamma_\mathrm{wall} \cup \Gamma_\mathrm{obs},\\
		\partial_n u - pn &= 0 \quad &&\text{ on } \Gamma_\mathrm{out},
	\end{alignedat}
\end{equation}
where $u$ and $p$ are the fluid velocity and pressure, respectively. Note that we use the usual do-nothing condition (see, e.g., \cite{John2016Finite}) at the outlet $\Gamma_\mathrm{out}$, where $n$ denotes the outer unit normal vector on $\Gamma$ and $\partial_n u = Du\ n$ is the normal derivative of the velocity. For the geometrical constraints, we fix the volume and barycenter of the flow domain, which is equivalent to fixing those of the obstacle. The volume and barycenter of $\Omega$ are given by
\begin{equation}
	\label{eq:notation_volume}
	\mathrm{vol}(\Omega) = \integral{\Omega} 1 \dmeas{x} \qquad \text{ and } \qquad \mathrm{bc}(\Omega) = \frac{1}{\mathrm{vol}(\Omega)} \integral{\Omega} x \dmeas{x}.
\end{equation}
Hence, the shape optimization problem can be formulated as
\begin{equation}
	\label{eq:shape_stokes}
	\min_{\Omega \in \admissiblegeom} \mathcal{J}(\Omega, u) = \integral{\Omega} Du : Du \dmeas{x} + \frac{\nu_\mathrm{vol}}{2} \left( \mathrm{vol}(\Omega) - \mathrm{vol}(\Omega^0) \right)^2 + \frac{\nu_\mathrm{bc}}{2} \abs{ \mathrm{bc}(\Omega) - \mathrm{bc}(\Omega^0) }^2 \quad \text{ s.t. } \eqref{eq:stokes},
\end{equation}
where $\Omega^0$ denotes the initial geometry. Note that the geometrical constraints are regularized with a quadratic penalty function as in \cite{Blauth2021Nonlinear}. The set of admissible geometries for this problem is given by
\begin{equation*}
	\mathcal{A} = \Set{\Omega \subset \R^d | \Omega \subset \holdall, \Gamma_\mathrm{in} = \Gamma_\mathrm{in}^0, \Gamma_\mathrm{wall} = \Gamma_\mathrm{wall}^0, \Gamma_\mathrm{out} = \Gamma_\mathrm{out}^0}
\end{equation*}
for some initial flow domain $\Omega^0$ with boundary $\Gamma^0$, i.e., only the boundary $\Gamma_\mathrm{obs}$ is deformable and subjected to the optimization. The hold-all domain $\holdall$ is given by $\Omega \cup \Omega_\mathrm{obs}$, where $\Omega_\mathrm{obs}$ is the domain of the obstacle. For the shape derivative and adjoint system of this problem, we refer the reader to \cite{Blauth2021Adjoint}. Note that the solution of the Stokes system \eqref{eq:stokes} is visualized in Figure~\ref{fig:comparison_stokes_velocity} where the velocity magnitude is depicted on the initial and optimized geometry.

\begin{figure}[b]
	\centering
	\begin{subfigure}{0.475\textwidth}
		\centering
		\includegraphics[width=\textwidth]{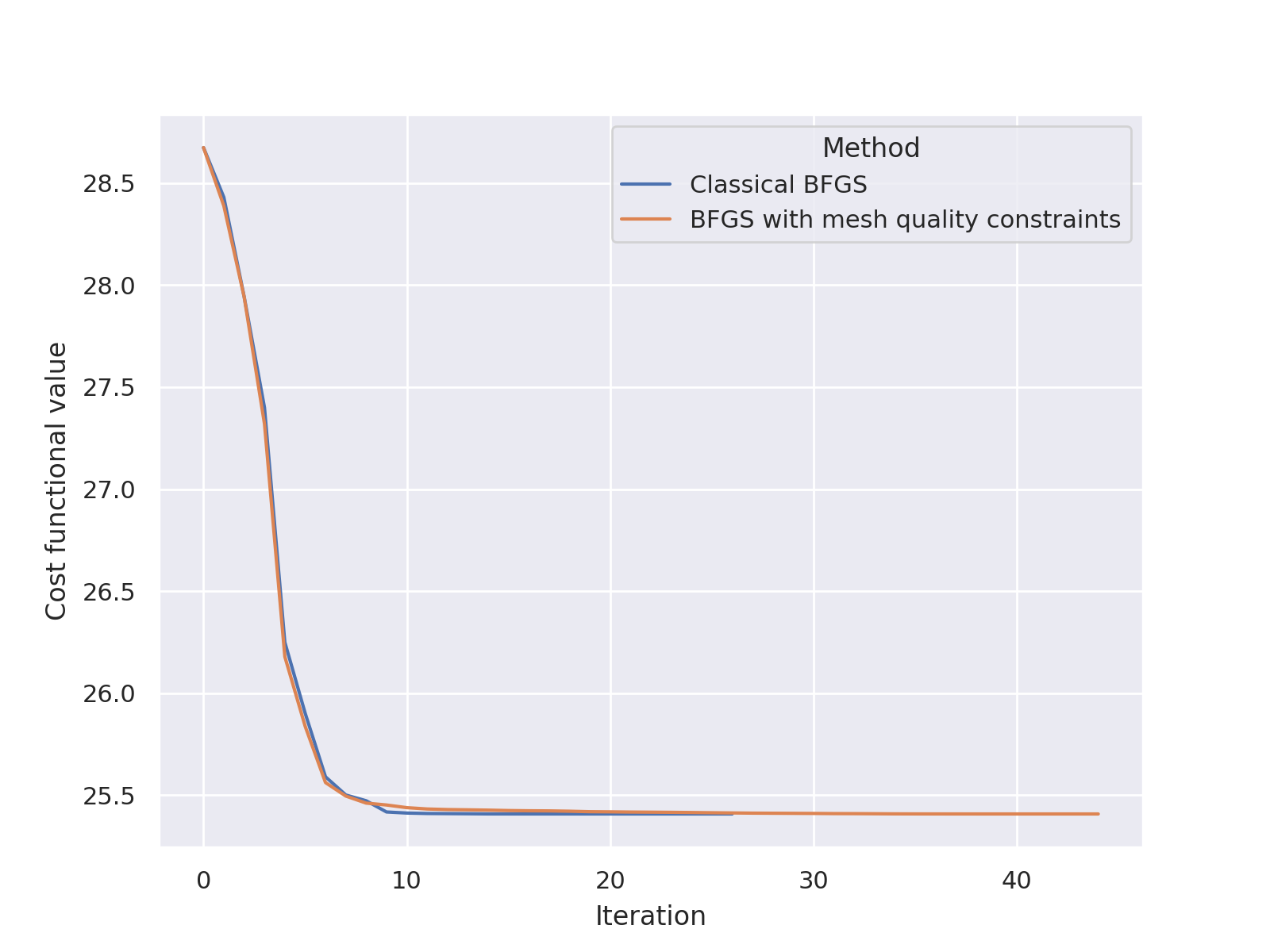}
		\caption{Cost functional value.}
		\label{sfig:stokes_cost_functional}
	\end{subfigure}
	\hfil
	\begin{subfigure}{0.475\textwidth}
		\centering
		\includegraphics[width=\textwidth]{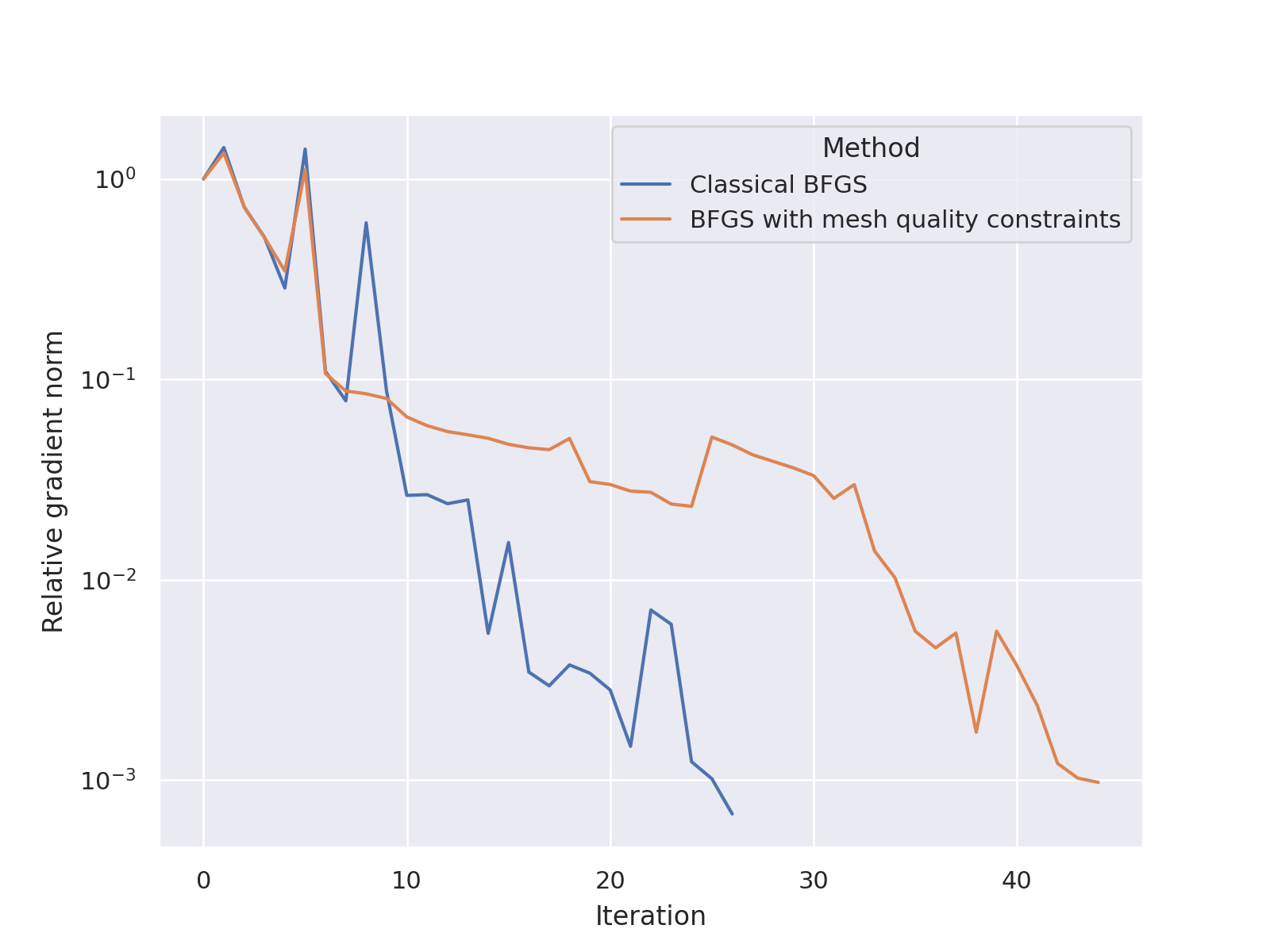}
		\caption{Relative norm of the gradient deformation.}
	\end{subfigure}
	\caption{History of the optimization algorithm for problem \eqref{eq:shape_stokes}.}
	\label{fig:history_stokes}
\end{figure}

We consider this problem in two dimensions and use $\holdall = (-3, 3) \times (-2, 2)$, $\Omega^0 = \holdall \setminus \Omega_\mathrm{obs}^0$, where the initial obstacle $\Omega_\mathrm{obs}^0$ is a circle with center $(0,0)$ and radius 0.5, in analogy to \cite{Blauth2021Nonlinear}. The domain is discretized with a non-uniform mesh generated with Gmsh \cite{Geuzaine2009Gmsh}. The mesh consists of \num{3535} nodes and \num{6674} triangles and is refined near $\Gamma_\mathrm{obs}$. As we employ our optimization software cashocs \cite{Blauth2021cashocs,Blauth2023Version}, we use the finite element software FEniCS \cite{Alnaes2015FEniCS, Logg2012Automated} for the discretization of the involved PDEs. In particular, for the discretization of the state and adjoint systems the LBB-stable Taylor-Hood elements are used, i.e., quadratic Lagrange elements are used for the velocity and linear Lagrange elements are used for the pressure. The inlet velocity is given by the parabolic profile $u_\mathrm{in}(x) = \nicefrac{1}{4} \left( 2-x_2 \right) ( 2 + x_2 )$. Finally, the parameters for the quadratic penalty functions are chosen as $\nu_\mathrm{vol} = \num{1e3}$ and $\nu_\mathrm{bc} = \num{1e5}$, which ensures that the geometrical constraints are satisfied with an accuracy below \qty{0.1}{\percent}. We solve this problem with the BFGS method implemented in cashocs (cf.~\cite{Schulz2016Efficient}) with and without the proposed mesh quality constraints. The linear elasticity equations \eqref{eq:linear_elasticity} are used for computing the gradient deformation with parameters $\mu_\mathrm{elas} = 1$, $\lambda_\mathrm{elas} = 0$, and $\delta_\mathrm{elas} = 0$. The optimization is considered to be converged if the norm of the gradient deformation is below a relative tolerance of \num{1e-3}. For the mesh quality constraints, we choose a minimum angle threshold of $\alpha_\mathrm{thr} = \qty{0.436}{\radian}$, which corresponds to a minimum angle of \qty{25}{\degree}, and a numerical tolerance of \num{1e-2}.

\begin{figure}[t]
	\centering
	\begin{subfigure}{0.6\textwidth}
		\centering
		\includegraphics[width=\textwidth]{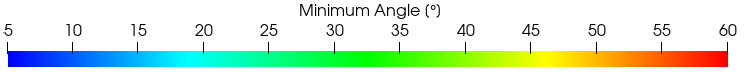}
	\end{subfigure}\\
	\begin{subfigure}[t]{0.475\textwidth}
		\centering
		\includegraphics[width=\textwidth]{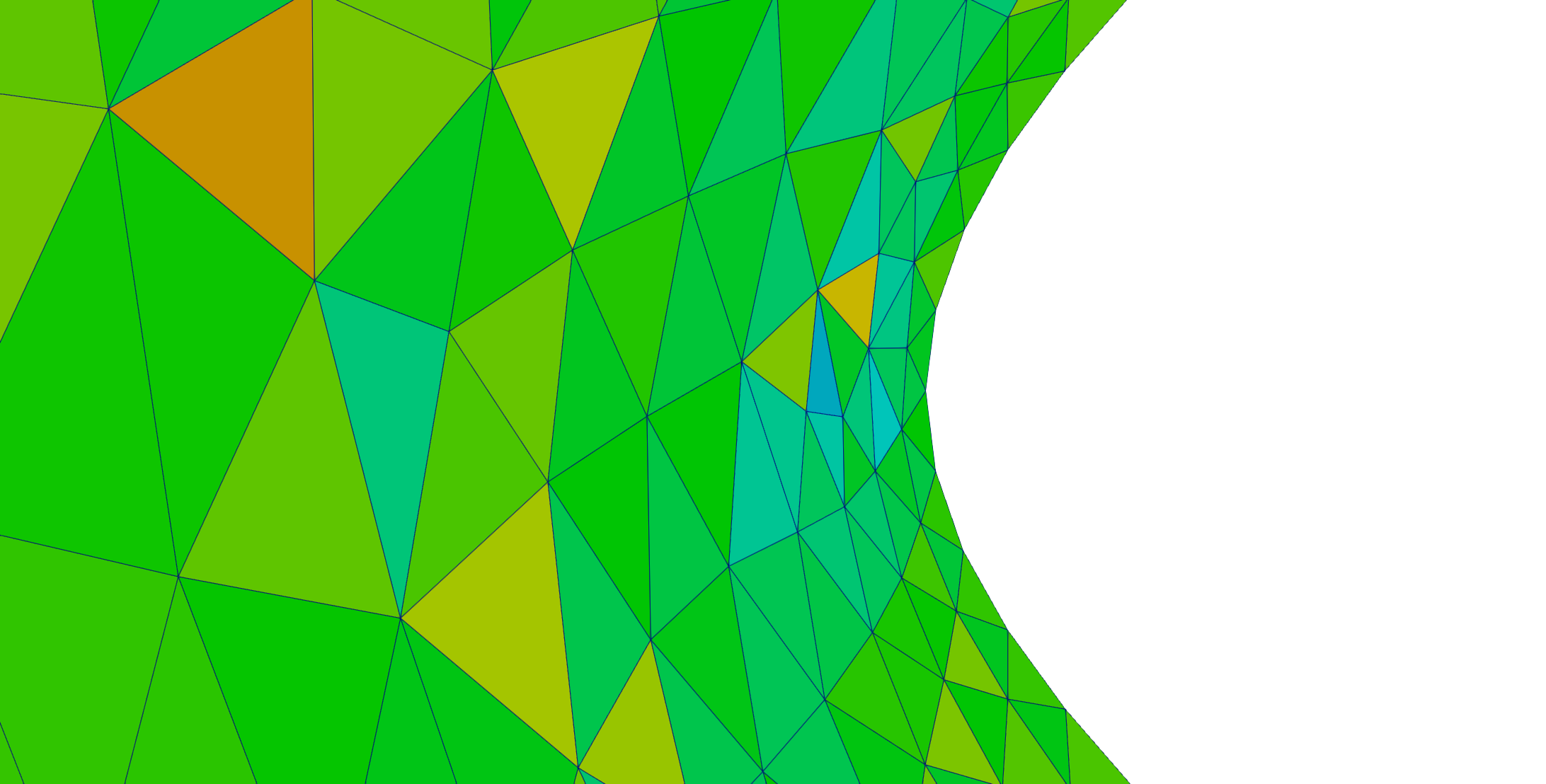}
		\caption{Classical BFGS method -- after 10 iterations.}
		\label{sfig:stokes_classical_10}
	\end{subfigure}
	\hfill
	\begin{subfigure}[t]{0.475\textwidth}
		\centering
		\includegraphics[width=\textwidth]{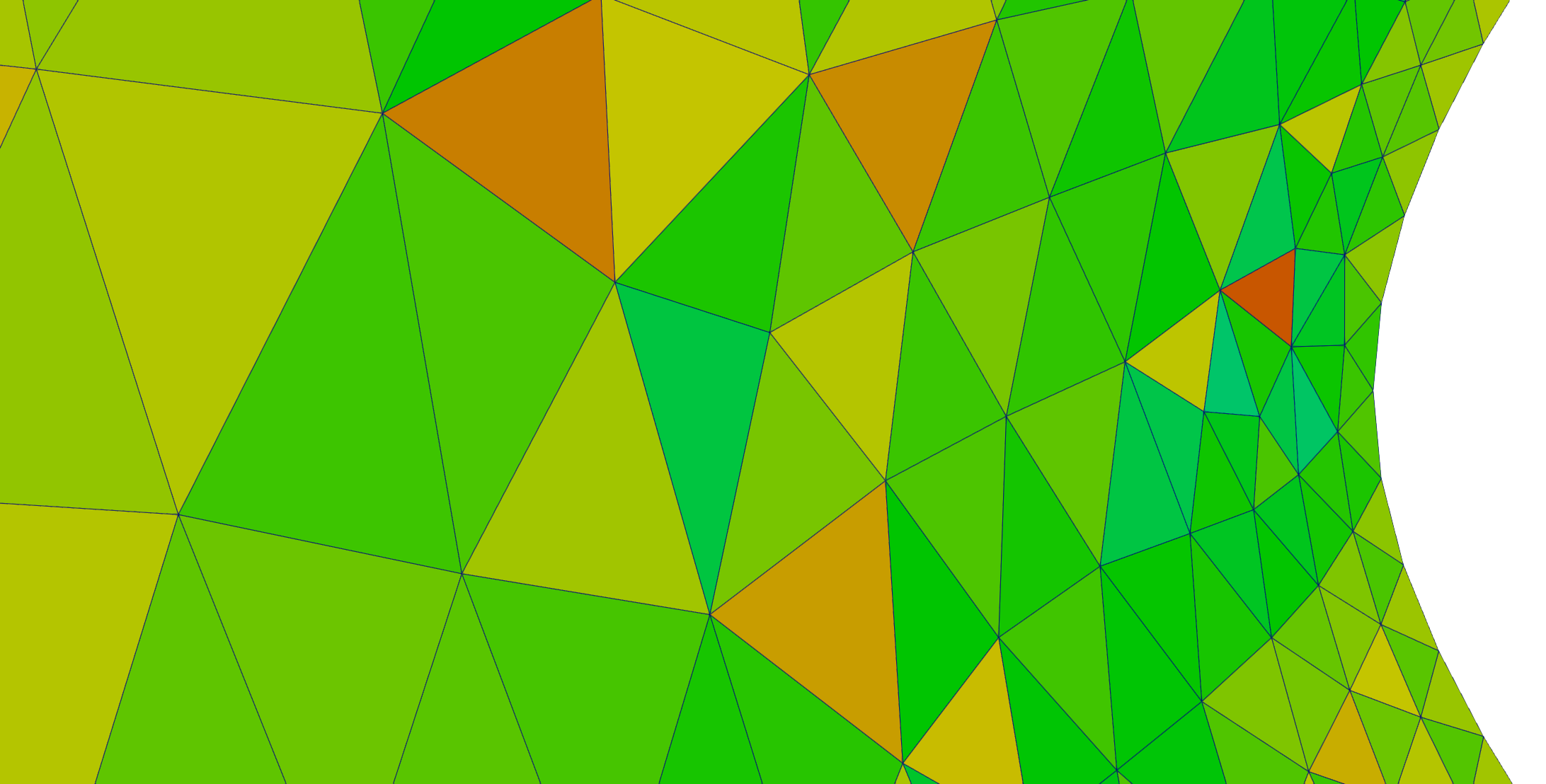}
		\caption{BFGS method with mesh quality constraints -- after 10 iterations.}
		\label{sfig:stokes_mqc_10}
	\end{subfigure}\\
	\begin{subfigure}[t]{0.475\textwidth}
		\centering
		\includegraphics[width=\textwidth]{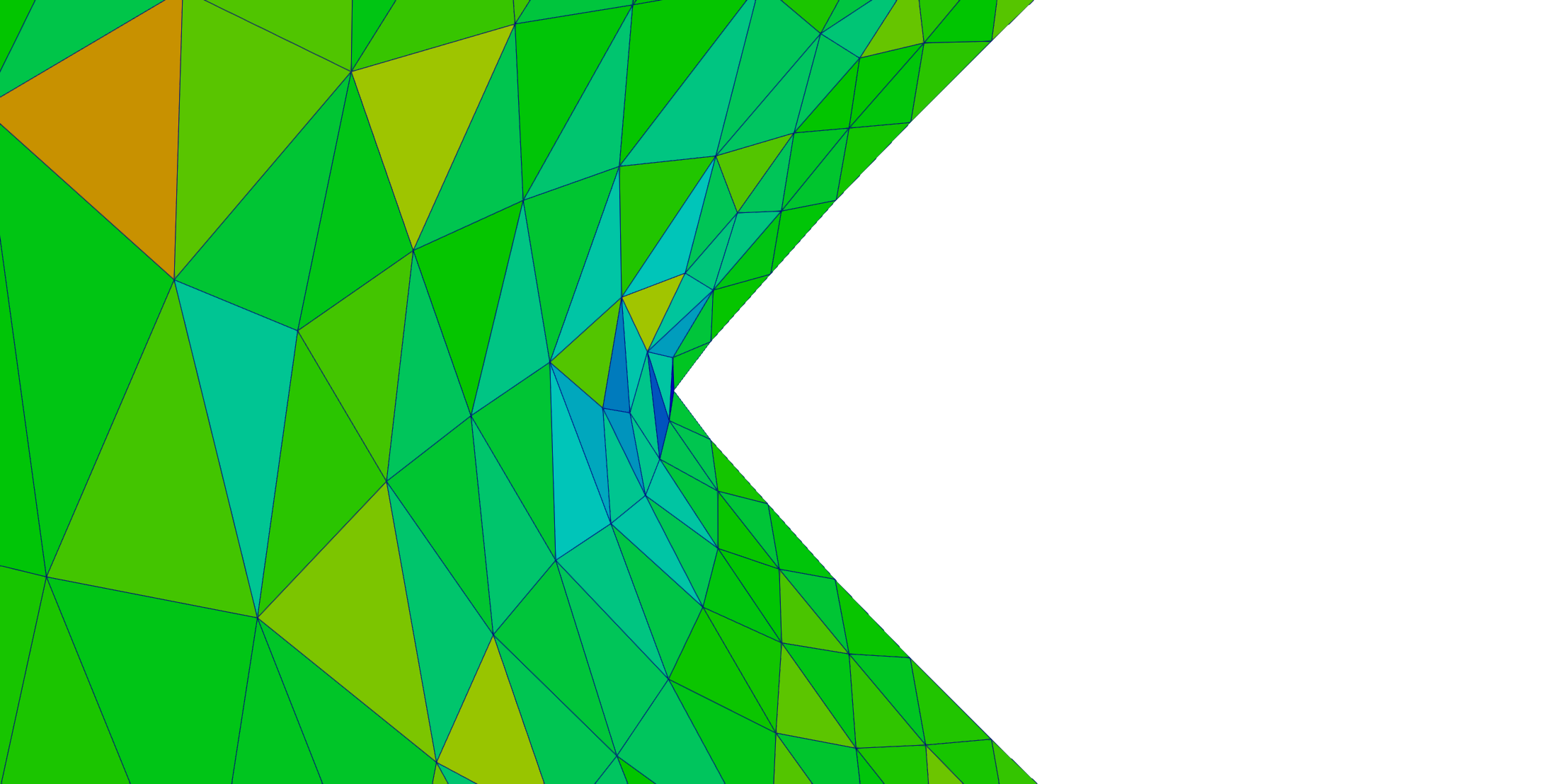}
		\caption{Classical BFGS method -- optimized shape.}
		\label{sfig:stokes_classical}
	\end{subfigure}
	\hfil
	\begin{subfigure}[t]{0.475\textwidth}
		\centering
		\includegraphics[width=\textwidth]{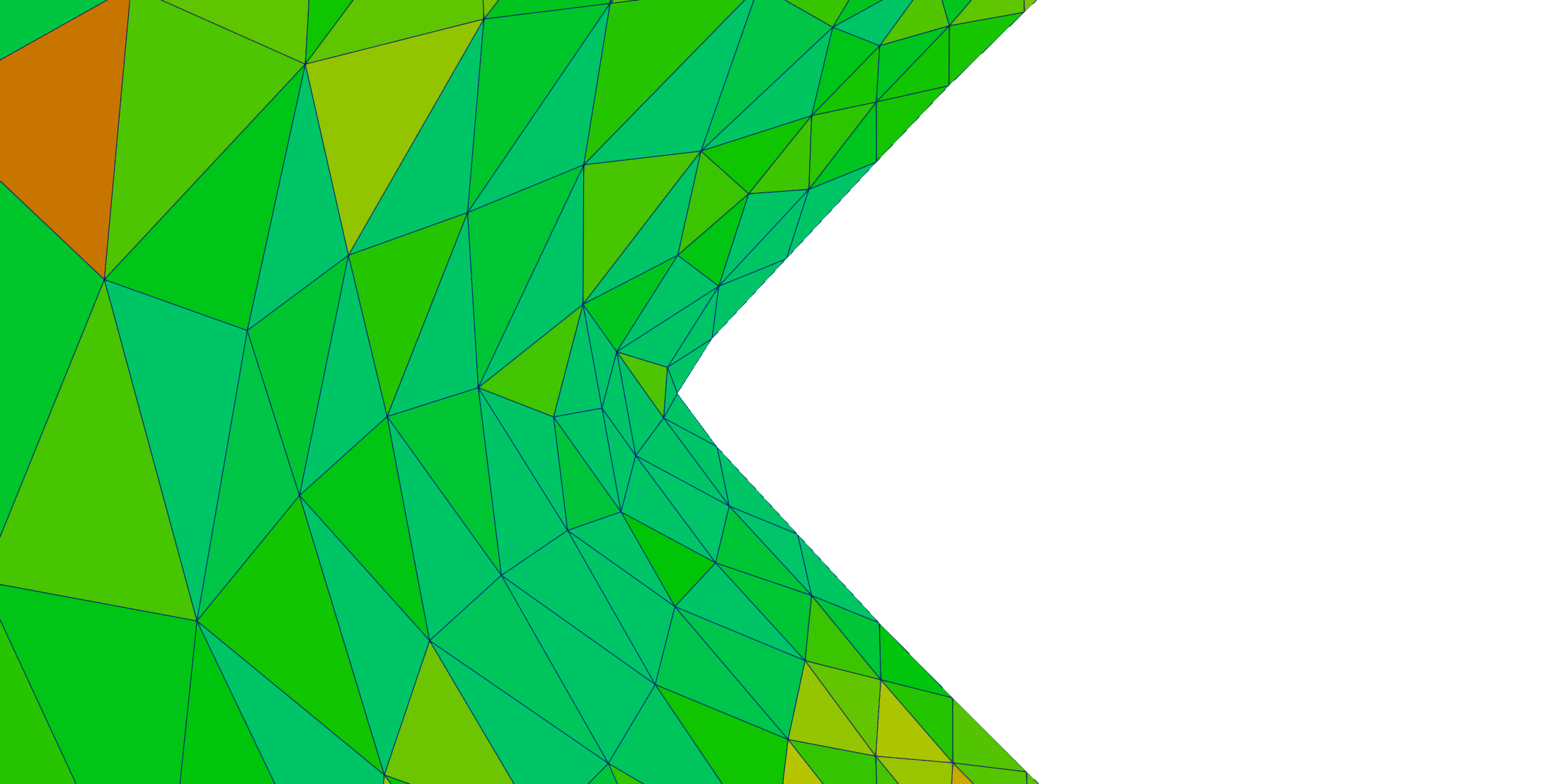}
		\caption{BFGS method with mesh quality constraints -- optimized shape.}
		\label{sfig:stokes_mqc}
	\end{subfigure}
	\caption{Optimized geometry for problem \eqref{eq:shape_stokes} -- close up of the front tip.}
	\label{fig:quality_stokes}
\end{figure}

\begin{figure}[b]
	\centering
	\begin{subfigure}[t]{0.475\textwidth}
		\centering
		\includegraphics[width=\textwidth]{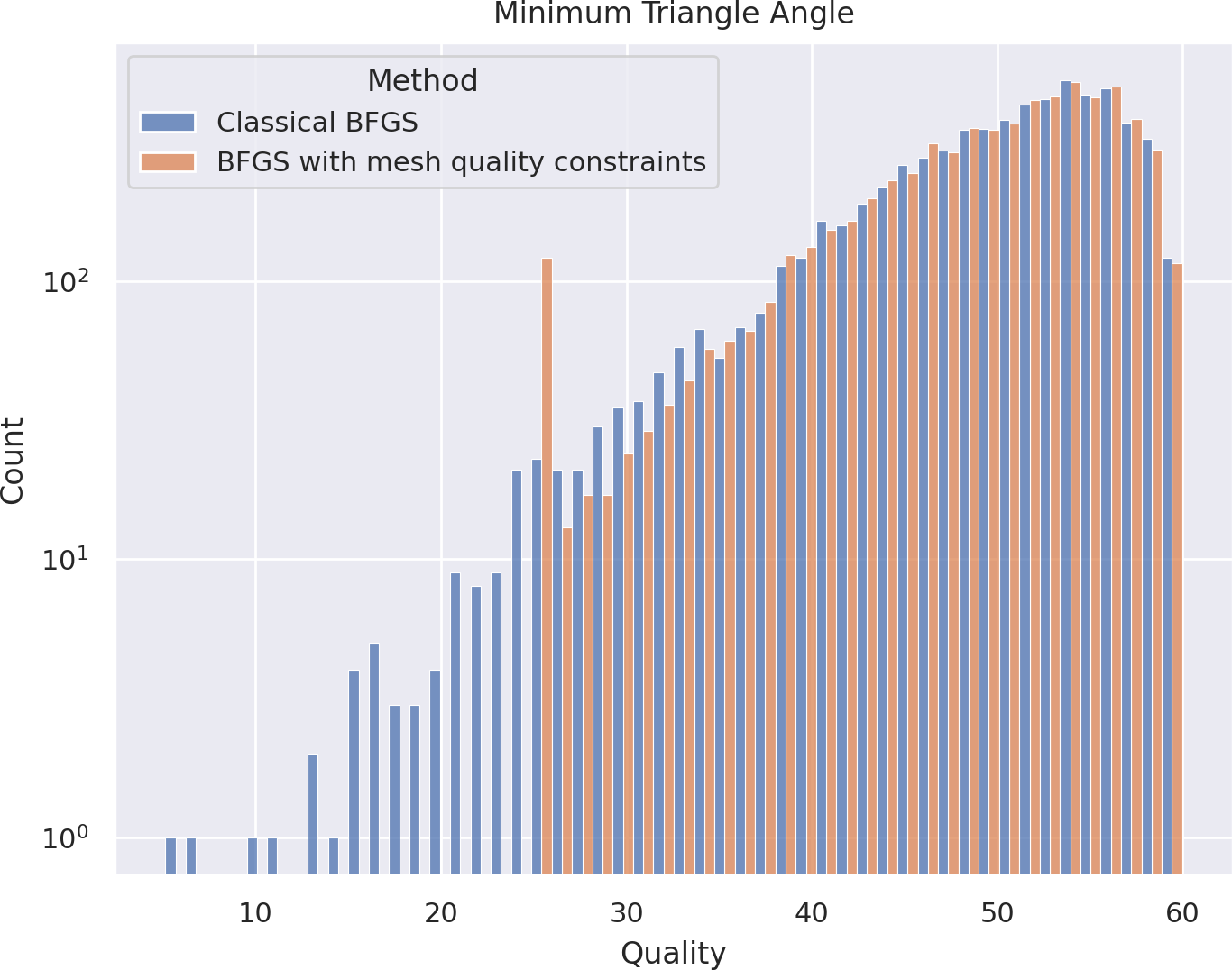}
		\caption{Minimum triangle angle in degrees.}
		\label{fig:dist_stokes_angle}
	\end{subfigure}
	\hfil
	\begin{subfigure}[t]{0.475\textwidth}
		\centering
		\includegraphics[width=\textwidth]{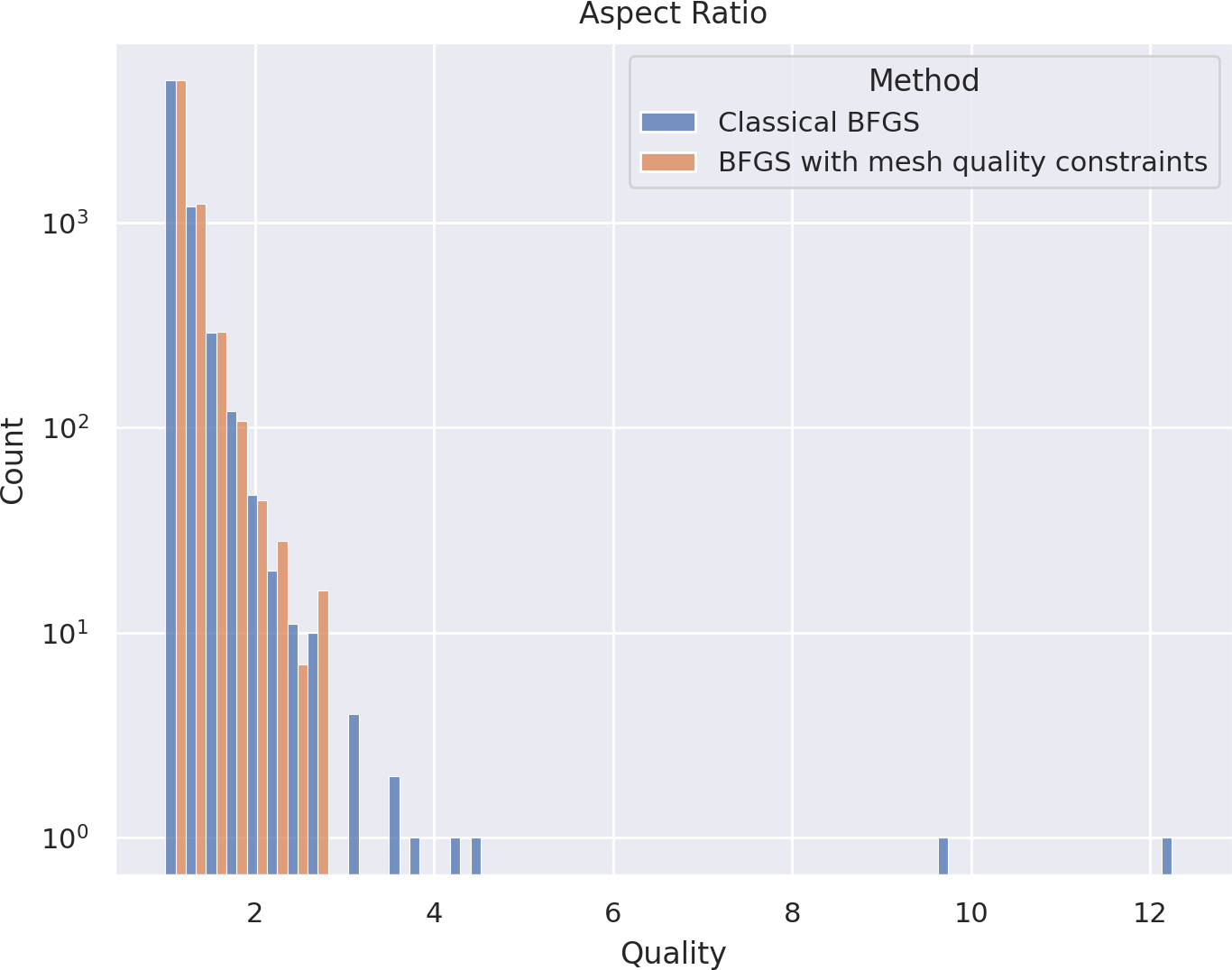}
		\caption{Aspect ratio.}
		\label{fig:dist_stokes_aspect_ratio}
	\end{subfigure}
	\caption{Histograms showing the distribution of the mesh quality after the shape optimization for problem \eqref{eq:shape_stokes}.}
	\label{fig:dist_stokes}
\end{figure}

The history of the cost functional and the gradient norm over the course of the optimization is visualized in Figure~\ref{fig:history_stokes}. We see that the BFGS method converged in \num{26}~iterations for the classical approach without mesh quality constraints and that it needed \num{44} iterations to converge with the approach proposed in this paper. Overall, the cost functional value for both approaches is nearly indistinguishable, but we observe that the norm of the gradient deformation reduces faster for the classical approach.
The resulting meshes together with the minimum angle of the mesh cells are shown in Figure~\ref{fig:quality_stokes}, where a close up around the front tip of the obstacle is shown. Note that Figures~\ref{sfig:stokes_classical_10} and \ref{sfig:stokes_mqc_10} show the geometry after 10 iterations, whereas Figures~\ref{sfig:stokes_classical} and \ref{sfig:stokes_mqc} show the optimized geometries obtained with the classical BFGS and BFGS method with mesh quality constraints, respectively.
In particular, we observe that the classical approach without mesh quality constraints approximates the optimal shape very well, however, doing so results in several mesh cells with very acute angles and, thus, a very low quality of the mesh near the tip. On the other hand, our proposed approach also approximates the optimal shape very well, but without any deterioration of the mesh quality. In particular, it can be seen that there are no bad cells in the optimized mesh close to the tip for our approach. 

\begin{Remark}
	We note that when considering only the evolution of the cost functional in Figure~\ref{sfig:stokes_cost_functional}, it appears that both approaches do not make much further progress after about 10 iterations. However, in Figures~\ref{sfig:stokes_classical_10} and \ref{sfig:stokes_mqc_10} the meshes for both approaches are shown after 10 iterations. Comparing these meshes with the optimized ones shown in Figures~\ref{sfig:stokes_classical} and \ref{sfig:stokes_mqc}, we observe that the shape has not yet converged and that the front tip is still very rounded for both approaches. In addition, we observe that the geometry produced with the classical BFGS method is a bit closer to the optimized one compared to the one produced by the BFGS method with mesh quality constraints. For these reasons, we observe that the norm of the gradient deformation is, indeed, the appropriate stopping criterion.
\end{Remark}

\begin{figure}[t]
	\centering
	\begin{subfigure}[t]{0.475\textwidth}
		\centering
		\includegraphics[width=\textwidth]{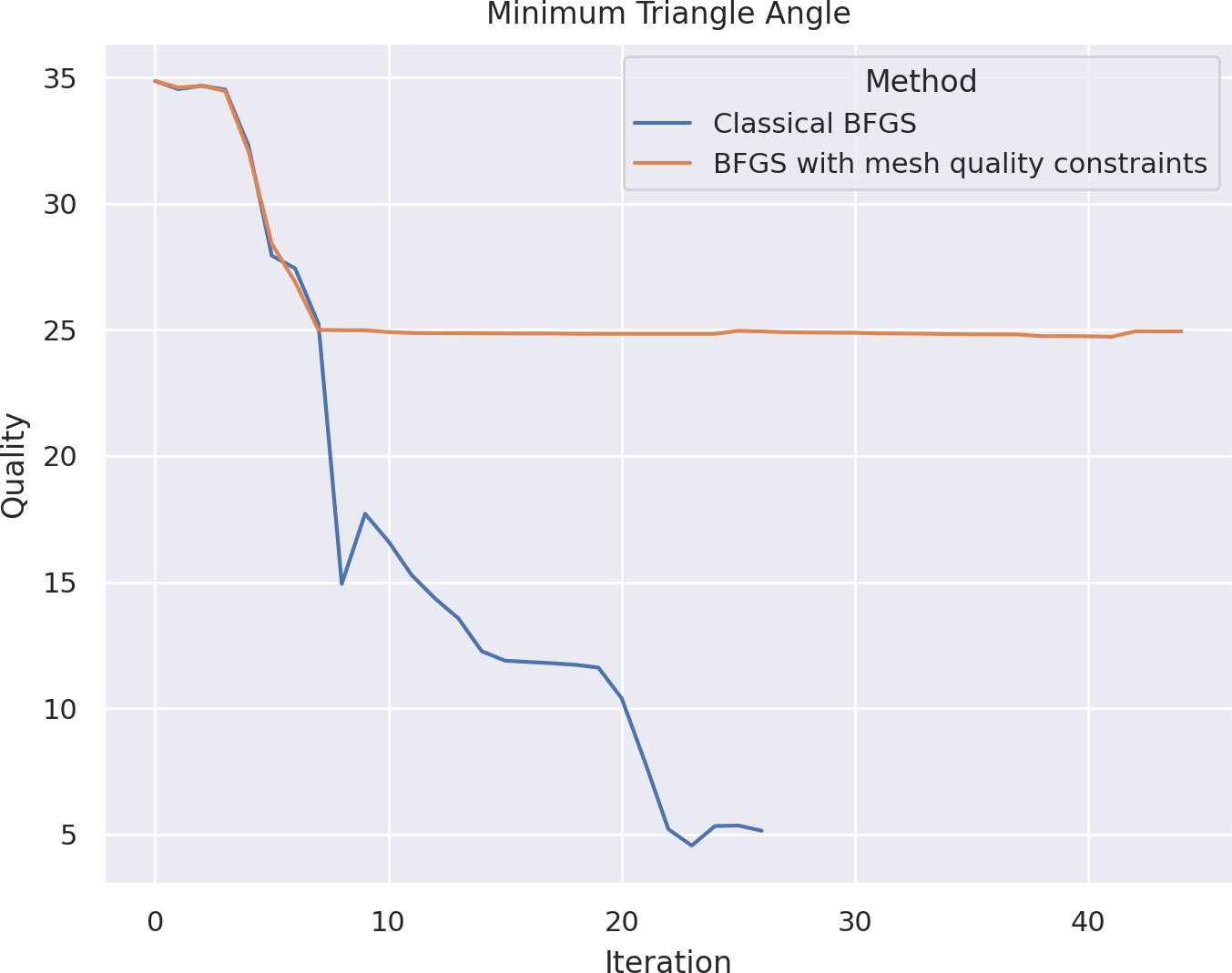}
		\caption{Minimum triangle angle in degrees.}
		\label{fig:iterations_stokes_angle}
	\end{subfigure}
	\hfil
	\begin{subfigure}[t]{0.475\textwidth}
		\centering
		\includegraphics[width=\textwidth]{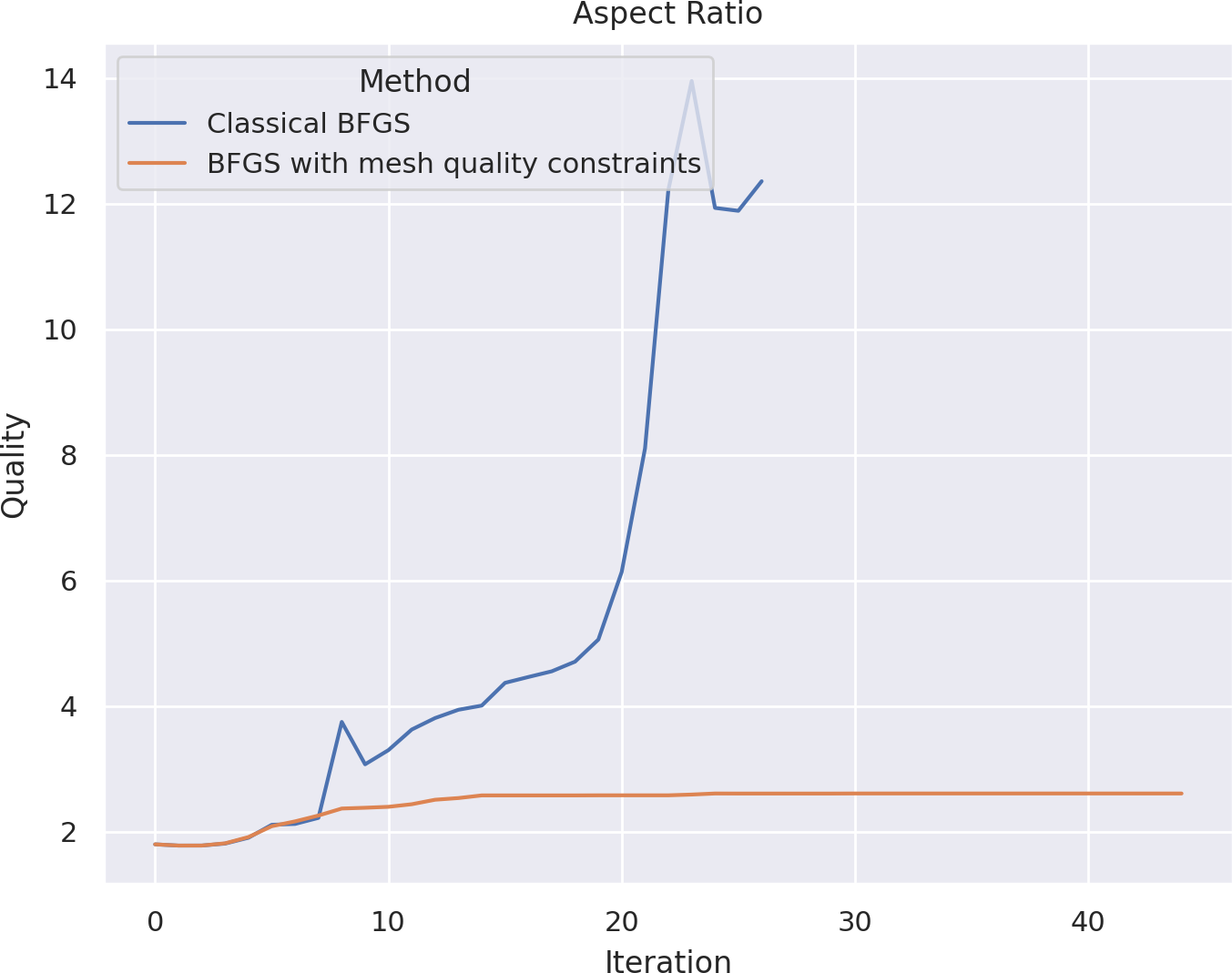}
		\caption{Aspect ratio.}
		\label{fig:iterations_stokes_aspect_ratio}
	\end{subfigure}
	\caption{Evolution of the mesh quality during the shape optimization for problem \eqref{eq:shape_stokes}.}
	\label{fig:iterations_stokes}
\end{figure}

The observation that our proposed method yields an overall improved mesh quality is reinforced in Figure~\ref{fig:dist_stokes}, where a histogram of the triangle angles and the aspect ratio of the optimized mesh with and without mesh quality constraints is shown. Here, the aspect ratio is is defined as $\frac{\abs{t}_\infty}{2 \sqrt{3} r}$, where $\abs{t}_\infty$ is the longest edge of the triangle and $r$ denotes its inradius. Note that the aspect ratio is normalized so that it has a value of \num{1} for a regular triangle and tends to $+\infty$ for a degenerating cell. From this figure, we can clearly see the superior behavior of our proposed method. In particular, the minimum angle is bounded from below by $\alpha_\mathrm{thr} = \qty{25}{\degree}$ for our proposed method, whereas it is significantly lower for the classical approach, where several badly shaped mesh cells exist. This is similar for the aspect ratio, where our proposed method ensures its boundedness, so that our method leads to higher quality meshes overall. Due to the logarithmic scaling of the histograms' $y$-axis, we can also observe a typical feature of the mesh quality in the context of shape optimization: Usually, there are only a couple of mesh cells that are deformed in a bad way and whose quality deteriorates over the course of the optimization. In this example, there are only about \num{10} cells which have a higher aspect ratio compared to the cells of the mesh obtained with our proposed approach. Moreover, we remark that in the final iteration of our proposed method only \num{130} of the \num{20022} mesh quality constraints are active together with an additional \num{160} equality constraints for the fixed boundaries, so that the numerical effort for projecting the gradient deformation is very small.

\begin{table}[b]
	\centering
	{\footnotesize
		\caption{Comparison of the mesh quality of the optimized mesh for problem \eqref{eq:shape_stokes}.}
		\label{tab:shape_stokes}
		\setlength{\tabcolsep}{1em}
		\begin{tabular}{l S S}
			\toprule
			\mbox{} & {Classical BFGS} & {BFGS with mesh quality constraints} \\
			\midrule
			Minimum triangle angle [\unit{\degree}] & 5.141 & 24.929 \\
			Maximum aspect ratio & 12.358 & 2.605 \\
			%
			%
			Active mesh quality constraints & {-} & 130 \\
			Equality constraints for fixed boundaries & {-} & 160 \\
			\bottomrule
		\end{tabular}
	}
\end{table}

In Figure~\ref{fig:iterations_stokes} the evolution of these mesh quality measures is shown over the course of the optimization, where the quality of the mesh in each iteration is defined as the quality of the worst mesh cell. For the first couple of iterations, we observe that the steps of the optimization algorithms coincide for both approaches as no mesh quality constraints are active. However, once the mesh quality constraints become active, we observe that the mesh quality stays bounded for our proposed method, whereas it becomes increasingly bad for the classical approach. In particular, we observe in Figure~\ref{fig:iterations_stokes_angle}, that the minimum angle of all cells in the mesh is, indeed, bounded from below by \qty{25}{\degree} during the entire optimization, at least up to the numerical tolerance. This, in turn, also leads to bounds on the aspect ratio of the mesh during the optimization. Hence, this example shows that it is, in fact, sufficient to consider the angles of the mesh cells to guarantee a bound of the mesh quality. Note that the corresponding mesh quality measures for the optimized geometries are summarized in Table~\ref{tab:shape_stokes}. However, we remark that the greatly increased mesh quality comes with a moderate increase in computational resources as our proposed method required nearly twice as many iterations to converge compared to the classical approach. This is mainly due to the fact that during the line search of the gradient projection method the stepsize has to be sufficiently small so that no previously inactive constraints are violated (cf.~Section~\ref{ssec:rosen_finite}).

\subsection{Shape Optimization of a Pipe}

\begin{figure}[b]
	\centering
	\begin{subfigure}{0.6\textwidth}
		\centering
		\includegraphics[width=\textwidth]{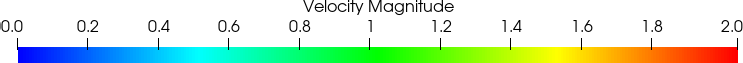}
	\end{subfigure}\\
	\vspace{0.25cm}
	\begin{subfigure}{0.475\textwidth}
		\centering
		\includegraphics[width=\textwidth]{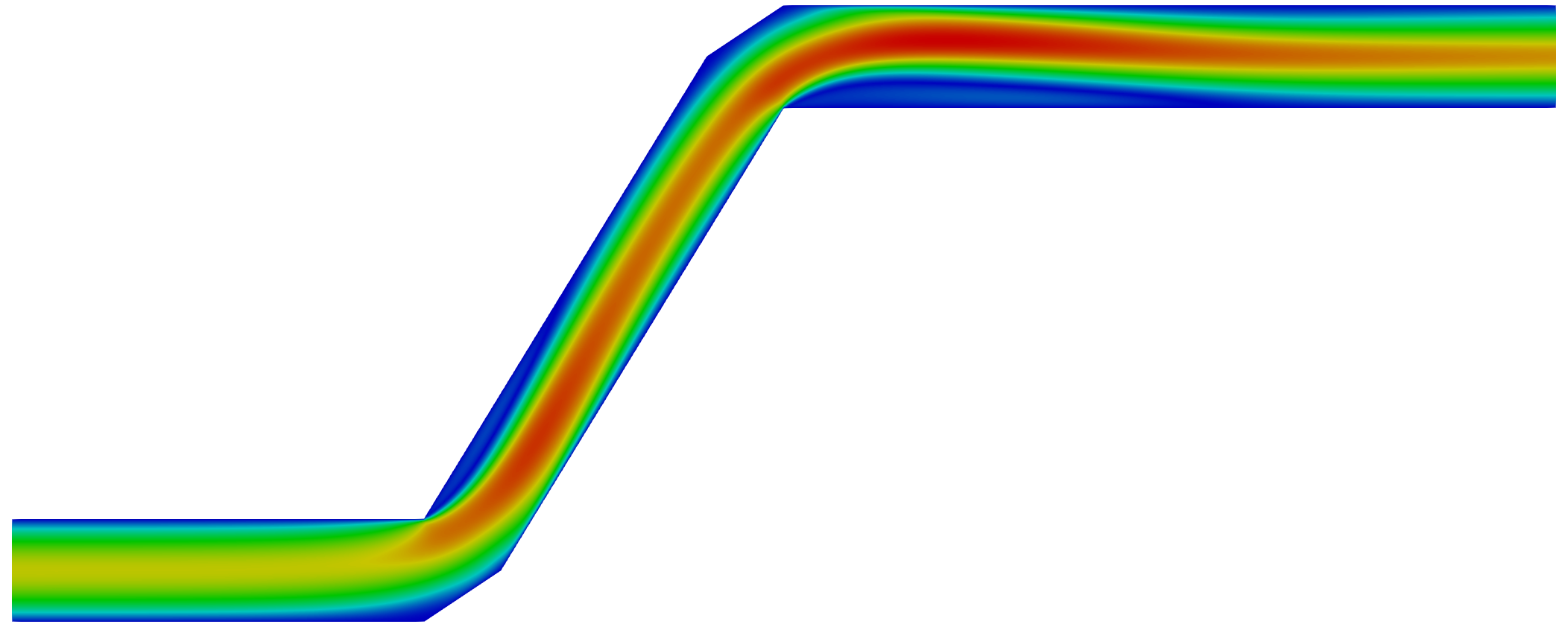}
		\caption{Initial geometry.}
	\end{subfigure}
	\hfil
	\begin{subfigure}{0.475\textwidth}
		\centering
		\includegraphics[width=\textwidth]{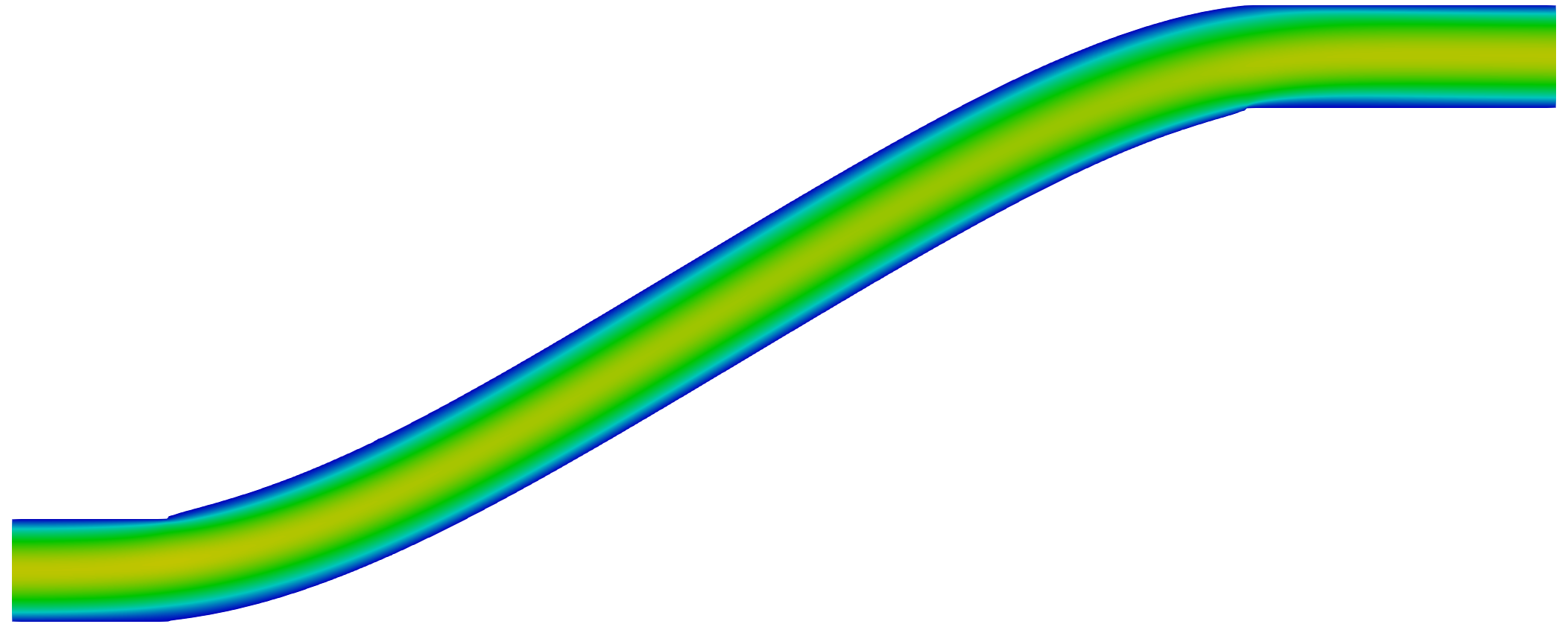}
		\caption{Optimized geometry.}
	\end{subfigure}
	\caption{Velocity magnitude for the Navier-Stokes problem \eqref{eq:pipe}, computed with the gradient descent method with mesh quality constraints.}
	\label{fig:comparison_pipe_velocity}
\end{figure}

As second model example, we consider the optimization of a two-dimensional pipe which has been also considered in, e.g., \cite{Blauth2021Nonlinear, Schmidt2010Efficient, Ham2019Automated}. Here, we again consider the minimization of the fluid's energy dissipation but now consider the incompressible Navier-Stokes equations as model for the flow in the pipe. Similarly to the previous example, we also use a volume constraint to ensure a non-trivial solution to the optimization problem. Let the flow domain be again denoted by $\Omega$ and its boundary $\Gamma = \partial \Omega$ be divided into the inlet $\Gamma_\mathrm{in}$, the wall $\Gamma_\mathrm{wall}$, and the outlet $\Gamma_\mathrm{out}$. For the shape optimization, only one part of the boundary, $\Gamma_\mathrm{def} \subset \Gamma_\mathrm{wall}$, is deformable, whereas the remaining wall boundary, $\Gamma_\mathrm{fix}$, and the in- and outlet remain fixed. The non-dimensionalized steady incompressible Navier-Stokes system is given by
\begin{equation}
	\label{eq:pipe}
	\begin{alignedat}{2}
		-\frac{1}{\mathrm{Re}} \Delta u + (u\cdot \nabla) u + \nabla p &= 0 \quad &&\text{ in } \Omega,\\
		\nabla \cdot u &= 0 \quad &&\text{ in } \Omega,\\
		u &= u_\mathrm{in} \quad &&\text{ on } \Gamma_\mathrm{in},\\
		u &= 0 \quad &&\text{ on } \Gamma_\mathrm{wall},\\
		\frac{1}{\mathrm{Re}}\partial_n u - p n &= 0 \quad &&\text{ on } \Gamma_\mathrm{out}.
	\end{alignedat}
\end{equation}
Here, $u$ and $p$ are, again, the fluid's velocity and pressure, and $\mathrm{Re} > 0$ is the Reynolds number. Using the notation from \eqref{eq:notation_volume}, we can formulate the shape optimization problem as
\begin{equation}
	\label{eq:shape_pipe}
	\min_{\Omega \in \admissiblegeom} \mathcal{J}(\Omega, u) = \integral{\Omega} \frac{1}{\mathrm{Re}} Du : Du \dmeas{x} + \frac{\nu_\mathrm{vol}}{2} \left( \mathrm{vol}(\Omega) - \mathrm{vol}(\Omega^0) \right)^2 \quad \text{ s.t. } \eqref{eq:pipe},
\end{equation}
where set of admissible geometries for this problem is given by
\begin{equation*}
	\admissiblegeom = \Set{\Omega \subset \R^d | \Gamma_\mathrm{in} = \Gamma_\mathrm{in}^0, \Gamma_\mathrm{fix} = \Gamma_\mathrm{fix}^0, \Gamma_\mathrm{out} = \Gamma_\mathrm{out}^0 }.
\end{equation*}
Here $\Omega^0$ denotes the initial geometry and $\Gamma^0 = \partial \Omega^0$ is its boundary. For the shape derivative and corresponding adjoint system of this problem we refer the reader, again, to \cite{Blauth2021Adjoint}. Note that the solution of the Navier-Stokes equations \eqref{eq:pipe} can be seen in Figure~\ref{fig:comparison_pipe_velocity}, where the velocity magnitude is shown on the initial and optimized geometries.

We consider this problem numerically in two dimensions. For the finite element mesh we use a uniform mesh generated with Gmsh \cite{Geuzaine2009Gmsh} consisting of \num{14856} nodes and \num{28708} triangles. For the discretization of the Navier-Stokes system we use the LBB-stable Taylor-Hood elements as discussed in the previous section. The inlet velocity is given by the parabolic profile $u_\mathrm{in}(x) = 6~(1 - x_2)~(x_2 - 0)$. The Reynolds number for the flow is chosen as $\mathrm{Re} = 400$ and the parameter for the volume regularization is chosen as $\nu_\mathrm{vol} = 1$. We again use the equations of linear elasticity \eqref{eq:linear_elasticity} for computing the gradient deformation with parameters $\mu_\mathrm{elas} = 1$, $\lambda_\mathrm{elas} = 0$, and $\delta_\mathrm{elas} = 0$. To solve the problem numerically, we employ the gradient descent method implemented in our software cashocs. A relative stopping criterion with tolerance of \num{1e-3} is used to terminate the method once the norm of the gradient deformation falls below this threshold. Analogously to the previous example, for the mesh quality constraints we use a minimum angle threshold of $\alpha_\mathrm{thr} = \SI{0.436}{\radian}$ corresponding to $\SI{25}{\degree}$ and a numerical tolerance of \num{1e-2}.

\begin{figure}[t]
	\centering
	\begin{subfigure}{0.475\textwidth}
		\centering
		\includegraphics[width=\textwidth]{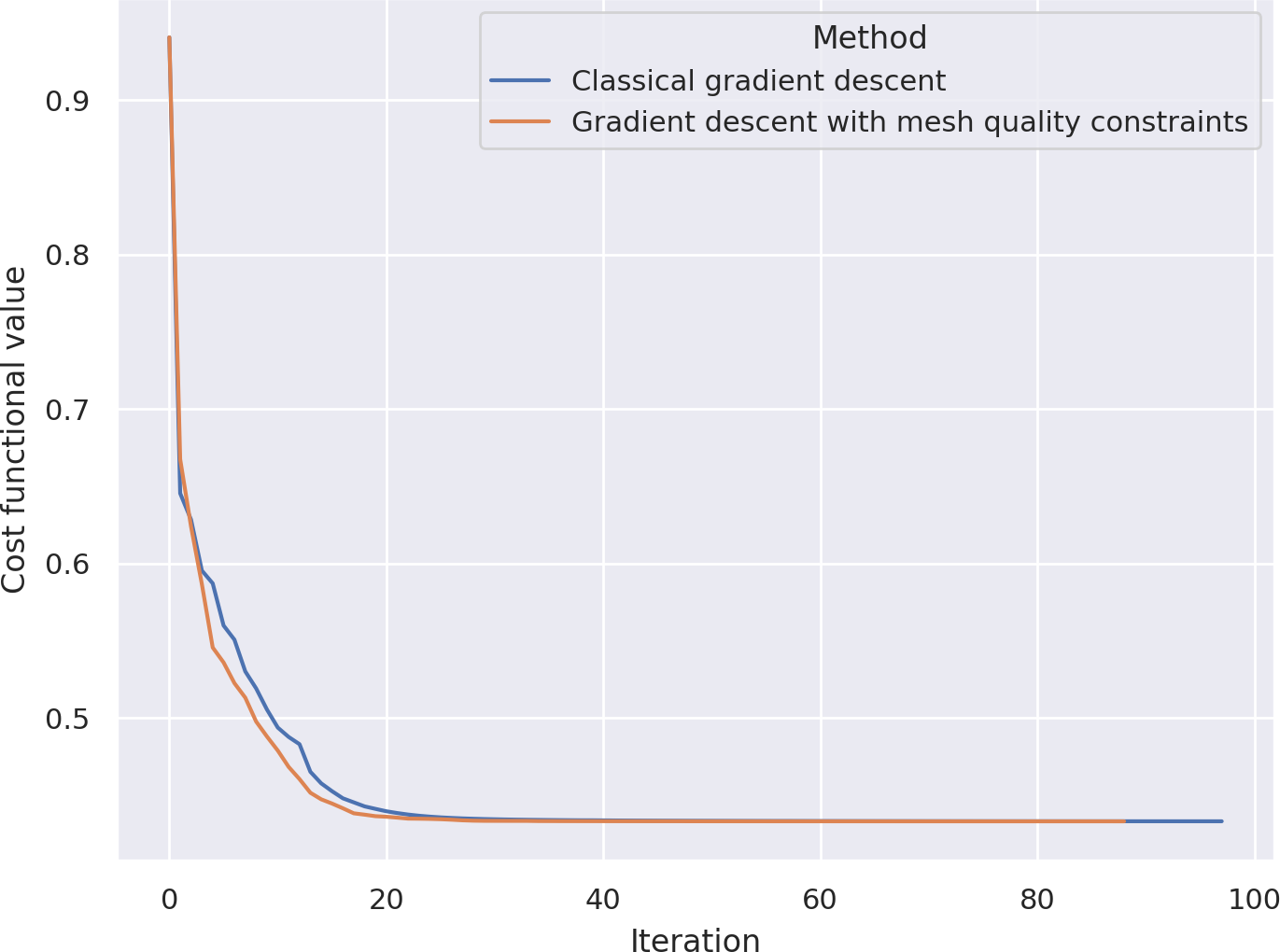}
		\caption{Cost functional value.}
	\end{subfigure}
	\hfil
	\begin{subfigure}{0.475\textwidth}
		\centering
		\includegraphics[width=\textwidth]{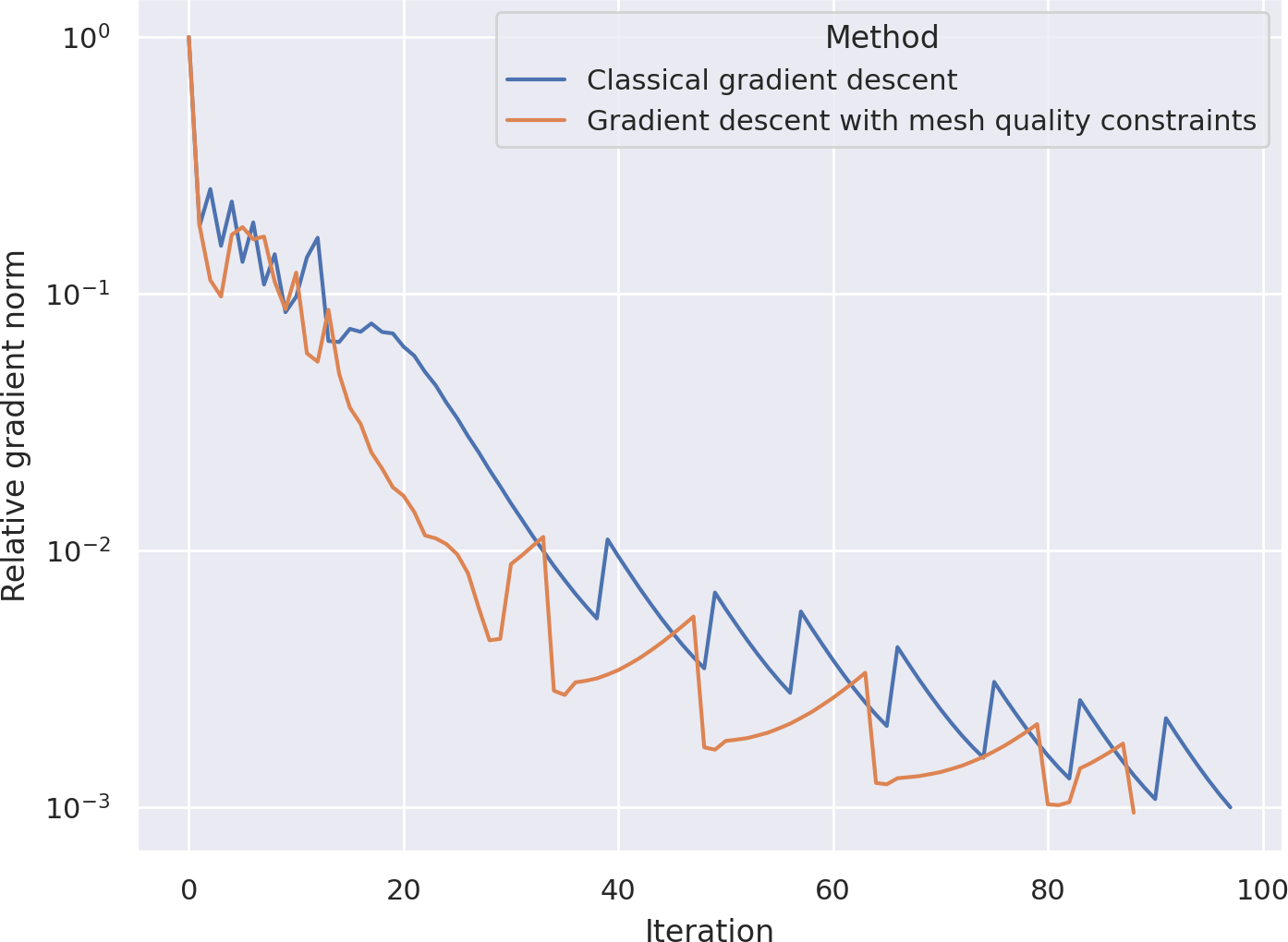}
		\caption{Relative norm of the gradient deformation.}
	\end{subfigure}
	\caption{History of the optimization algorithm for problem \eqref{eq:shape_pipe}.}
	\label{fig:history_pipe}
\end{figure}

\begin{figure}[b]
	\centering
	\begin{subfigure}{0.6\textwidth}
		\centering
		\includegraphics[width=\textwidth]{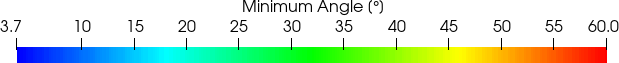}
	\end{subfigure}\\
	\begin{subfigure}[t]{0.475\textwidth}
		\centering
		\includegraphics[width=\textwidth]{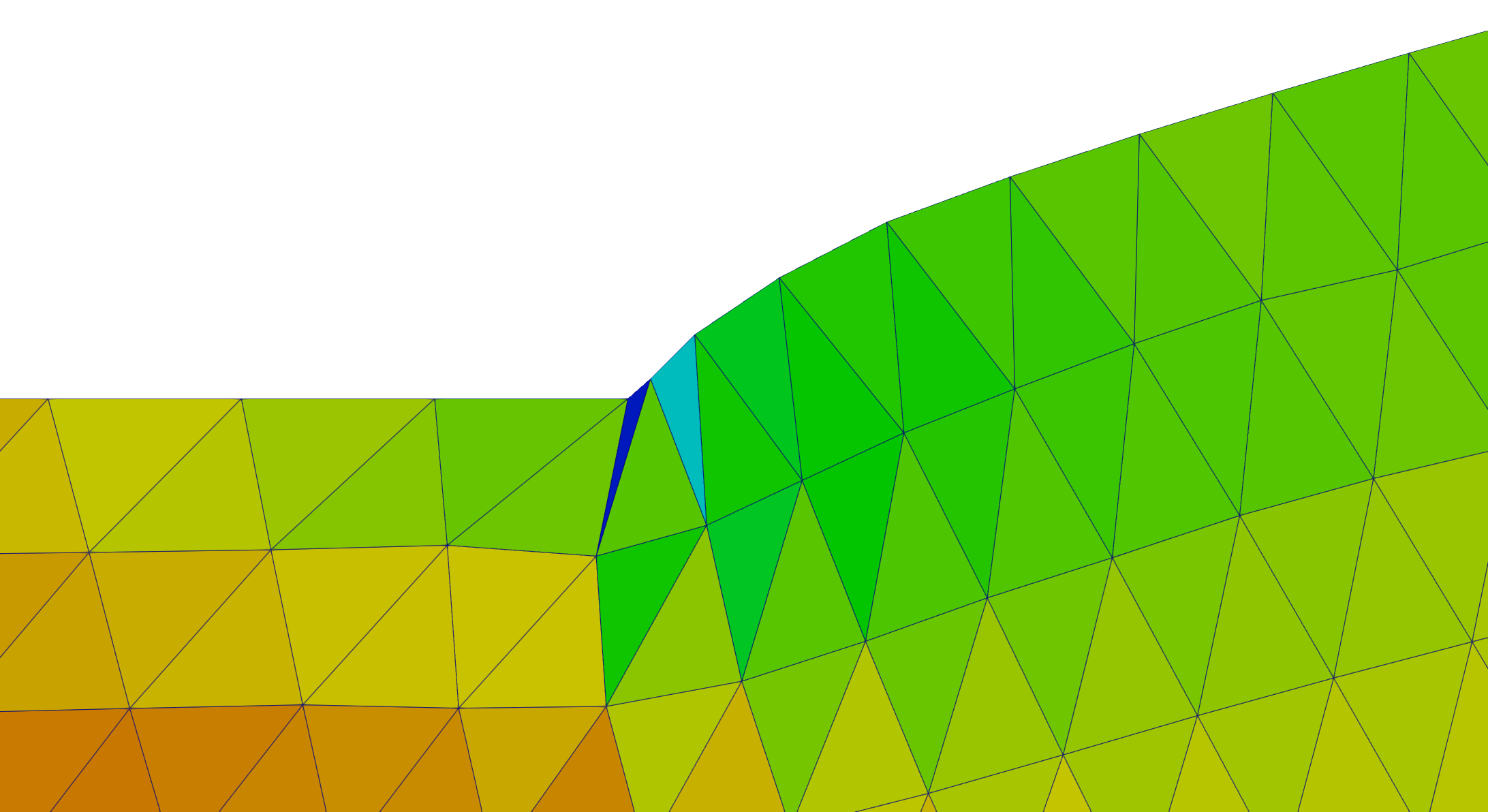}
		\caption{Classical gradient descent method.}
	\end{subfigure}
	\hfil
	\begin{subfigure}[t]{0.475\textwidth}
		\centering
		\includegraphics[width=\textwidth]{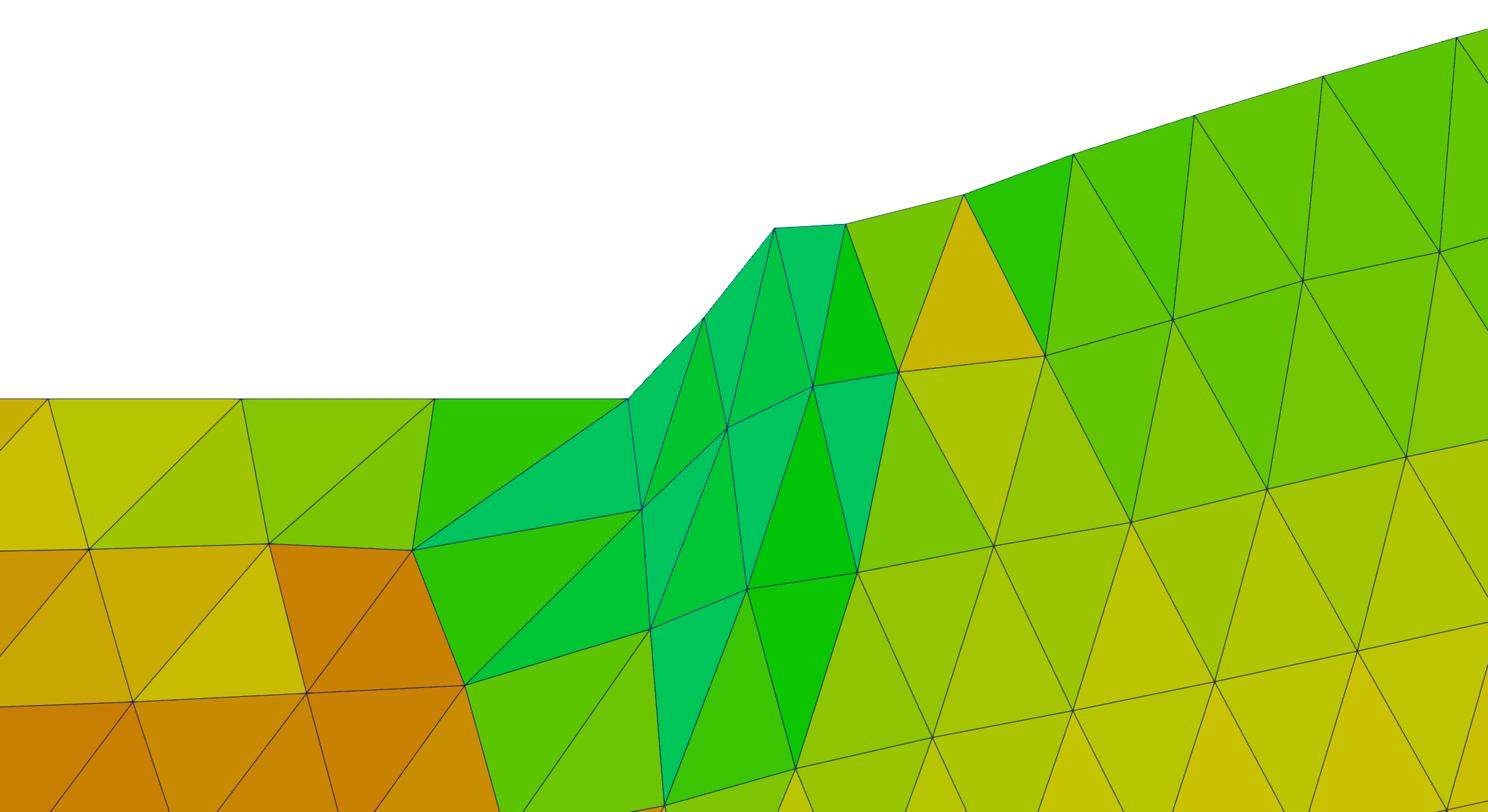}
		\caption{Gradient descent method with mesh quality constraints.}
	\end{subfigure}
	\caption{Optimized geometry with mesh and mesh quality for problem \eqref{eq:shape_pipe} -- close up of the lower kink.}
	\label{fig:quality_pipe}
\end{figure}

For this problem, the classical gradient descent method needed 97 iterations to converge, whereas the gradient descent method with mesh quality constraints propsed in this paper only required 88 iterations. The history of the objective function as well as the norm of the gradient deformation can be found in Figure~\ref{fig:history_pipe}. There, we can see that the cost functional value is nearly identical for both approaches and that our proposed method only has a slightly smaller cost functional value during the first 20 iterations. After about 20 iterations, both approaches reach the final cost functional value and no further decrease is happening for either method. For the norm of the gradient deformation these findings are similar, however, we observe that our approach converged slightly faster than the classical one for this particular problem.

In Figure~\ref{fig:quality_pipe} the computational mesh of the optimized geometry is investigated detailedly near the lower kink of the geometry. Note that the kinks occur at the interface between deformable and fixed boundary and that they are the most problematic zones for the mesh quality (cf.\ \cite{Mueller2023Shape}). We observe that the classical approach indeed has some problematic mesh cells near these kinks where the mesh quality deteriorates. For our proposed approach, on the other hand, the mesh quality does not deteriorate in these areas. The optimized geometry resembles a kind of bulge for both approaches. However, it can be clearly seen that the bulge is significantly more pronounced with the approach using the mesh quality constraints. Apart from that small and localized difference, the optimized geometries for both considered approaches coincide.

\begin{Remark}
	We note that the kinks shown in Figure~\ref{fig:quality_pipe} are undesirable from the practical point of view, particularly the bulge generated by our approach which considers the mesh quality constraints. To obtain smoother shapes for practical applications, a curvature-based regularization could be applied which would penalize the occurrence of kinks. We refer the reader to, e.g., \cite{Sundaramoorthi2009New, Dogan2012First, Strandmark2011Curvature}, where curvature regularization and the corresponding shape functionals are used and investigated.
\end{Remark}

\begin{figure}[t]
	\centering
	\begin{subfigure}[t]{0.475\textwidth}
		\centering
		\includegraphics[width=\textwidth]{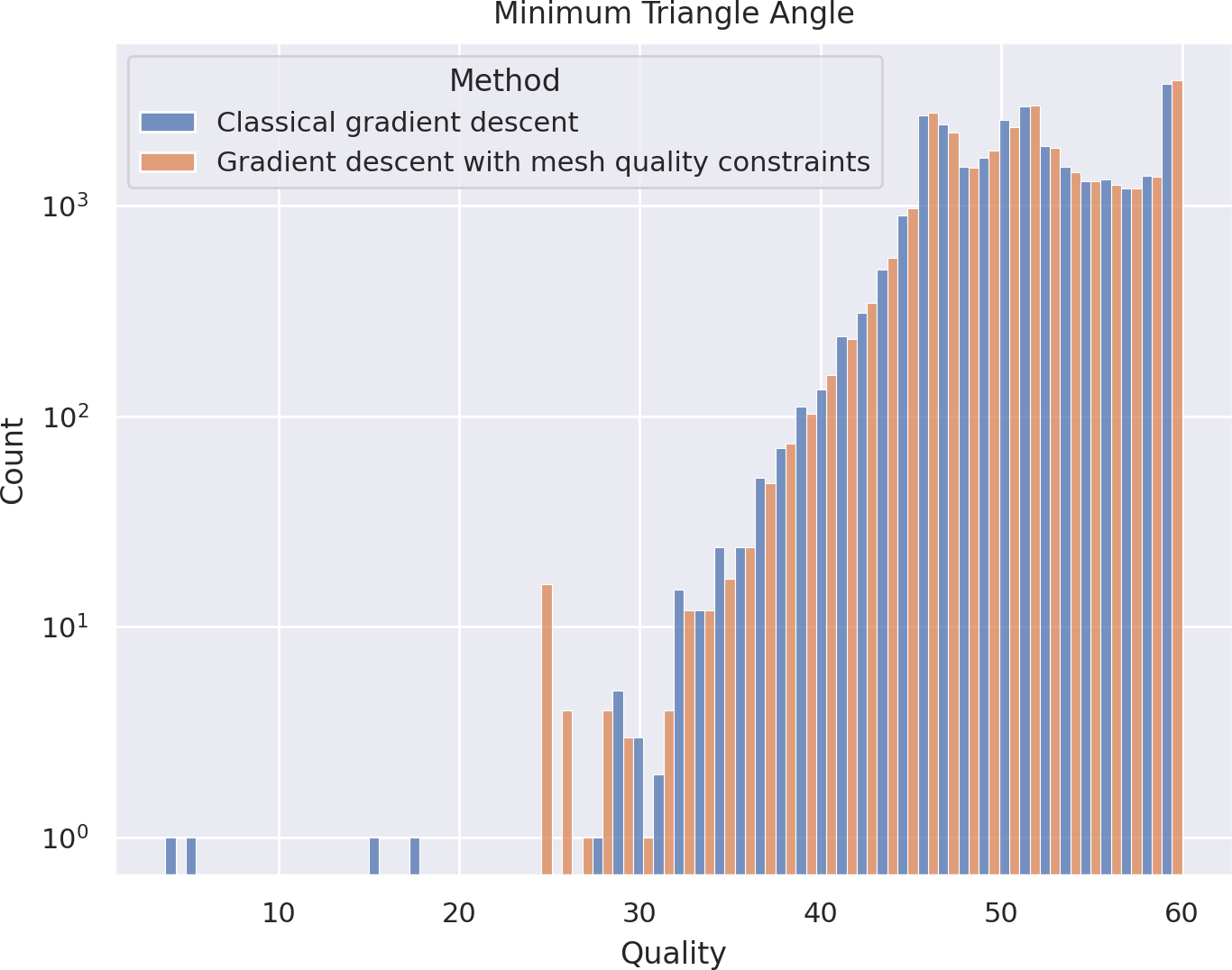}
		\caption{Minimum triangle angle in degrees.}
		\label{fig:dist_pipe_angle}
	\end{subfigure}
	\hfil
	\begin{subfigure}[t]{0.475\textwidth}
		\centering
		\includegraphics[width=\textwidth]{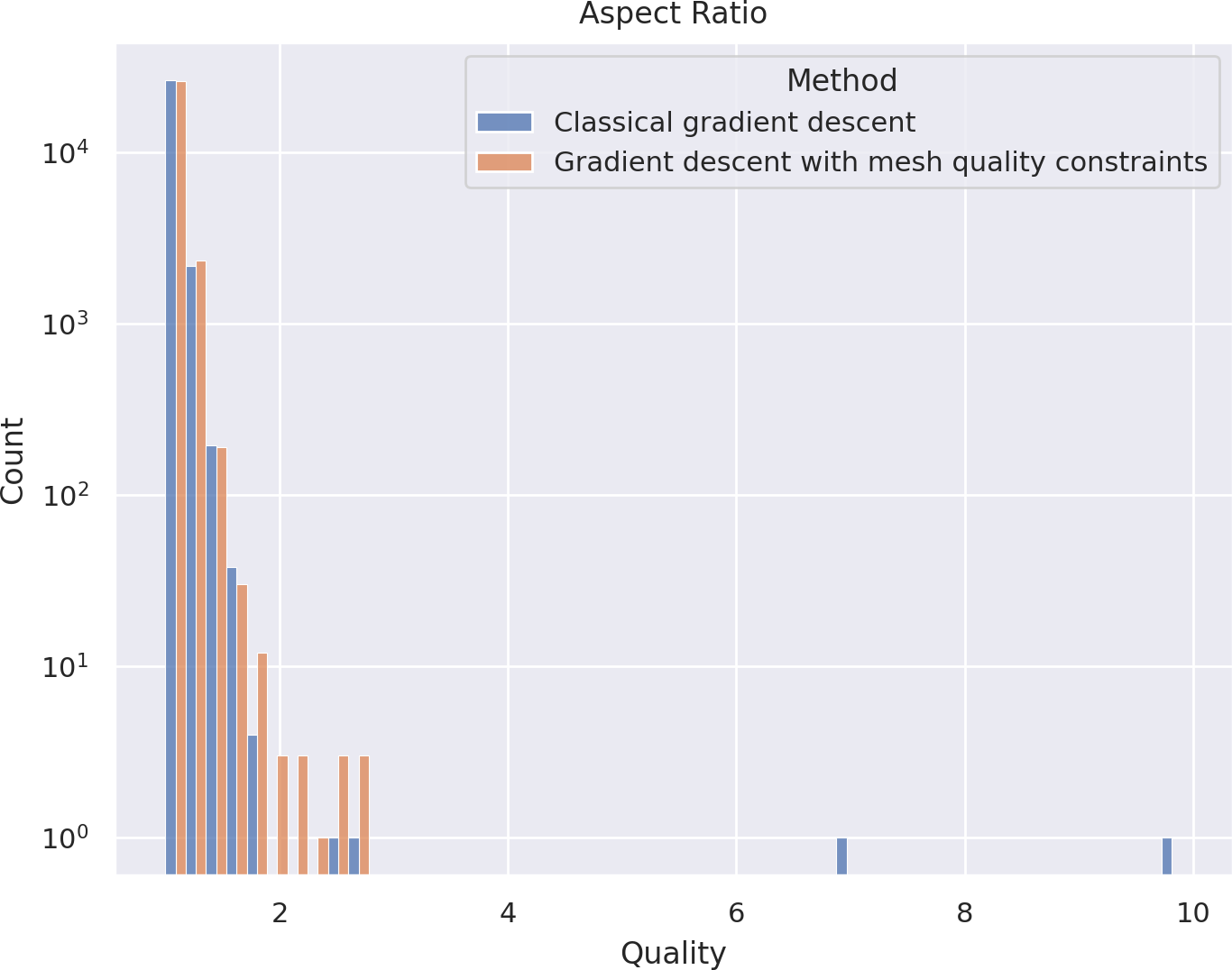}
		\caption{Aspect ratio.}
		\label{fig:dist_pipe_aspect_ratio}
	\end{subfigure}
	\caption{Histograms showing the distribution of the mesh quality after the shape optimization for problem \eqref{eq:shape_pipe}.}
	\label{fig:dist_pipe}
\end{figure}

\begin{figure}[b]
	\centering
	\begin{subfigure}[t]{0.475\textwidth}
		\centering
		\includegraphics[width=\textwidth]{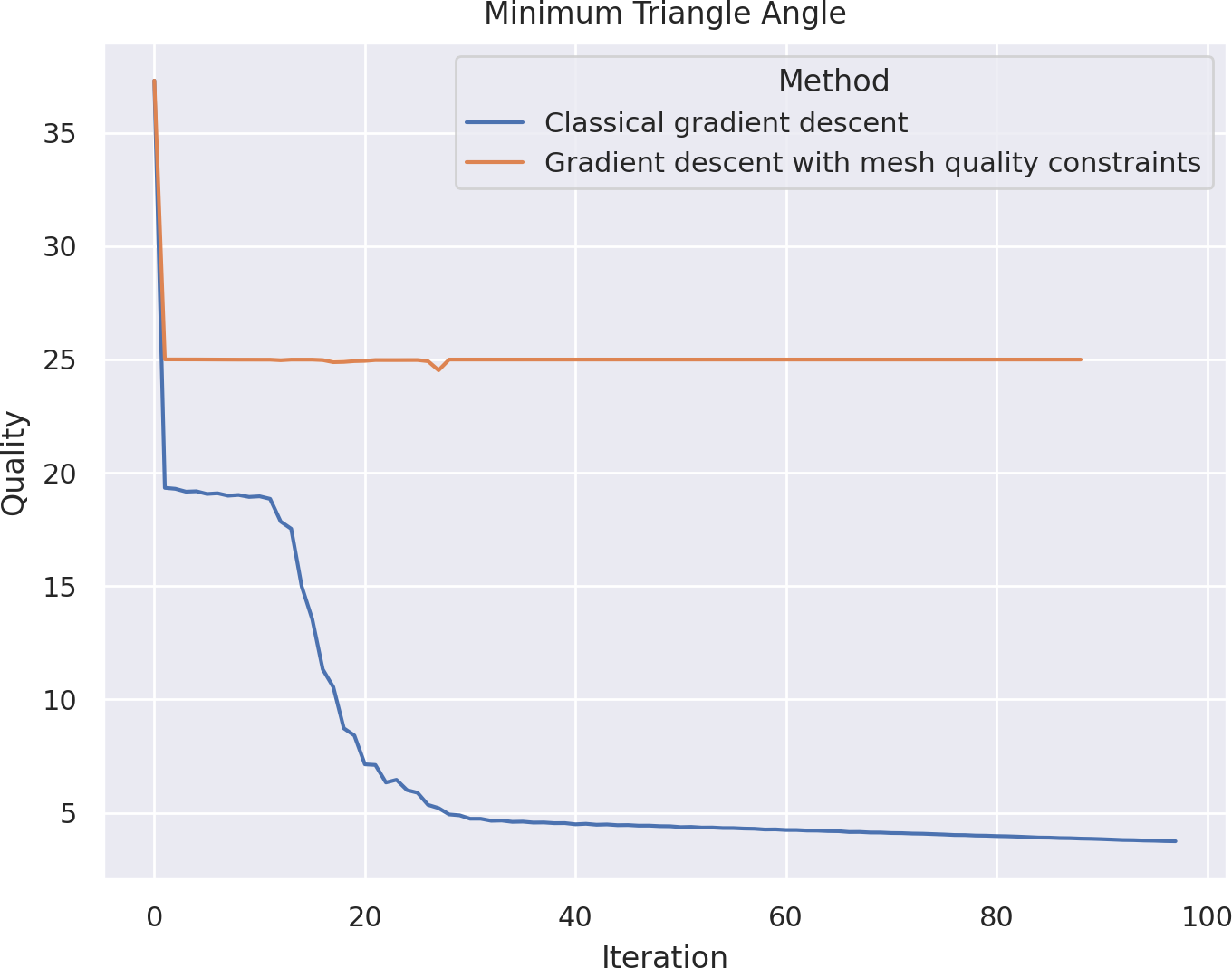}
		\caption{Minimum triangle angle in degrees.}
		\label{fig:iterations_pipe_angle}
	\end{subfigure}
	\hfil
	\begin{subfigure}[t]{0.475\textwidth}
		\centering
		\includegraphics[width=\textwidth]{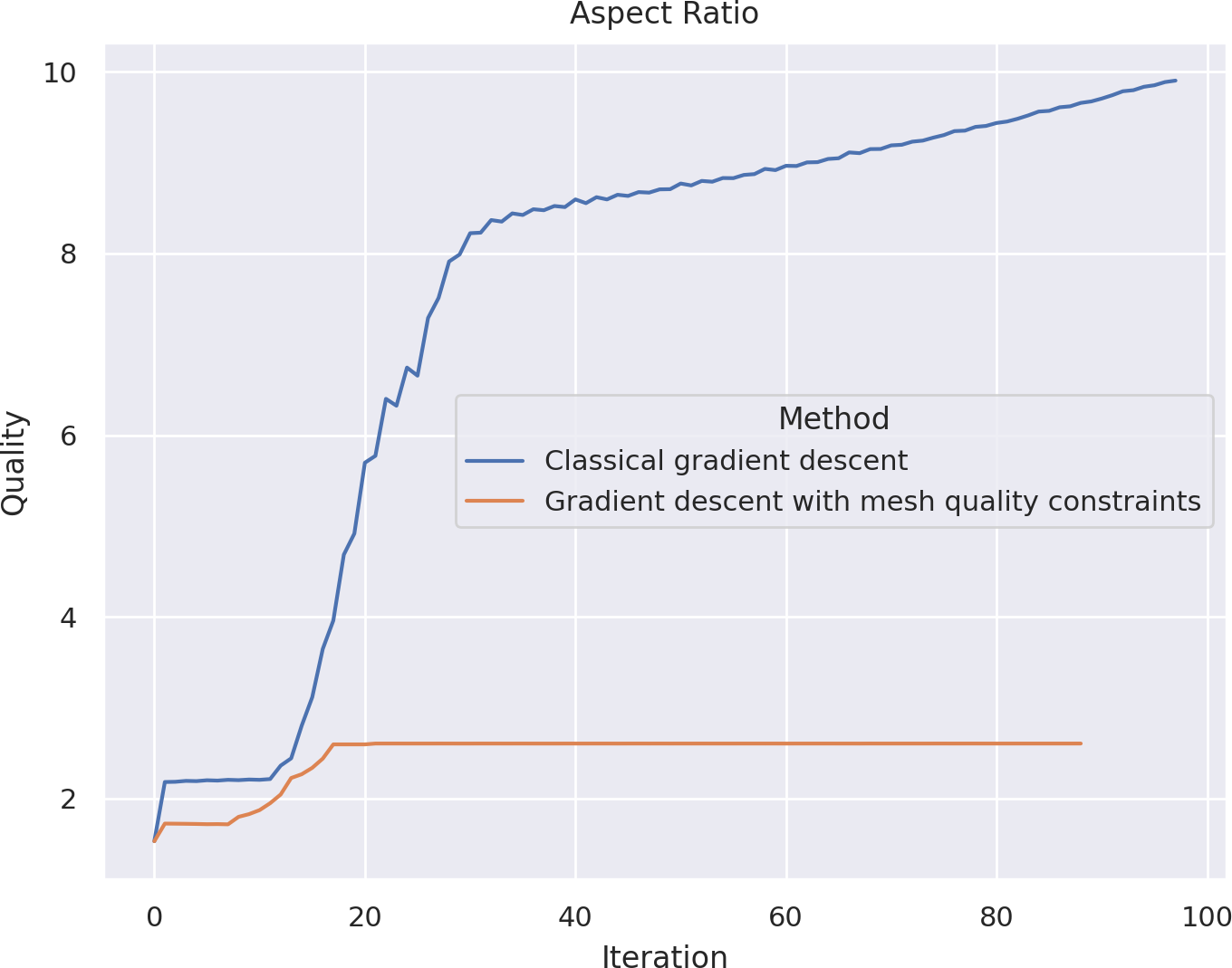}
		\caption{Aspect ratio.}
		\label{fig:iterations_pipe_aspect_ratio}
	\end{subfigure}
	\caption{Evolution of the mesh quality during the shape optimization for problem \eqref{eq:shape_pipe}.}
	\label{fig:iterations_pipe}
\end{figure}

Histograms of the mesh quality for the optimized designs for both approaches are shown in Figure~\ref{fig:dist_pipe} for the minimum triangle angle and the aspect ratio of the mesh cells. Analogously to our previous example, we observe that our proposed method is able to enforce bounds on the mesh quality. In fact, the prescribed lower bound of \SI{25}{\degree} for the minimum angle is attained with our method. In contrast, for the classical approach we see that there are two problematic mesh cells with a minimum angle below \SI{10}{\degree}, which are located at the kinks of the optimized geometry. For these, both the minimum triangle angle and aspect ratio are significantly worse compared to the method with mesh quality constraints. Note that to achieve the bounds on the mesh quality, only 26 of the total \num{86124} mesh quality constraints were active in the final iteration of our method. Additionally, \num{592} equality constraints were considered for the fixed boundaries.

In Figure~\ref{fig:iterations_pipe} the evolution of the mesh quality over the course of the shape optimization is visualized. Here, the mesh quality is, again, defined as the quality of the worst mesh cell. We observe that the mesh quality directly falls off for the classical gradient descent method and continues to deteriorate over the course of the optimization. For our proposed method, on the other hand, the minimum angle remains bounded from below by the specified minimum angle of \SI{25}{\degree}. Consequently, also the aspect ratio of the mesh remains bounded for our proposed approach, which is not the case for the classical one. Finally, the mesh quality measures for the optimized geometries are summarized in Table~\ref{tab:shape_pipe}. 

\begin{table}[t]
	\centering
	{\footnotesize
		\caption{Comparison of the mesh quality of the optimized mesh for problem \eqref{eq:shape_pipe}.}
		\label{tab:shape_pipe}
		\setlength{\tabcolsep}{1em}
		\begin{tabular}{l S S}
			\toprule
			\mbox{} & {Classical gradient descent} & {Gradient descent with mesh quality constraints} \\
			\midrule
			Minimum triangle angle [\unit{\degree}] & 3.745 & 24.994 \\
			Maximum aspect ratio & 9.903 & 2.604 \\
			%
			%
			Active mesh quality constraints & {-} & 26 \\
			Equality constraints for fixed boundaries & {-} & 592 \\
			\bottomrule
		\end{tabular}
	}
\end{table}

\subsection{Shape Optimization of a Structured Packing for Distillation}

Let us now consider a large-scale, three-dimensional shape optimization problem which highlights the capabilities and practical importance of our proposed method. Note that this example is taken from \cite{Blauth2025CFD} and considers the shape optimization of a structured packing for distillation. We refer the reader to our previous work \cite{Blauth2025CFD}, where this problem and the shape optimization results are discussed in greater detail, as this is beyond the scope of this paper. In this paper, we focus on the aspect how the proposed mesh quality constraints enable the numerical solution of this problem. 

\begin{figure}[t]
	\centering
	\begin{subfigure}[t]{0.5\textwidth}
		\centering
		\includegraphics[height=0.225\textheight]{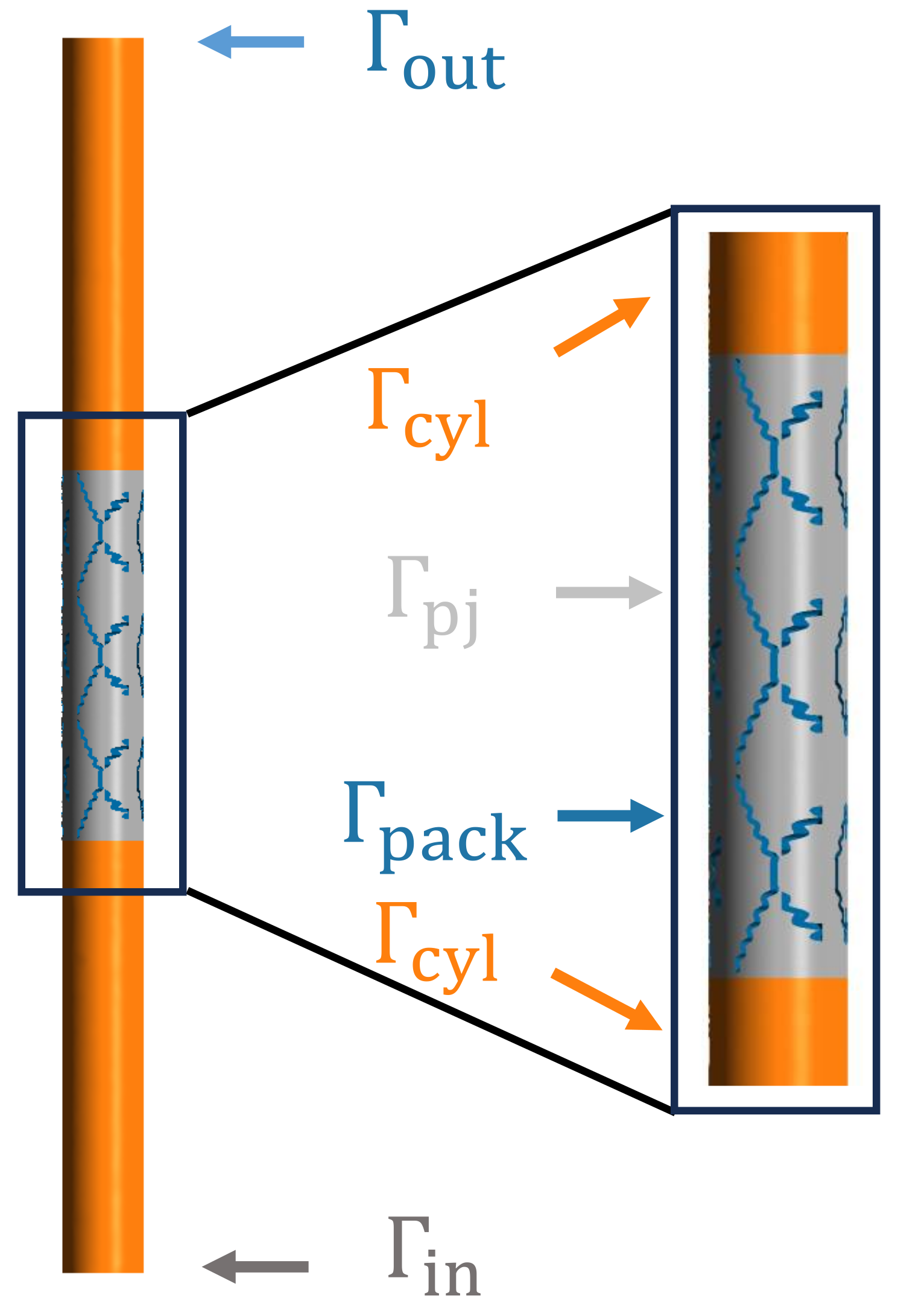}
		\caption{Geometrical description of the boundaries.}
	\end{subfigure}%
	\begin{subfigure}[t]{0.5\textwidth}
		\centering
		\includegraphics[height=0.225\textheight]{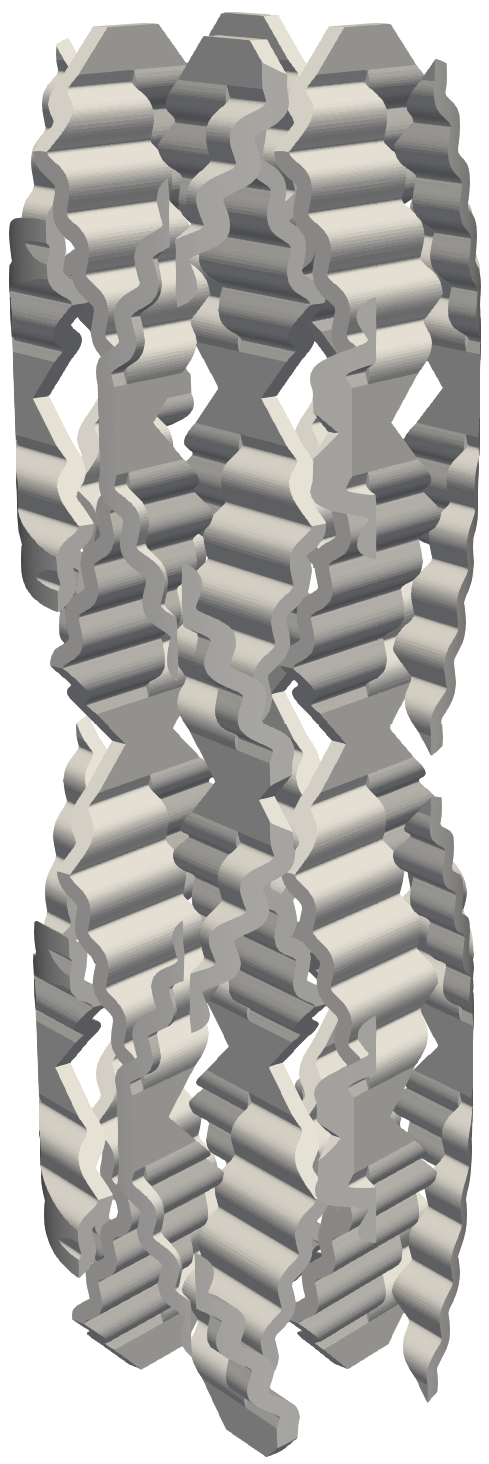}
		\caption{CAD image of the initial RP9M-3D packing.}
		\label{fig:rp9m_3d}
	\end{subfigure}
	\caption{Geometrical setup for the shape optimization of a structured packing for distillation taken from \cite{Blauth2025CFD}.}
	\label{fig:geometry_packing}
\end{figure}

We consider the following simplified single phase model of the mass transfer in the distillation column from \cite{Neukaeufer2022Development}, which consists of the incompressible Navier-Stokes equations and a convection-diffusion equation for the gas phase flow. The former is given by
\begin{equation}
	\label{eq:navier_stokes}
	\begin{alignedat}{2}
		-\mu \Delta u + \rho (u\cdot \nabla) u + \nabla p &= 0 \quad &&\text{ in } \Omega,\\
		\nabla \cdot u &= 0 \quad &&\text{ in } \Omega, \\
		u &= u_\mathrm{in} \quad &&\text{ on } \Gamma_\mathrm{in},\\
		u &= 0 \quad &&\text{ on } \Gamma_\mathrm{wall},\\
		\mu \partial_n u - pn &= 0 \quad &&\text{ on } \Gamma_\mathrm{out},
	\end{alignedat}
\end{equation}
where $u$ is the fluid velocity, $p$ is its pressure, $\mu$ is the dynamic viscosity of the fluid, and $\rho$ is its density. The boundary of the fluid domain $\Omega$ is divided into inlet $\Gamma_\mathrm{in}$, outlet $\Gamma_\mathrm{out}$, and wall boundary $\Gamma_\mathrm{wall}$. The latter consists of the cylinder jacket $\Gamma_\mathrm{cyl}$, packing jacket $\Gamma_\mathrm{pj}$, and packing $\Gamma_\mathrm{pack}$. The geometrical setup of this problem is shown in Figure~\ref{fig:geometry_packing}. To model the mass transfer in the distillation column, the following convection-diffusion equation is used
\begin{equation}
	\label{eq:concentration}
	\begin{alignedat}{2}
		- \nabla \cdot \left( \tilde{D} \nabla c \right) + u \cdot \nabla c &= 0 \quad &&\text{ in } \Omega,\\
		c &= c_\mathrm{in} \quad &&\text{ on } \Gamma_\mathrm{in},\\
		c &= c_\mathrm{pack} \quad &&\text{ on } \Gamma_\mathrm{pj} \cup \Gamma_\mathrm{pack},\\
		D\partial_n c &= 0 \quad &&\text{ on } \Gamma_\mathrm{cyl} \cup \Gamma_\mathrm{out},
	\end{alignedat}
\end{equation}
where $c$ is a fictitious concentration and $\tilde{D}$ is the diffusion coefficient. For a detailed discussion of the simplified model and the simplifying assumptions used for deriving it we refer the reader to \cite{Blauth2025CFD,Neukaeufer2022Development}. As model parameters, we use the ones from \cite{Blauth2025CFD}, i.e., $\mu = \qty{1.728e-5}{\pascal \second}$, $\rho = \qty{1.138}{\kilogram \per \cubic \meter}$ and $\tilde{D} = \qty{3.72e-6}{\square \meter \per \second}$, and for the boundary conditions we use $u_\mathrm{in} = \qty{0.933}{\meter \per \second}$, $c_\mathrm{in} = \qty{100}{\mol \per \cubic \meter}$, and $c_\mathrm{pack} = \qty{1}{\mol \per \cubic \meter}$. 

We consider the following shape optimization problem from \cite{Blauth2025CFD}
\begin{equation}
	\label{eq:shape_distillation}
	\max_{\Omega \in \mathcal{A}} \mathcal{J}(\Omega, c) = \beta(c) = \frac{\dot{V}}{A_\mathrm{geo}} \log \left( \frac{c_\mathrm{pack} - c_\mathrm{in}}{c_\mathrm{pack} - c_\mathrm{out}} \right) \qquad \text{ s.t. } \eqref{eq:navier_stokes} \text{ and } \eqref{eq:concentration}.
\end{equation}
Here, $\beta$ is the logarithmic mass transfer coefficient which depends on the volumetric flow rate of the fluid $\dot{V}$, the surface area of the packing and packing wall $A_\mathrm{geo}$, and the flow-averaged outlet concentration $c_\mathrm{out}$, which is defined as
\begin{equation*}
	c_\mathrm{out} = \frac{\integral{\Gamma_\mathrm{out}} \rho (u \cdot n) c \dmeas{s} }{\integral{\Gamma_\mathrm{out}} \rho u \cdot n \dmeas{s}}.
\end{equation*}
The goal of the optimization problem \eqref{eq:shape_distillation} is to increase the logarithmic mass transfer coefficient $\beta$, which has been successfully used in \cite{Neukaeufer2022Development} as a qualitative criterion for the separation efficiency of the distillation column. For the shape optimization, we consider the following set of admissible domains
\begin{equation*}
	\mathcal{A} = \Set{\Omega \subset \R^3 | \Omega \subset \holdall, \Gamma_\mathrm{in} = \Gamma_\mathrm{in}^0, \Gamma_\mathrm{cyl} = \Gamma_\mathrm{cyl}^0, \Gamma_\mathrm{out} = \Gamma_\mathrm{out}^0},
\end{equation*}
where the hold-all domain $\holdall$ is given by the cylinder shown in Figure~\ref{fig:geometry_packing} and $\Omega^0$ is the initial geometry with boundary $\Gamma^0$, i.e., only the shape of the packing jacket and packing is allowed to change. This is further restricted by only allowing geometrical changes of the packing jacket in $z$-direction, so that the overall cylindrical shape is maintained.

The domain is meshed with Ansys\textsuperscript{\textregistered} Meshing and the resulting mesh consists of about \num{13.7}~million tetrahedrons and \num{2.35}~million nodes. Note that the initial geometry contains some regions where the packing is nearly tangential to the packing jacket. These regions could not be meshed with high quality cells so that there are some cells in the initial mesh which do have a comparatively low quality. 

\begin{Remark}
	\label{rem:gmsh}
	Usually, the meshes for shape optimization with cashocs are generated with Gmsh \cite{Geuzaine2009Gmsh} as cashocs has an automated remeshing workflow for such meshes. For this problem, it was not possible to generate a mesh with a sufficiently high quality with Gmsh so that the commercial software Ansys\textsuperscript{\textregistered} Meshing was used to generate a finite element mesh. However, for this approach, no automated remeshing is available with cashocs.
\end{Remark}

To discretize the state equations we proceed as follows: The Navier-Stokes equations are discretized with linear Lagrange elements for both velocity and pressure. As this choice of elements is not LBB-stable, a pressure-stabilized Petrov-Galerkin (PSPG) stabilization is used. Moreover, a streamline-upwind Petrov-Galerkin (SUPG) stabilization is used for stabilizing the effects of convection and a least squares incompressibility constraint (LSIC) stabilization is used to enhance the conservation of mass. For a detailed discussion of these stabilization approaches we refer the reader, e.g., to \cite{John2016Finite, Elman2014Finite}. For the discretization of the convection-diffusion equation \eqref{eq:concentration}, we use linear Lagrange elements. To stabilize the equation, both an SUPG and crosswind stabilization are used \cite{Codina1993discontinuity}. For more details regarding the nonlinear and linear solvers used for the numerical solution of this problem, the reader is referred to \cite{Blauth2025CFD}.

\begin{figure}[b]
	\centering
	\begin{subfigure}{0.475\textwidth}
		\centering
		\includegraphics[width=\textwidth]{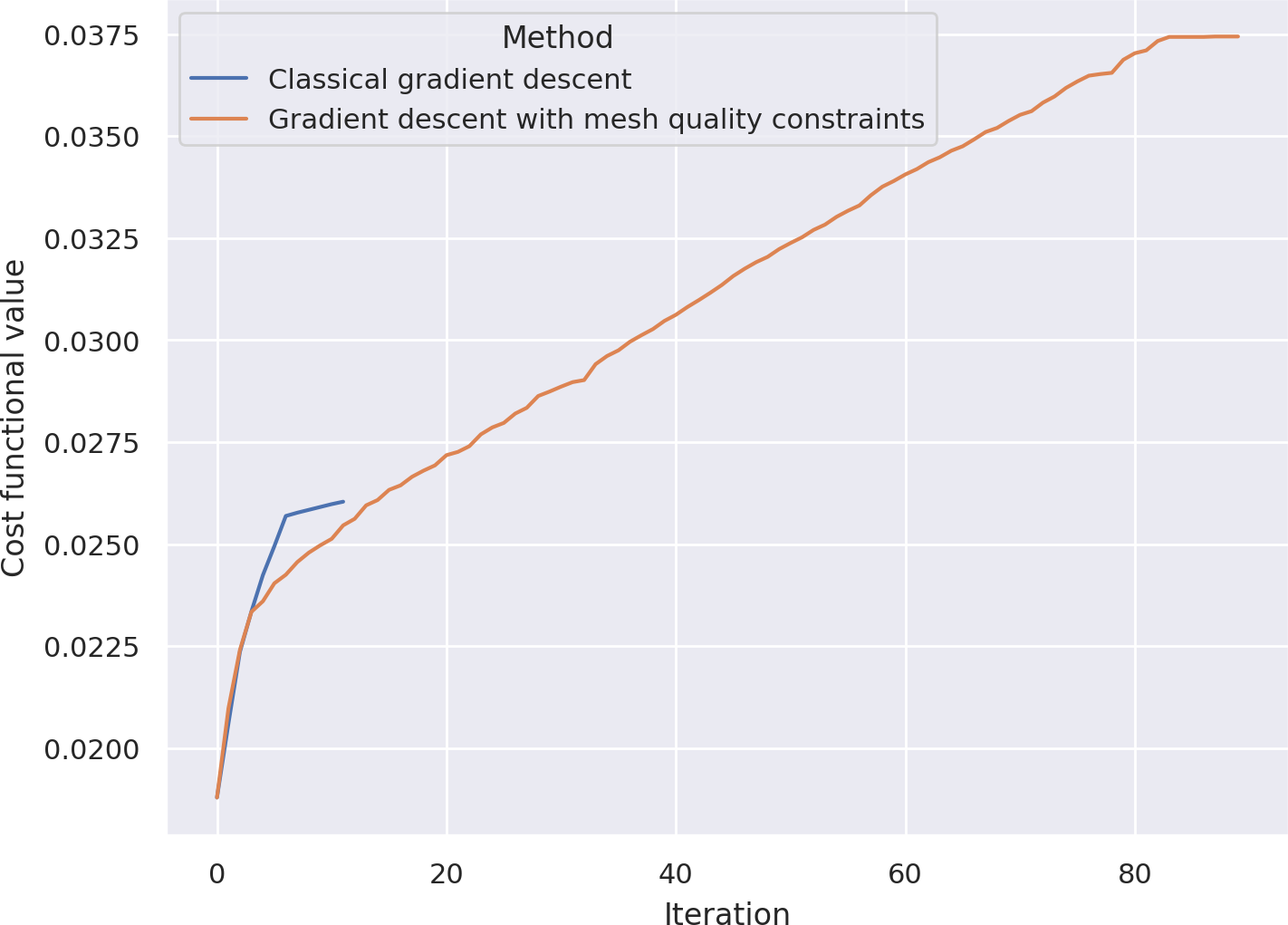}
		\caption{Cost functional value.}
	\end{subfigure}
	\hfil
	\begin{subfigure}{0.475\textwidth}
		\centering
		\includegraphics[width=\textwidth]{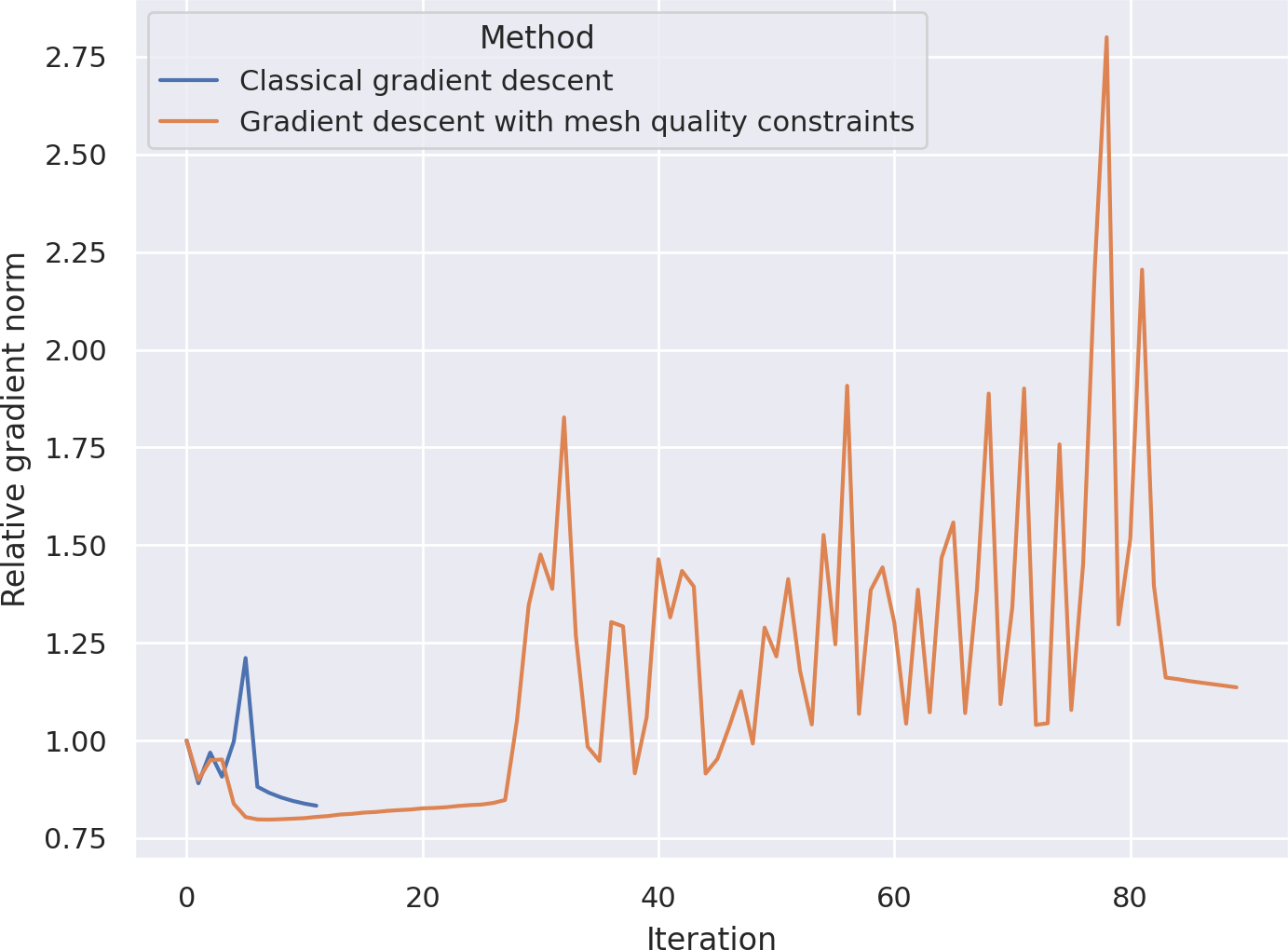}
		\caption{Relative norm of the gradient deformation.}
	\end{subfigure}
	\caption{History of the optimization algorithm for problem \eqref{eq:shape_distillation}.}
	\label{fig:history_packing}
\end{figure}

For the numerical solution of the shape optimization problem, we use the shape of the RP9M-3D packing \cite{Neukaeufer2022Development} as initial geometry, which is shown in Figure~\ref{fig:rp9m_3d}. We employ the gradient descent method for solving the shape optimization problem \eqref{eq:shape_distillation} numerically and use a fixed amount of 90 iterations. To solve this large-scale problem numerically, we use \num{32} cores of an Intel\textsuperscript{\textregistered} Xeon\textsuperscript{\textregistered} Gold 6240R CPU. To prevent mesh zones from overlapping, the domain of the packing, which is not part of the flow domain, is also meshed and the state equations are extended by zero in this region. However, the computational domain for the linear elasticity equations used to determine the gradient deformation also includes the packing domain. Our proposed mesh quality constraints ensure that the cells which discretize the packing cannot be compressed arbitrarily so that the flow regions cannot overlap. For the computation of the gradient deformation the linear elasticity equations are used with $\mu_\mathrm{elas} = 1$, $\lambda_\mathrm{elas} = 0$, and $\delta_\mathrm{elas} = 0$. These parameters are weighted w.r.t.\ the mesh element size as discussed in \cite{Blauth2021Model}, i.e., the (artificial) stiffness of the mesh cells is weighted by one over the volume of each element. This leads to an inhomogeneous stiffness of the mesh cells which ensures that large elements have a lower stiffness than small ones as they tend to \qe{absorb} deformations better. For the mesh quality constraints, we use the second approach discussed in Section~\ref{ssec:choice_of_angle}, i.e., the minimum solid angle threshold for this problem is different for each cell and a relative tolerance of $\nu = 0.25$ is used here. This approach is used due to the fact that there are some bad mesh cells in the initial mesh, so that using a single global minimum angle threshold is not appropriate (cf.~Section~\ref{ssec:choice_of_angle}). Moreover, we use a numerical tolerance of \num{2.5e-3} for deciding which constraints are active.

We solve this problem with the classical gradient descent method, which does not consider the mesh quality constraints, and also use our novel approach proposed in this paper. Note that, for this problem the gradient descent method proved to be more robust than the BFGS method, which is why we only report the results obtained with the former method in the following. The evolution of the cost functional and the norm of the gradient deformation over the course of the optimization is shown in Figure~\ref{fig:history_packing}. Here, we observe that the classical approach failed after only 11 iterations, whereas our proposed method was able to perform the desired 90 iterations. The reason for the failure of the classical approach for this example is that the employed iterative linear solvers failed to converge due to the deterioration of the mesh quality, which will be discussed more detailedly below. Unfortunately, as discussed in Remark~\ref{rem:gmsh}, there was no possibility to perform a remeshing, so that the shape optimization failed for the classical approach. From the evolution of the cost functional, we can observe that the classical approach is a bit more efficient in increasing the objective functional during the first iterations compared to the approach with mesh quality constraints. However, as the classical approach failed after 11 iterations, our proposed method yields a significantly higher cost functional value after 90 iterations, so that the optimization with our approach is more successful. Regarding the norm of the gradient deformation we note that this is not reduced during the optimization for either approach. However, for the purposes of this engineering application, the optimization with mesh quality constraints is deemed successful as the cost functional significantly increased.

\begin{figure}[t]
	\centering
	\includegraphics[width=0.4\textwidth]{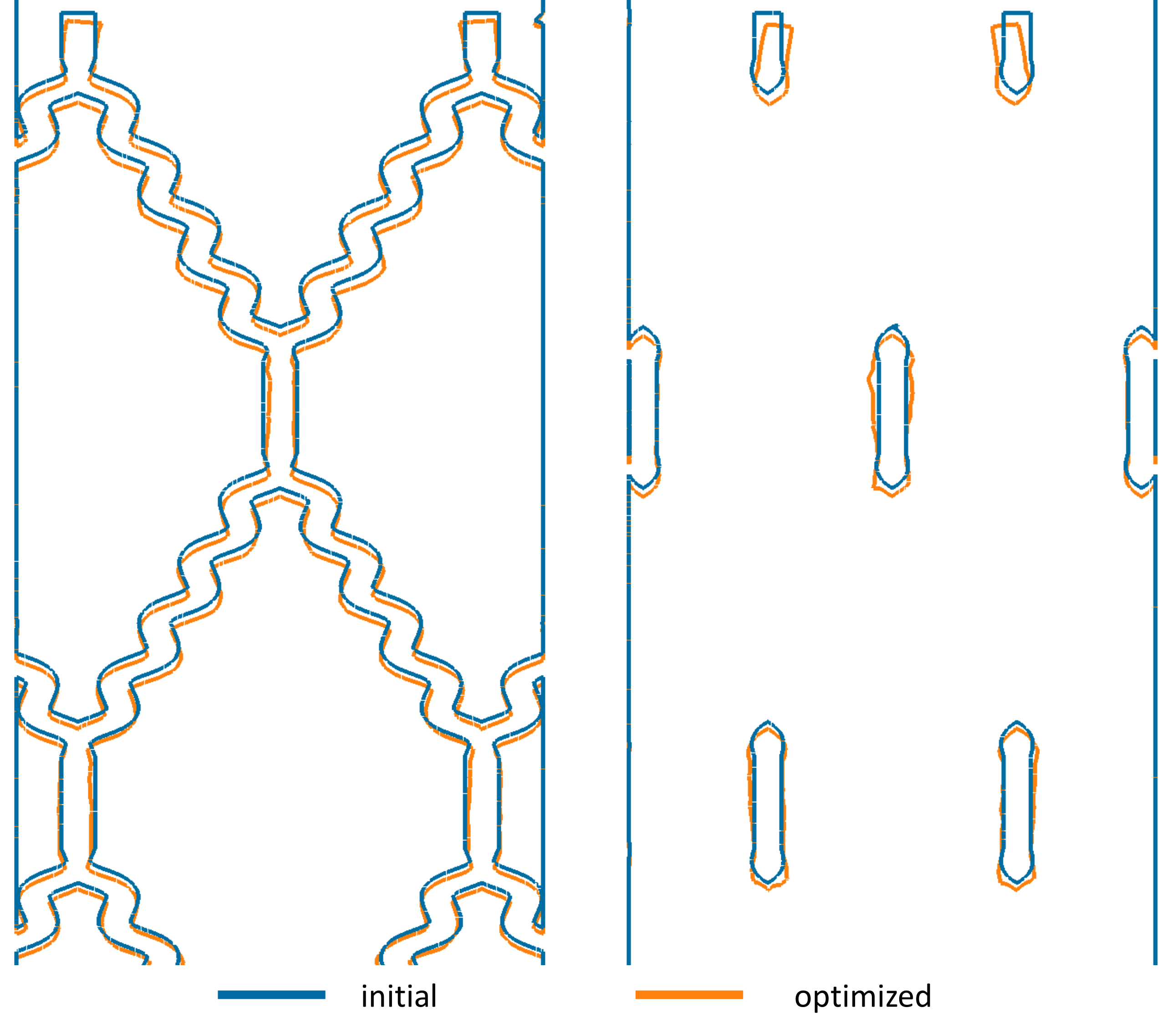}
	\caption{Comparison of the initial and optimized geometries on the $xy$-plane (left) and $yz$-plane (right) for problem \eqref{eq:shape_distillation}, taken from \cite{Blauth2025CFD}.}
	\label{fig:packing_geometry}
\end{figure}

\begin{figure}[b]
	\centering
	\begin{subfigure}[t]{0.475\textwidth}
		\centering
		\includegraphics[width=\textwidth]{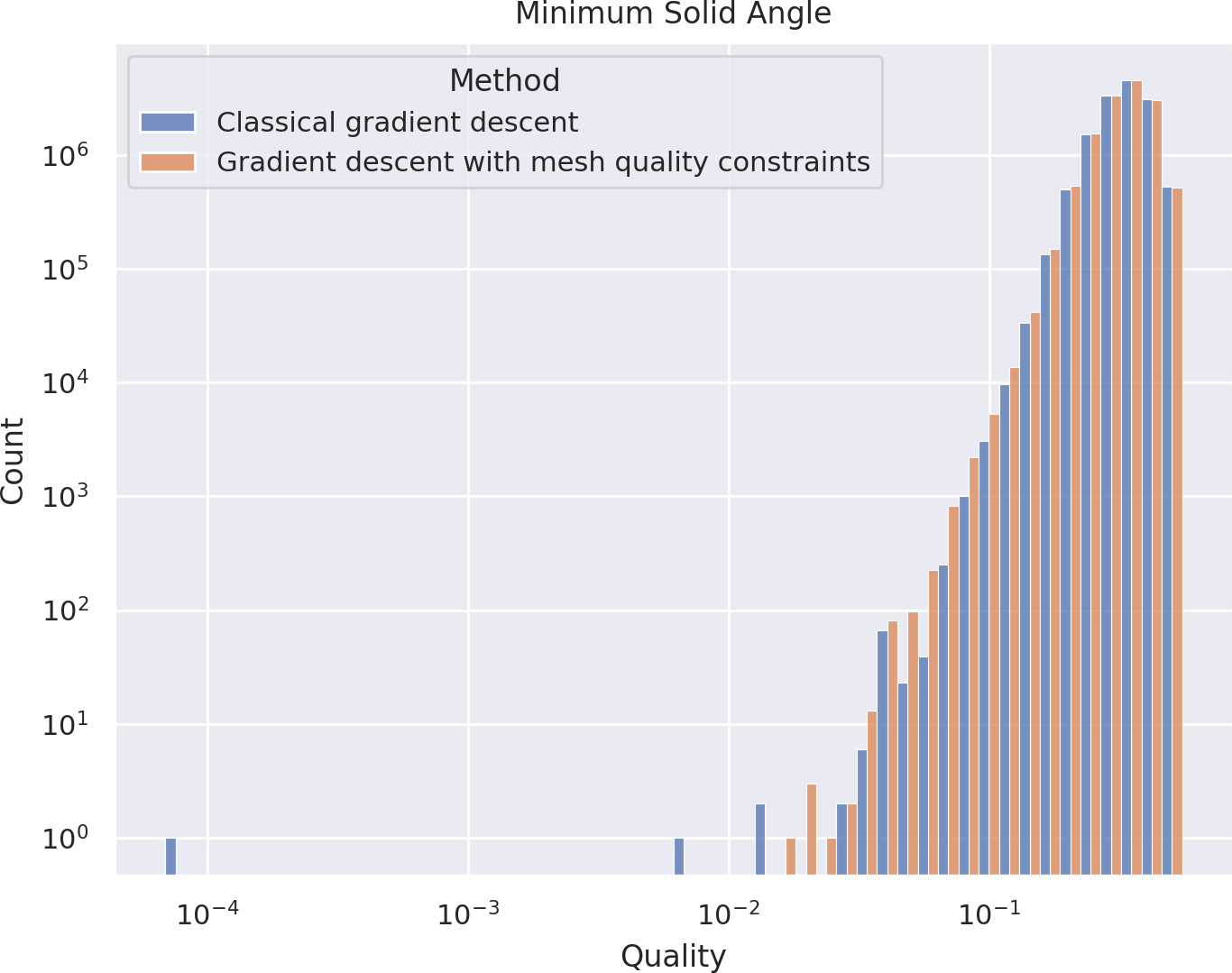}
		\caption{Minimum solid angle in steradian.}
		\label{fig:dist_packing_angle}
	\end{subfigure}
	\hfil
	\begin{subfigure}[t]{0.475\textwidth}
		\centering
		\includegraphics[width=\textwidth]{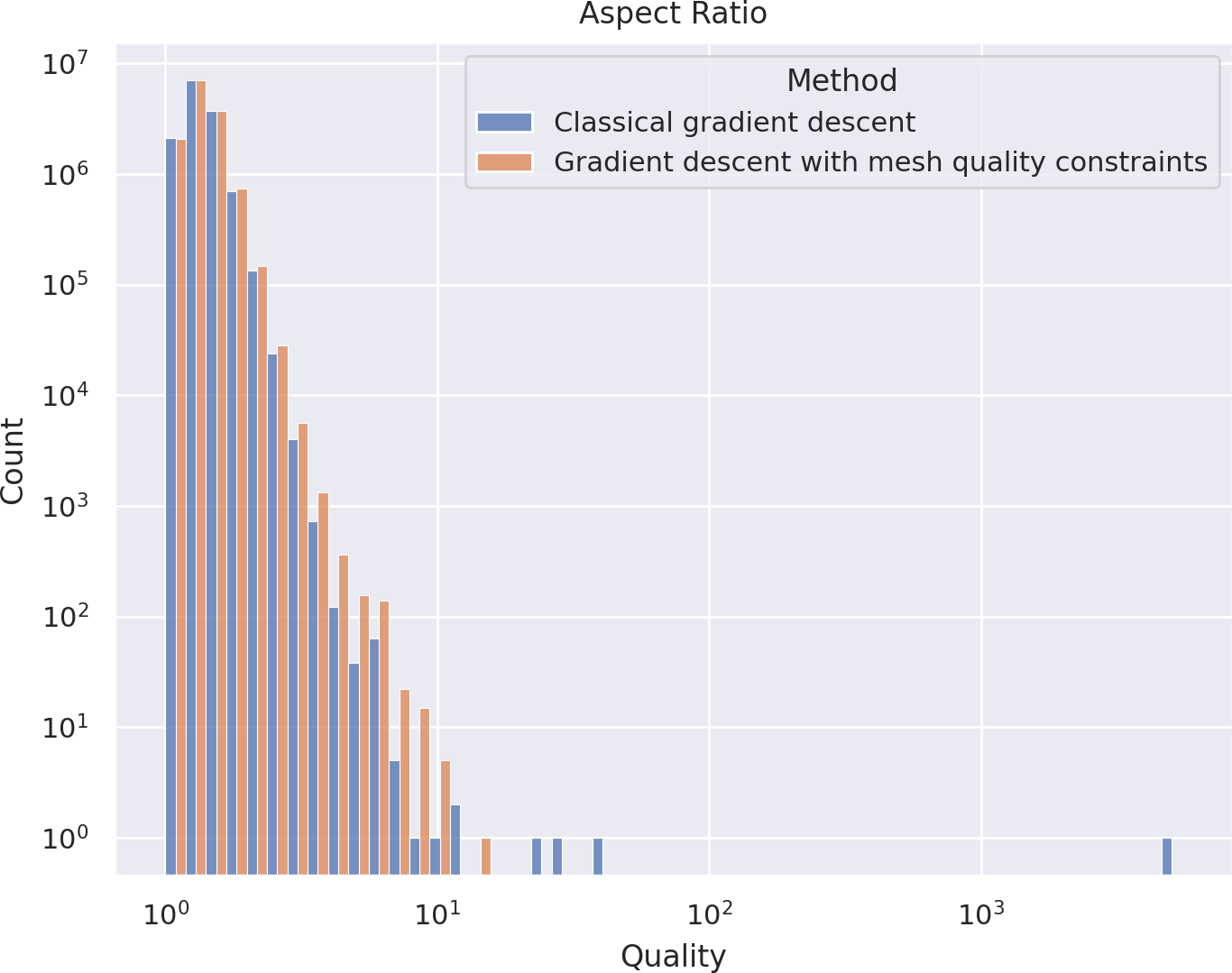}
		\caption{Aspect ratio.}
		\label{fig:dist_packing_aspect_ratio}
	\end{subfigure}
	\caption{Histograms showing the distribution of the mesh quality after the shape optimization for problem \eqref{eq:shape_distillation}.}
	\label{fig:dist_packing}
\end{figure}

In our previous publication \cite{Blauth2025CFD} it is shown that this shape optimization approach is very successful in practice. Our simulation results predict an increase of about \qty{20}{\percent} for the mass transfer coefficient $\beta$. Note that in Figure~\ref{fig:packing_geometry} a comparison of the initial and optimized geometries is shown in the $xy$- and $yz$-planes. 
For the validation of the shape optimization results, both the initial and optimized design were additively manufactured and investigated experimentally in \cite{Blauth2025CFD}. The experimental results were in great agreement with the simulation results and showed that the optimized design has a significantly enhanced separation efficiency, which was increased by about \qty{20}{\percent}. For the sake of brevity, we refer the reader to our publication \cite{Blauth2025CFD}, where the shape optimization results, geometrical changes, and experimental validation are discussed in greater detail.

\begin{figure}[t]
	\centering
	\begin{subfigure}[t]{0.475\textwidth}
		\centering
		\includegraphics[width=\textwidth]{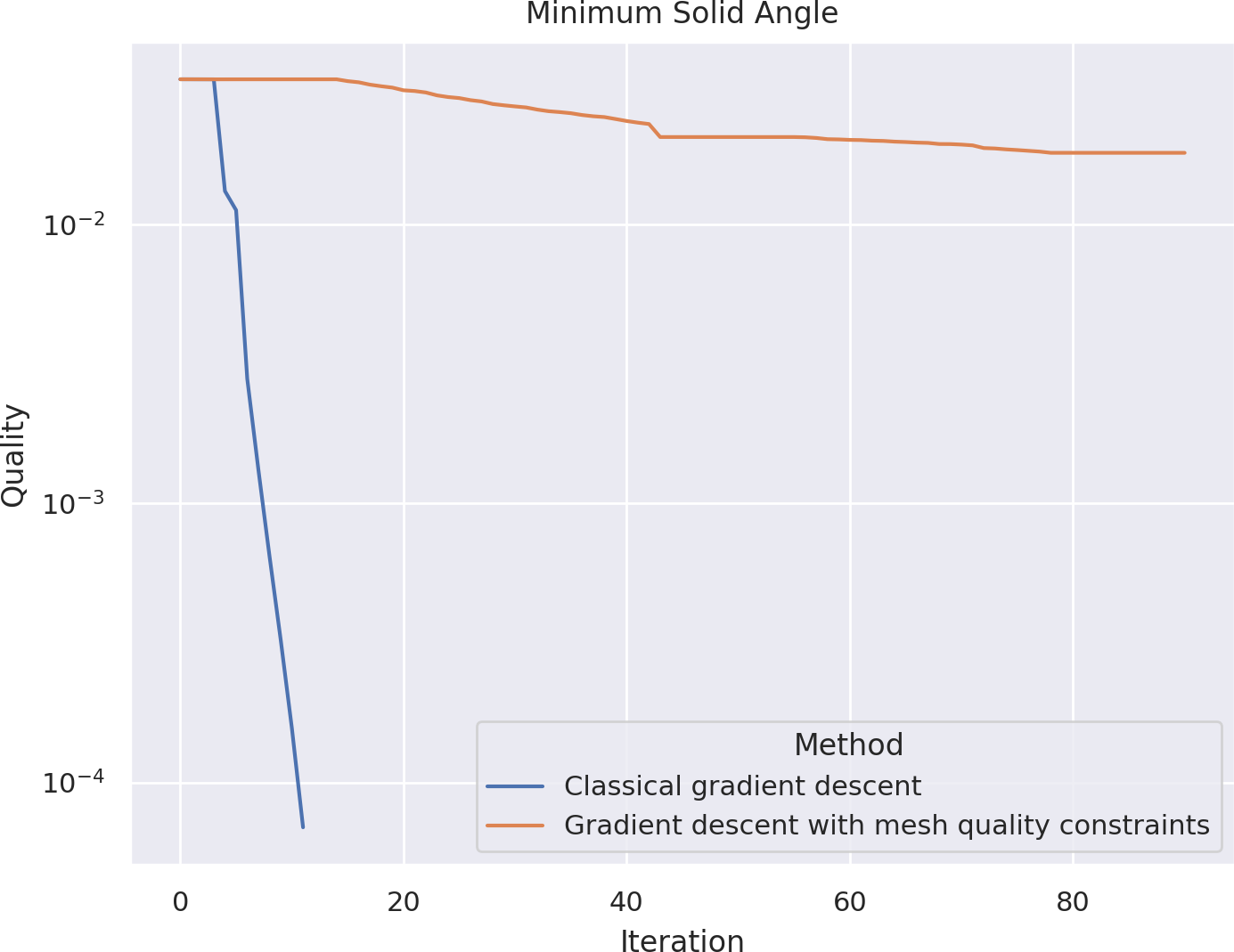}
		\caption{Minimum solid angle in steradian.}
	\end{subfigure}
	\hfil
	\begin{subfigure}[t]{0.475\textwidth}
		\centering
		\includegraphics[width=\textwidth]{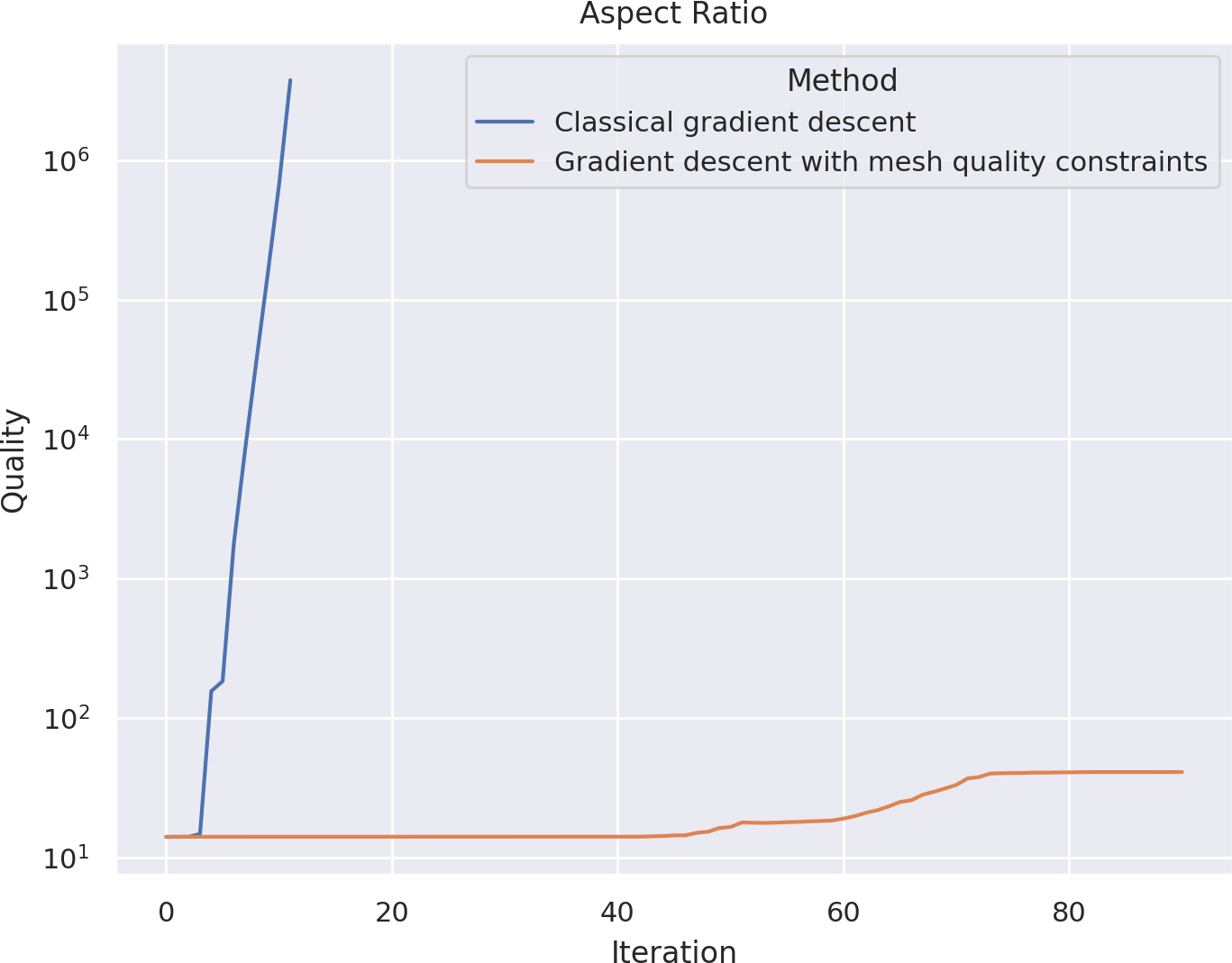}
		\caption{Aspect ratio.}
	\end{subfigure}%
	\caption{Evolution of the mesh quality for problem \eqref{eq:shape_distillation}.}
	\label{fig:iterations_packing}
\end{figure}

As stated previously, our proposed approach of using mesh quality constraints for shape optimization is detrimental for the successful shape optimization of this problem as the classical approach could not change the geometry significantly and failed after only 11 iterations. The reason for the failure of the classical gradient descent method is that the mesh quality degenerated due to the shape optimization which lead to failures of the iterative linear solvers.

In Figure~\ref{fig:dist_packing} histograms of the mesh quality are shown for the minimum solid angle and the aspect ratio of the mesh cells. Here, the aspect ratio of a tetrahedron is defined as $\frac{\abs{K}_\infty}{2 \sqrt{6} r}$, where $\abs{K}_\infty$ is the longest edge of the tetrahedron and $r$ is its inradius. We note that the $x$-axes of these histograms also have a logarithmic scaling as the quality of the worst cells is several magnitudes worse compared to the quality of the vast majority of the cells. It can be seen that the classical approach generates a mesh which has one particularly bad cell whose quality is deteriorating significantly and three other cells which also have a very bad quality. These few mesh cells are responsible for the breakdown of the classical gradient descent method after only 11 iterations. On the other hand, we observe that the quality of the worst mesh cells is substantially better with our proposed approach. We note that while there are some cells with a bad mesh quality for our approach, too, these cells already had a low quality in the initial mesh, as discussed above. Our approach ensures that the quality of these cells does not deteriorate much further and, for this, reason, enables the shape optimization of this complex example in the first place.

\begin{table}[t]
	\centering
	{\footnotesize
		\caption{Comparison of mesh quality of the optimized mesh for problem \eqref{eq:shape_distillation}.}
		\label{tab:shape_packing}
		\setlength{\tabcolsep}{1em}
		\begin{tabular}{l S S}
			\toprule
			\mbox{} & {Classical gradient descent} & {Gradient descent with mesh quality constraints} \\
			\midrule
			Minimum solid angle [\unit{\steradian}] & 6.904e-5 & 1.803e-02 \\
			Maximum aspect ratio & 5487.516 & 14.618 \\
			%
			%
			Active mesh quality constraints & {-} & 1584 \\
			Equality constraints for fixed boundaries & {-} & 195512 \\
			\bottomrule
		\end{tabular}
	}
\end{table}

These findings are confirmed by the results shown in Figure~\ref{fig:iterations_packing}, where the evolution of the minimum solid angle and the aspect ratio is depicted over the course of the optimization. As before, the mesh quality is defined to be the quality of the worst mesh cell. Again, we observe that our approach works very well and establishes a bound on the mesh quality during the entire optimization for both considered quality measures. In contrast, one can clearly observe that the mesh quality deteriorates very quickly for the classical approach and is about three orders of magnitude worse than the mesh quality obtained with our new approach. Finally, the mesh quality of the worst cells are also summarized in Table~\ref{tab:shape_packing}. Here, we can also see that only \num{1584} of the \num{54.8}~million mesh quality constraints are active during the final iteration of our proposed method, so that these only affect a tiny fraction of all mesh cells. Additionally, \num{195512} equality constraints for fixed boundaries were considered throughout the entire optimization.



\section{Conclusion and Outlook}
\label{sec:conclusion}

In this paper, we have presented a novel method for incorporating mesh quality constraints into shape optimization problems based on the gradient projection method. We recalled the gradient projection method as well as basic results from shape calculus and presented a general gradient-based shape optimization algorithm. To introduce our new approach, we investigated the discretization of the shape optimization algorithm with the finite element method. Afterwards we formulated constraints on the angles of triangular and solid angles of tetrahedral mesh cells, respectively, which bound these from below and above. This, consequently, bounds the quality of all mesh cells during the shape optimization. We discussed our implementation of the method based on our software cashocs \cite{Blauth2021cashocs}, focusing on additional constraints required for the case that some boundaries are fixed, the choice of the minimum angle threshold for the constraints, and the numerical solution of the additional equations for projecting the search directions. We investigated our proposed method for the drag minimization of an obstacle in Stokes flow, the optimization of the flow in a pipe governed by the Navier-Stokes equations, and the optimization of structured packings for distillation from \cite{Blauth2025CFD}. The numerical results highlight the great performance of our method as it leads to a significantly higher mesh quality than classical approaches with only a slight to moderate increase of computational cost. Moreover, our method enables the shape optimization of challenging large-scale problems with complicated geometries, which was not feasible previously, and preserves the mesh quality even in such settings.

For future work, a thorough numerical investigation and comparison of the method proposed in this paper with other approaches for preserving the mesh quality in shape optimization is of great interest. For this, we note that our method is, in principle, compatible with any approach that modifies the computation of search directions for shape optimization such as the ones presented in, e.g., \cite{Etling2020First, Mueller2021novel, Iglesias2018Two} and we expect that our approach will, with moderate additional numerical cost, enhance the performance of all of these approaches. Moreover, an extension of the method to other cell types, such as quadrilateral meshes for two-dimensional and hexahedral meshes for three-dimensional problems is of practical importance. Finally, the framework presented in this paper can also be extended to treat additional constraints on the mesh. One example for this would be to also consider the discrete curvature of the mesh and to bound this in order to ensure, e.g., the manufacturability of the optimized design.


\bibliographystyle{siamplain}
\bibliography{/p/tv/blauths/literature_db.bib}

\end{document}